\documentclass[leqno,10pt]{amsart}
\usepackage{amsmath,mathrsfs}
\usepackage{amssymb, amsmath, color}
\usepackage{graphicx}
\usepackage{dutchcal}
\usepackage{pdfsync}
\usepackage[colorlinks,
linkcolor=red,
anchorcolor=blue,
citecolor=green
]{hyperref}
\usepackage{colonequals}

\def\inte#1{
	\displaystyle\mathop{#1\kern0pt}^\circ }

\let\pa=\partial

\let\f=\frac

\def\pa{\partial}

\def\virgp{\raise 2pt\hbox{,}}
\def\cdotpv{\raise 2pt\hbox{;}}

\def\C{\mathop{\mathbb C\kern 0pt}\nolimits}
\def\DD{\mathop{\mathbb D\kern 0pt}\nolimits}
\def\EE{\mathop{{\mathbb E \kern 0pt}}\nolimits}
\def\K{\mathop{\mathbb K\kern 0pt}\nolimits}
\def\N{\mathop{\mathbb N\kern 0pt}\nolimits}
\def\Q{\mathop{\mathbb Q\kern 0pt}\nolimits}
\def\R{\mathop{\mathbb R\kern 0pt}\nolimits}
\def\SS{\mathop{\mathbb S\kern 0pt}\nolimits}
\def\ZZ{\mathop{\mathbb Z\kern 0pt}\nolimits}
\def\TT{\mathop{\mathbb T\kern 0pt}\nolimits}
\def\P{\mathop{\mathbb P\kern 0pt}\nolimits}

\newcommand{\lr}[1]{\langle #1 \rangle}

\def\na{\nabla}

\def\eps{\epsilon}

\def\hphi{\hat{\phi}}

\newcommand{\beq}{\begin{equation}}
	\newcommand{\eeq}{\end{equation}}
\newcommand{\ben}{\begin{eqnarray}}
	\newcommand{\een}{\end{eqnarray}}
\newcommand{\beno}{\begin{eqnarray*}}
	\newcommand{\eeno}{\end{eqnarray*}}

\DeclareMathOperator*{\esssup}{ess\,sup}

\newtheorem{defi}{Definition}[section]
\newtheorem{thm}{Theorem}[section]
\newtheorem{lem}{Lemma}[section]
\newtheorem{rmk}{Remark}[section]

\newtheorem{prop}{Proposition}[section]

\begin{document}
	
	\title[On semi-classical limit of homogeneous quantum Boltzmann equation]
	{On semi-classical limit of spatially homogeneous quantum Boltzmann equation: asymptotic expansion
	}
	
	\author[L.-B. He, X. Lu, M. Pulvirenti and Y.-L. Zhou]{Ling-Bing He, Xuguang Lu, Mario Pulvirenti and Yu-Long Zhou}
	\address[L.-B. He]{Department of Mathematical Sciences, Tsinghua University\\
		Beijing 100084,  P. R.  China.} \email{hlb@mail.tsinghua.edu.cn}
	\address[X. Lu]{Department of Mathematical Sciences, Tsinghua University\\
		Beijing 100084,  P. R.  China.} \email{xglu@mail.tsinghua.edu.cn}
	\address[M. Pulvirenti]{Dipartimento di Matematica, Universite di Roma La Sapienza, Piazzale Aldo Moro 5,
		00185 Rome, Italy;
		International Research Center M $\&$ MOCS, Universite dell'Aquila, Palazzo Caetani, Cisterna di Latina,
		(LT) 04012, Italy}  \email{pulviren@mat.uniroma1.it}
	\address[Y.-L. Zhou]{School of Mathematics, Sun Yat-Sen University, Guangzhou, 510275, P. R.  China.} \email{zhouyulong@mail.sysu.edu.cn}

	\begin{abstract}  We continue our previous work \cite{HLP1} on the limit of the spatially homogeneous quantum Boltzmann equation as the Planck constant $\eps$ tends to zero, also known as the semi-classical limit. For general interaction potential, we prove the following:
		(i). The spatially homogeneous quantum Boltzmann equations are locally well-posed in some weighted Sobolev spaces with  quantitative estimates uniformly in $\eps$.
		(ii). The semi-classical limit can be further  described by 
		the following asymptotic expansion formula:
		$$ f^\eps(t,v)=f_L(t,v)+O(\eps^{\vartheta}).$$
		This holds locally in time in Sobolev spaces. Here
		$f^\eps$ and $f_L$ are solutions to the quantum Boltzmann equation and the Fokker-Planck-Landau equation with the same initial data.  The convergent rate $0<\vartheta \leq 1$ depends on the integrability   of the Fourier transform of the particle interaction potential. Our new ingredients lie in a detailed analysis of the Uehling-Uhlenbeck operator from both angular cutoff and non-cutoff perspectives.
	\end{abstract}
	
	\maketitle
	
	\setcounter{tocdepth}{1}
	\tableofcontents
	
	\setcounter{equation}{0}
	
	\section{Introduction}
	The quantum Boltzmann equations for  Fermi-Dirac and
	Bose-Einstein statistics   proposed
	by Uehling and Uhlenbeck in \cite {UU} (after Nordheim \cite {N}) should be derived from the evolution of real Fermions   and Bosons  in the so called weak-coupling limit (see \cite{Bal} and \cite {BCEP}). Since Fokker-Planck-Landau equation is the effective equation associated with a dense and weakly interacting gas of classical particles (see \cite {ICM}, \cite {BPS}), it is not surprising that the semi-classical limits of the solutions to quantum Boltzmann equations are expected to be solutions to the Fokker-Planck-Landau equation.
	
	The weak convergence of the limit is justified in the paper \cite{HLP1}. The main purpose of this article is to provide a detailed asymptotic expansion formula to describe the semi-classical limit in the classical sense.

	\subsection{Setting of the problem}  The Cauchy problem of the spatially homogeneous quantum Boltzmann equation reads
	\ben\label{BUU} \pa_t f
	=Q_{UU}^{\eps}(f), \quad f|_{t=0} =f_0, \een 
	which describes time evolution of the gas given initial datum $f_0$.
	Here the solution
	$f=f(t, v)\ge 0$ is the density of the gas. The  Uehling-Uhlenbeck operator 
	$Q_{UU}^{\eps}$ in the weakly coupled regime is defined by 
	\ben\quad Q_{UU}^{\eps}(f) = \int_{\R^3\times\SS^2}  B^{\eps}(|v-v_*|,\cos\theta)
	\big(f_*'f'(1\pm\eps^3f_*) (1\pm\eps^3f) -f_*f(1\pm\eps^3f_*')(1\pm\eps^3f')\big)\mathrm{d}\sigma \mathrm{d}v_*,\label{OPQUU} 
	\een
	where 
	\ben \label{DefB}  
	B^{\eps}(|v-v_*|,\cos\theta) \colonequals \eps^{-4} |v-v_*| \left(\hat{\phi}\left(\eps^{-1}|v-v_*|\sin (\theta/2)\right) \pm \hat{\phi}\left(\eps^{-1}|v-v_*|\cos (\theta/2)\right) \right)^2.
	\een
	
	$\bullet$ {\it Some explanation on the model.}  On the derivation of \eqref{OPQUU} and \eqref{DefB} in the weak-coupling limit, we refer to \cite{BCEP,BCEP04,BCEP08,ESY04}. To make	\eqref{OPQUU} and \eqref{DefB} clear, we have the following remarks:
	\begin{enumerate}
		\item  The parameter $\eps$ is the Planck constant $\hbar$.  Note that for simplicity, we already drop the factor $2\pi$ appeared in \cite{HLP1}.  As our goal is to study the semi-classical limit, we always assume $0<\eps<1$.
		
		\item
		The sign $"+"$  and the sign $"-"$   correspond to Bose-Einstein particles and Fermi-Dirac particles respectively.
		
		\item The real-valued function
		$\hphi$ is the Fourier transform of the  particle interaction potential $\phi(|x|)$.

		\item The deviation angle $\theta$ is defined through $\cos\theta\colonequals\f{v-v_*}{|v-v_*|}\cdot \sigma$. Thanks to the symmetric property of the collision kernel, we can assume that $\theta\in [0,\pi/2]$.
		
		\item We   use the standard shorthand $h=h(v)$, $g_*=g(v_*)$,
		$h'=h(v')$, $g'_*=g(v'_*)$ where $v'$, $v_*'$ are given by
		\beno
		v'=\frac{v+v_{*}}{2}+\frac{|v-v_{*}|}{2}\sigma, \quad 
		v'_{*}=\frac{v+v_{*}}{2}-\frac{|v-v_{*}|}{2}\sigma, \quad  \sigma\in\SS^{2}.
		\eeno
	\end{enumerate}
	
	\smallskip
	
	$\bullet$ {\it Basic assumptions on the potential function $\phi$.}  
	For $a \geq 0$, let
	\ben \label{notation-I-Iprime}
	I_{a} \colonequals \int_0^\infty \hphi^2(r) r^a \mathrm{d}r, \quad 
	I^{\prime}_{a} \colonequals \int_0^\infty |r (\hphi)^{\prime}(r)|^2 r^{a} \mathrm{d}r.
	\een 
	Our basic assumptions on $\hphi$ are 
	\beno &&{\bf (A1)} \quad I_{0} + I_{3} + I^{\prime}_{3} < \infty. \quad \quad \quad \quad \quad \quad \quad \quad \quad \quad \quad \quad \quad \quad \quad \quad \quad \quad \quad \quad \quad \quad \quad \quad \quad \quad \quad \quad
	\\ &&{\bf(A2)} \quad I_{3+\vartheta} + I^{\prime}_{3+\vartheta} <\infty \text{ for some }\vartheta\in(0,1]. \quad \quad \quad \quad \quad \quad \quad \quad \quad \quad \quad \quad \quad \quad \quad \quad \quad \quad \quad \quad \quad \quad \quad \quad \quad \quad \quad \quad
	\eeno
	Several remarks are in order:
	\begin{enumerate} \item The condition $I_{0}<\infty$ in  {\bf(A1)} is used to bound $\int_{\SS^2} B^\epsilon (|v-v_*|,\cos\theta)\mathrm{d}\sigma<\eps^{-3}I_{0}$, 
		see \eqref{upper-bound-of-A-eps} for details. This is the key point to prove the {\it global existence}  of the mild solution for the Fermi-Dirac particles. 
		 In the weak coupling regime \eqref{DefB},  finiteness of the $\sigma$-integral holds even for some inverse power law potentials. Indeed, taking $\phi(|x|)=|x|^{-p} (0<p<3)$,
		 one can check that  $\int_{\SS^2} B^\epsilon (|v-v_*|,\cos\theta)\mathrm{d}\sigma$ is finite
		 for $p>2$ and infinite for $p \leq 2$, see \cite{Z-2023} for more details on this. For the infinite case, one may need some angular cutoff.
		 The condition $I_{0}<\infty$ is reminiscent of Grad's angular cutoff assumption for inverse power law potentials in the low density regime,
		 which allows one separate the Boltzmann operator into gain and loss terms. 
From now on, we will call such mathematical treatment as "angular cutoff" view. 
However, such treatment is not enough since the upper bound blows up as $\eps \to 0$.

		\item   $I_{3}$ is derived by computing   the momentum transfer which is defined as follows:
		\ben\label{MometumTrans} M_o^\epsilon (|v-v_*|) \colonequals
		\int_{\SS^2} B^\epsilon (|v-v_*|,\cos\theta)(1-\cos\theta)\mathrm{d}\sigma.
		\een
		We will show in the later that $M_o^\epsilon\sim   I_{3}|v-v_*|^{-3}$ when $\epsilon$ is sufficiently small. This
		is also related to the diffusion coefficient of Fokker-Plank-Landau collision operator(see \eqref{OPQL} and \eqref{def-a}). 
		The condition $I_{3}<\infty$
		is reminiscent of angular non-cutoff kernels for inverse power law potentials as one always relies an additional order-2 $\theta^2$ to deal with angular singularity.
		From now on, we will call such mathematical treatment as "angular non-cutoff" view.

		\item The condition $I_{3} + I^{\prime}_{3} < \infty$ in
		{\bf(A1)} allows us to derive the cancellation lemma from the point view of angular non-cutoff. It plays the essential role to get the  uniform-in-$\eps$ estimate.
		
		\item Assumption {\bf(A1)}  is used to prove uniform-in-$\eps$  local well-posedness and propagation of regularity. 
		To get the asymptotic expansion for the semi-classical limit, technically we  need assumption
		{\bf(A2)}.
		
		\item We do not impose any point-wise condition on $\hat{\phi}$. All the constants depending on $\phi$ in this article are through the quantities in {\bf(A1)} and  {\bf(A2)}.
	\end{enumerate}

	\medskip
	
	By formal computation(see \cite{BP}), the Cauchy problem \eqref{BUU} will converge to that of the Fokker-Planck Landau equation 
	\ben\label{landau}
	\partial_t f = Q_L(f,f), \quad f|_{t=0} =f_0.
	\een 
	Here the Fokker-Planck Landau operator $Q_L$ reads
	\ben \label{OPQL}
	Q_L(g,h)(v) \colonequals \nabla \cdot \int_{\R^3} a(v-v_*) \{ g(v_*) \nabla h(v) - \nabla g(v_*) h(v)     \} \, \mathrm{d}v_*,
	\een
	where $a$ is a matrix-valued function that is symmetric, (semi-definite) positive.
	It depends on the interaction potential between particles, and is defined by (for $i,j =1, 2, 3$)
	\ben\label{def-a}
	a_{ij}(z) =2\pi I_{3}|z|^{-1}\, \Pi_{ij} (z), \quad \Pi_{ij}(z)= \delta_{ij} - \frac{z_i z_j}{|z|^2},
	\een
	where $I_{3}$ is defined in \eqref{notation-I-Iprime}.
	\smallskip
	
	Our goal is to study the semi-classical limit from 
	\eqref{BUU} to \eqref{landau} in some weighted Sobolev space. To do that, we seperate our proof into two parts: well-posedness results for \eqref{BUU} with uniform-in-$\epsilon$ estimates  and the asymptotic expansion formula with explicit error estimates.
	
	\subsection{Main results}  Before introducing the main results, we list some facts on the notations.

	\noindent$\bullet$  As usual, $a\lesssim b$ is used to denote that there is a
	universal constant $C$ such that $a \leq C b$.
	The notation $a\sim b$ means $a\lesssim b$ and $b\lesssim
	a$. We denote by $C(\lambda_1,\lambda_2,\cdots, \lambda_n)$ or $C_{\lambda_1,\lambda_2,\cdots, \lambda_n}$ some constant depending on parameters $\lambda_1,\lambda_2,\cdots, \lambda_n$. The notation $a\lesssim_{\lambda_1,\lambda_2,\cdots, \lambda_n} b$ is interpreted as $a \leq C_{\lambda_1,\lambda_2,\cdots, \lambda_n} b$.

	\noindent$\bullet$ We recall the $L^{p}$ space for $1 \leq p \leq \infty$ through the norm
	\beno 
	\|f\|_{L^p}	\colonequals  \left(\int_{\R^3}  |f(v)|^p \mathrm{d}v \right)^{1/p} \text{ for } 1 \leq p < \infty; \quad  \|f\|_{L^\infty}	\colonequals  \esssup_{v \in \R^3} |f(v)|.
	\eeno
	Denote the weight function by $W_{l}(v) \colonequals  (1+|v|^2)^{\frac{l}{2}}$ for $l \in \R$ and write $W = W_{1}$ for simplicity. Then the weighted $L^{p}_{l}$ space is defined through the norm $\|f\|_{L^p_l} \colonequals  \|W_l f\|_{L^p}$. We denote the multi-index $\alpha =(\alpha _1,\alpha _2,\alpha _3) \in \N^3$ with
	$|\alpha |=\alpha _1+\alpha _2+\alpha _3$.
	For up to order $N \in \N$ derivatives,
	the weighted Sobolev space $W^{N,p}_l$ on $\R^3$ with $p\in [1, \infty], l \in \R$ is defined  through the following norm	\beno
	\|f\|_{W^{N,p}_l}
	\colonequals
	\sum_{|\alpha|\le N}  \|\partial^\alpha f\|_{L^p_l}.
	\eeno
	If $p = 2$, denote by $H^{N}_l$ the Hilbert space  with
	$\|f\|_{H^{N}_l} = \|f\|_{W^{N,2}_l}$. 
	
	\noindent$\bullet$ For simplicity, for $T>0$, let $\mathcal{A}_{T}  \colonequals L^{\infty}([0,T]; L^{1}(\R^3))$ associated with the norm $\|f\|_{T} \colonequals \sup_{0 \leq t \leq T} \|f(t)\|_{L^{1}}$. Let $\mathcal{A}_{\infty}  \colonequals L^{\infty}([0,\infty); L^{1}(\R^3))$ associated with the norm $\|f\|_{\infty} \colonequals \sup_{t \geq 0} \|f(t)\|_{L^{1}}$.

	\noindent$\bullet$  Given a non-negative initial datum $f_0 \in L^{1}(\R^3) \cap L^{\infty}(\R^3)$. We consider the initial value problem \eqref{BUU} for Bose-Einstein particles with $0<\eps<1$ and for Fermi-Dirac particles with $0< \eps \leq \min \{ 1, \|f_0\|_{L^{\infty}}^{-1/3}\}$. Set $\|f_0\|_{L^{\infty}}^{-1/3} = \infty$ for the  case $\|f_0\|_{L^{\infty}} = 0$ where the problem trivially has a zero solution.
	
	\noindent$\bullet$ For simplicity, we will use the shorthand $\int (\cdots) \mathrm{d}V =\int_{v, v_* \in \R^3, \sigma \in \SS^2, (v - v_*) \cdot \sigma \geq 0} (\cdots) \mathrm{d}v \mathrm{d}v_* \mathrm{d}\sigma$. We drop integration domain in most of the integrals if there is no confusion.
	
	\medskip
	
	Next we introduce the definition of mild solution to \eqref{BUU}.
	
	\begin{defi} \label{def-mild-solution}	
		For $T>0$, set $\mathcal{A}_{T} \colonequals  L^{\infty}([0,T]; L^{1}(\R^3) \cap L^{\infty}(\R^3)$. 
		A measurable non-negative function $f \in \mathcal{A}_{T}$ on $[0, T] \times  \R^3$ is called a local(or global)
		mild solution of the initial value problem \eqref{BUU} if  $T < \infty$(or $T = \infty$)
		there is a null set $Z \subset \R^3$ s.t., for all $t \in  [0, T]$ and $v \in \R^3 \setminus Z $, 
		\beno 
		f(t,v) = f_{0}(v) + \int_{0}^{t} Q_{UU}^{\eps}(f)(\tau, v) \mathrm{d}\tau,
		\eeno
		and additionally for Fermi-Dirac particles, it holds that
		\ben \label{upper-bound}
		\|f(t)\|_{L^{\infty}} \leq \eps^{-3}.
		\een
	\end{defi}

	Our first result is the global well-posedness and propagation of regularity of the Cauchy problem \eqref{BUU} for Fermi-Dirac particles.
	
	\begin{thm}[Fermi-Dirac particles]\label{gwp-F-D} Let $\hat{\phi}$ verify  \textbf{(A1)} and $0\le f_0 \in L^1 \cap L^{\infty}$. Suppose that $0< \eps \leq \min \{1, \|f_0\|_{L^{\infty}}^{-1/3}\}$.
		\begin{enumerate}
			\item{\bf (Global well-posedness)}	The Cauchy problem \eqref{BUU} for Fermi-Dirac particles admits a unique global mild solution $f^\eps$.
			
			\item{\bf (Propagation of regularity uniformly in $\eps$)}
			If $f_0 \in H^{N}_{l}$ for $N, l \geq 2$, 
			there exists a lifespan $T^*=T^*(N, l, \phi, \|f_0\|_{H^{N}_{l}})>0$ independent of $\eps$,
			such that the family of solution $\{f^{\eps}\}_{\eps}$ is uniformly bounded in
			$L^{\infty}([0, T^*]; H^{N}_{l}) \cap C([0, T^*]; H^{N-2}_{l})$.  More precisely, uniformly in $\eps$, 
			\ben \label{bounded-by-initial-data}
			\sup_{t\in[0,T^*]} \|f^{\eps}(t)\|_{H^{N}_{l}} \le 2\|f_0\|_{H^{N}_{l}}, \quad 
			\een
			and for $0 \leq t_1 \leq t_2 \leq T^*$,
			\ben \label{continuous-in-time}
			\|f^{\eps}(t_2) - f^{\eps}(t_1)\|_{H^{N-2}_{l}} \le C(N, l, \phi, \|f_0\|_{H^{N}_{l}}) (\|f_0\|_{H^{N}_{l}}^2 + \|f_0\|_{H^{N}_{l}}^3) (t_2 - t_1). 
			\een
		\end{enumerate}	 
	\end{thm}

	Our second result is the local well-posedness and propagation of regularity $H^{N}_{l}$ of the Cauchy problem \eqref{BUU} for Bose-Einstein particles.
	\begin{thm}[Bose-Einstein particles]\label{lwp-B-E} Let $\hat{\phi}$ verify  \textbf{(A1)}. Let $0 \leq f_0 \in H^{N}_{l}$ for $N, l \geq 2$, then for any $0<\eps<1$, the Cauchy problem \eqref{BUU} for  Bose-Einstein particles admits a unique local  mild solution $f^\eps \in L^{\infty}([0, T^*]; H^{N}_{l}) \cap C([0, T^*]; H^{N-2}_{l})$ where $T^*=T^*(N, l, \phi, \|f_0\|_{H^{N}_{l}})>0$ is independent of $0< \eps <1$. Moreover, the family of solution $\{f^{\eps}\}_{0< \eps < 1}$ satisfies \eqref{bounded-by-initial-data} and \eqref{continuous-in-time} uniformly in $\eps$.
	\end{thm}

	\begin{rmk}
		Note that \eqref{continuous-in-time} implies \eqref{BUU} holds in the space $H^{N-2}_{l}$ for almost all $t \in [0, T]$.
		Thanks to the weak convergence result in \cite{HLP1}, 
		similar local well-posedness also holds for the Landau equation \eqref{landau} with estimates \eqref{bounded-by-initial-data} and \eqref{continuous-in-time}.  
	\end{rmk}
	
	Our last result is on the asymptotic expansion for the semi-classical limit.

	\begin{thm} [Semi-classical limit with convergence rate]
		\label{mainthm}  Let $0 \leq N \in \N, 2 \leq l \in \R$. 
		Suppose that 		
		(i). $\hat{\phi}$ satisfies \textbf{(A1)} and  \textbf{(A2)};
		(ii).
		$0\le f_0 \in H^{N+3}_{l+5}$;
		(iii).
		For Fermi-Dirac particles, 
		$0< \eps \leq \min \{1, \|f_0\|_{L^{\infty}}^{-1/3}\}$.
		Let $f^\eps$ and  $f_L$ be the solutions to \eqref{BUU} and \eqref{landau} respectively with the initial datum $f_0$ 
		on $[0, T^*]$ where $T^*=T^*(N, l, \phi, \|f_0\|_{H^{N+3}_{l+5}})$ given in Theorem \ref{gwp-F-D} and \ref{lwp-B-E}. Then  for $t\in [0,T^*]$, it holds that
		\ben\label{AsEx}f^\eps(t,v)=f_L(t,v)+\eps^{\vartheta} R^\eps(t,v), \een where
		\ben \label{error-estimate}
		\sup_{t\in[0,T^*]}  \|R^\eps\|_{H^{N}_{l}}
		\leq   C(\|f_0\|_{H^{N+3}_{l+5}}; N, l, \phi).  \een
	\end{thm}
	
	Some comments on these results are in order:
	
	\begin{rmk} \label{eps-dependent-propagation} Owing to the fact that Fermi-Dirac particles enjoy the $L^\infty$ upper bound \eqref{upper-bound}, indeed, we can prove the global propagation of regularity with the quantitative estimates as follows:
		
		\begin{itemize}\item{\bf (Global propagation of regularity)}	If 
			$f_0 \in L^2_{l}$ for $l \geq 2$, then
			for any $t \geq 0$,
			\ben \label{F-D-propgation-L2-l}
			\| f^\eps(t)\|_{L^{2}_l} \leq \| f_0\|_{L^{2}_l} \exp \left(t C_{\eps,\phi,l}(\|f_0\|_{L^1} + \eps^{-3}) \right).	
			\een 
			If 
			$f_0 \in L^1_l \cap L^{\infty}_l \cap H^1_l$ for $l \geq 2$, then
			for any $t \geq 0$,
			\ben \label{F-D-propgation-H1-l}
			\| f^\eps(t)\|_{L^1_l \cap L^{\infty}_l \cap H^1_l} \leq C(\| f_0\|_{L^1_l \cap L^{\infty}_l \cap H^1_l}, t; \eps, \phi, l).	
			\een
			If 
			$f_0 \in W^{1,1}_l \cap W^{1,\infty}_l \cap H^{N}_l$ for $N, l \geq 2$, then
			for any $t \geq 0$,
			\ben \label{F-D-propgation-HN-l}
			\| f^\eps(t)\|_{W^{1,1}_l \cap W^{1,\infty}_l \cap H^{N}_l} \leq C(\| f_0\|_{W^{1,1}_l \cap W^{1,\infty}_l \cap H^{N}_l}, t; \eps, \phi, l, N).
			\een \end{itemize}
		We cannot expect similar results for the B-E particles because of the  B-E condensation phenomenon.
		The proofs of these propagation results are based on the fact that derivatives and weights can be suitably distributed  across 
		the non-linear terms such that the targeting high-order norm grows at most linearly under the premise that the lower-order norms are already propagated globally. Since these results are somewhat deviated from the main purpose of this article, their proofs are not given for brevity.
	\end{rmk}
	
	\begin{rmk} 
		The propagation of regularity uniformly in $\eps$   in equation \eqref{bounded-by-initial-data} and the temporal continuity described in equation \eqref{continuous-in-time} are applicable to both Fermi-Dirac and Bose-Einstein particles. This enables us to delve deeper into exploring the semi-classical limit as presented in Theorem \ref{mainthm}.
	\end{rmk}
	
	\begin{rmk}\label{reason-of-uniform-propagation} To get the upper bound of the error term $R^\eps$ uniformly in $\eps$,  
		it is compulsory to impose high regularity on the solutions $f^\eps$ and $f_L$ since we need to estimate the error between  $Q_{UU}^{\eps}(f^\eps)$ and $Q_L(f_L)$, see \eqref{upper-bound-3-order-higher} for the error equation and
		Lemma \ref{Q1-and-QL} for the main estimate.
	\end{rmk}
	
	\begin{rmk}  We emphasize that the asymptotic expansion \eqref{AsEx} is sharp. This can be easily seen from the proof of Theorem \ref{mainthm}. Roughly speaking, to get the factor $\eps^{\vartheta}$, we need to kill the singularity which behaves like the Riesz potential $|x|^{-2-\vartheta}$. Obviously this singularity can be removed in the case of $\vartheta <1$. For the borderline case  $\vartheta=1$, we further check that  the corresponding  part in fact behaves   as  $\f{K(x)}{|x|^3}$ which is the kernel of the typical  Calderon-Zygmund operator. From this point, the expansion $f^\eps=f_L+O(\eps)$ is sharp for any smooth potential function $\phi$.
	\end{rmk}
	
		\begin{rmk}  The asymptotic formula	\eqref{AsEx} holds locally in time for any initial data in  $H^{N+3}_{l+5}$. We may  expect large time or even global-in-time result for some special 
			initial data,
			considering the recent progress on global well-posedness in the homogeneous \cite{LiLu} and inhomogeneous  \cite{bae2021relativistic,JXZ,ouyang2021quantum,Z-2023} case.
	\end{rmk}

	\subsection{Short review} Quantum Boltzmann equation has been widely investigated by many authors. In this subsection, we first give a short review on the existing results. Then we explain
	the main difficulty for the problem of semi-classical limit. 
	
	In most of the literature on the quantum Boltzmann equation, the authors usually take  $\eps=1$ in the definition of Uehling-Uhlenbeck operator \eqref{OPQUU}. In this situation, for Fermi-Dirac particles, we  have the a priori bound for the solution $f$, that is, $f\le 1$. We refer readers to \cite{Lu2001,Lu2008} for the existence result. For Bose-Einstein particles, we refer readers to \cite{Lu2004,Briant-Einav} for the existence of measure solution and the local well-posedness in weighted $L^\infty$ spaces. For the Bose-Einstein condensation at low temperature, we refer to \cite{EV1,EV2} and also \cite{CaiLu,LiLu,Lu2016} for the recent progress.

	As for semi-classical limit of \eqref{BUU}, in \cite{BP}, Benedetto and Pulvirenti proved the convergence of the
	operator. More precisely, for a suitable class of integrable functions $f$ and any Schwartz function $\varphi$,
	\ben\label{asyOP} \lim\limits_{\eps\to 0}\langle Q_{UU}^{\eps}(f), \varphi \rangle=\langle Q_{L}(f,f), \varphi\rangle 
	.
	\een
	The notation  $\langle f, g\rangle \colonequals \int_{\R^3}f(v)g(v)\mathrm{d}v$ is used to denote the inner product for $v$ variable.
	In \cite{HLP1}, under some assumptions on $\hphi$,  the following results are proved: (1). Starting from the Eq.(Fermi-Dirac),  up to subsequences,   the {\it isotropic weak solution} to  Eq.(Fermi-Dirac) will converge to  the {\it isotropic weak solution} to Eq.(Fokker-Planck-Landau); (2). Starting from the Eq.(Bose-Einstein),  up to subsequences,    the {\it measure-valued isotropic weak solution } to  Eq.(Bose-Einstein) will converge to  the {\it measure-valued isotropic weak solution } to Eq.(Fokker-Planck-Landau). Here {\it isotropic  solution} means that the solution $f(t,v)$ is a radial function with respect to $v$, that is, $f(t,v)=f(t,|v|)$. To achieve these results,  the main idea is to reformulate the equations in the isotropic sense and then make full use of the cancellation   hidden in the cubic terms.
	\smallskip

	\subsection{Difficulties, strategies and new ideas} The main difficulty is induced by the singular scaling factor in the Uehling-Uhlenbeck operator \eqref{OPQUU}. One may attempt to use some normalization technique to deal with the parameter $\eps$. For instance,
	  if
	$\tilde{f}(t,v) \colonequals \eps^3 f(\eps^6t,\eps v)$,  
	one can easily verify that   $\tilde{f}$ is a solution to following equation:
	\ben\label{BUU1}
	\pa_t \tilde{f} = \int_{\R^3\times\SS^2} B^{1}(|v-v_*|,\cos\theta) \big(\tilde{f}_*'\tilde{f}'(1\pm \tilde{f}_*) (1\pm \tilde{f}) -\tilde{f}_*\tilde{f}(1\pm \tilde{f}_*')(1\pm \tilde{f}')\big)\mathrm{d}\sigma \mathrm{d}v_*.
	\een  
	Now \eqref{BUU} is reduced to  \eqref{BUU1} with the initial data $\tilde{f}|_{t=0}\colonequals \eps^3 f_0(\eps v)$. The good side is that
	 the equation \eqref{BUU1} itself contains no $\eps$.
	However, the bad side is that the initial data $\tilde{f}|_{t=0}$ is sufficiently large in weighted Sobolev spaces when $\eps$ is sufficiently small. It is challenging to establish 
	uniform lifespan for the nonlinear equations like \eqref{BUU1} starting from arbitrarily large initial data. Usually, 
	 the lifespan vanishes as the initial data blows up.  Therefore, we will directly consider \eqref{BUU}.
	\smallskip
	
Let us explain our strategy from the analysis of the collision operator. It is easy to see that   we can decompose $Q_{UU}^{\eps}$ into two parts:
	\ben \label{OPUU}  Q_{UU}^{\eps}(f) = Q(f,f)+R(f,f,f), \een
	where $Q(f,f)$ contains the quadratic terms and $R(f,f,f)$ contains the cubic terms. More precisely, 
	\ben\label{OPQ}
	Q(g,h)\colonequals\int_{\R^3\times\SS^2}  B^{\eps}(|v-v_*|, \cos\theta) (g_*'h'-g_*h)\mathrm{d}\sigma \mathrm{d}v_*  =  \sum_{i=1}^3 Q_i(g,h);
	\\\label{OPR}
	R(g,h,\rho)\colonequals  \pm\eps^3\int_{\R^3\times\SS^2} B^{\eps}(|v-v_*|, \cos\theta) ( g_*'h'(\rho+\rho_*)- g_*h(\rho'+\rho_*'))\mathrm{d}\sigma \mathrm{d}v_*.\een
	Here in \eqref{OPQ} for $i=1,2,3$, $Q_{i}$ is defined by
	\beno 
	Q_{i}(g,h) \colonequals \int_{\R^3\times\SS^2}   B^{\eps}_i(|v-v_*|, \cos\theta) (g_*'h'-g_*h)\mathrm{d}\sigma \mathrm{d}v_*, 
	\eeno
	where $B^{\eps}_{i}$ is defined by
	\ben 
	\label{DefB1}  B^{\eps}_{1}(|v-v_*|,\cos\theta) &\colonequals& \eps^{-4} |v-v_*|\hat{\phi}^2\left(\eps^{-1}|v-v_*|\sin (\theta/2)\right),
	\\
	\label{DefB2}
	B^{\eps}_{2}(|v-v_*|,\cos\theta) &\colonequals& \pm 2\eps^{-4} |v-v_*|\hat{\phi}\left(\eps^{-1}|v-v_*|\sin (\theta/2)\right)\hat{\phi} \left(\eps^{-1}|v-v_*|\cos(\theta/2) \right),
	\\
	\label{DefB3}  B^{\eps}_{3}(|v-v_*|,\cos\theta) &\colonequals& \eps^{-4} |v-v_*|\hat{\phi}^2\left(\eps^{-1}|v-v_*|\cos (\theta/2)\right).
	\een

	\begin{rmk} \label{split-kernel}
	Note that $B^{\eps} =  B^{\eps}_1+B^{\eps}_2+B^{\eps}_3$ and
	$|B^{\eps}_2| = 2 \sqrt{B^{\eps}_1 B^{\eps}_3}$.
	The divergence in $B^{\eps}_1$
	 arises from different physical reasons: the intensity of the collisions increases together with effective domain of $\hat{\phi}^2(\cdot)$
	 due to the vanishing of the scattering angle. 
	More precisely,  $\eps^{-1}|v-v_*|\sin (\theta/2) 
	\in [0, \eps^{-1}|v-v_*|\sqrt{2}/2]$ contains the dominate part of $\hat{\phi}^2(\cdot)$ as $\eps \to 0$. While for $B^{\eps}_3$,  $\eps^{-1}|v-v_*|\cos (\theta/2) 
	\in [\eps^{-1}|v-v_*|\sqrt{2}/2, \eps^{-1}|v-v_*|]$
	 goes to infinity and plays minor roles as $\eps \to 0$ because of the integrability conditions \textbf{(A1)} and  \textbf{(A2)}
	 on $\hat{\phi}^2$. In a word, in the limiting process  $\eps \to 0$, 
	 $B^{\eps}_1$ is the dominant part. 
	 If $B^{\eps}$ is replaced by $B^{\eps}_1$ and the cubic terms are ignored we would expect the same limiting behavior. See also Remark 1.1 in \cite{HLP1} for another treatment of the kernel and relevant discussions.
	\end{rmk}

	 In what follows, we will show that  $B^{\eps}_1$ and $B^{\eps}_3$ should be treated in a different manner. Roughly speaking, we will treat $B^{\eps}_1$ from the non-cutoff view and  $B^{\eps}_3$ from the cutoff view. This can be seen easily by the following computations.
	
	\underline{(i).} It is not difficult to derive that
	\[ \int_{\SS^2} B^{\eps}_{3}(|v-v_*|,\cos\theta) \mathrm{d}\sigma\sim I_3|v-v_*|^{-3},\]
	where we use the facts that $\cos(\theta/2)\sim1$ and the change of variables from $\cos(\theta/2)$ to $r\colonequals\eps^{-1}|v-v_*|\cos (\theta/2)$. The strong singularity induced by the relative velocity $|v-v_*|^{-3}$ is consistence with the Landau collision operator \eqref{OPQL}. The singularity can be removed by sacrificing the slight regularity of the solution.
	
	\underline{(ii).} Again by the similar calculation, we have
	\[ \int_{\SS^2} B^{\eps}_{1}(|v-v_*|,\cos\theta) \mathrm{d}\sigma\sim \eps^{-2}I_1|v-v_*|^{-1}. \]
	To kill the singular factor $\eps^{-2}$, we resort to the momentum transfer \eqref{MometumTrans}. Technically if we expand the Taylor expansion of $f(v')-f(v)$ up to the second order, then we may arrive at(see \eqref{order-2-cancelation-B1}  for details) 
	\beno \langle Q_1(g,h),f\rangle&\sim& \int B_1^\eps(|v-v_*|,\cos\theta)g_*h\bigg((v'-v)\cdot \na_vf+(v'-v)\otimes(v'-v):\na_v^2f\bigg)\mathrm{d}\sigma dv_*dv
	\\&\sim & I_3\int |g_*h|(|\na f||v-v_*|^{-2}+|\na^2 f||v-v_*|^{-1})dv_*dv. \eeno
	On the one hand, this suggests that if we aim to obtain a uniform-in-$\eps$ estimate, we should deal with $Q_1$ from the non-cutoff view. However, on the other hand, this approach results in a loss of derivatives, particularly for the solution.
	
	Now we are in a position to state our main strategy and the key ideas to overcome the difficulties. The strategy can be outlined in three steps.
	
	\underline{\it Step 1.} We construct a local mild solution in $L^1\cap L^\infty$ space via the contraction mapping theorem. Here the lifespan $T_*$ depends heavily on the parameter $
	\eps$ since we deal with \eqref{BUU} from angular cutoff view. For Fermi-Dirac particles, we get the propagation of the $L^\infty$ upper bound that $f(t)\le \eps^{-3}$ for any $t\in[0,T_*]$.
	
	\underline{\it Step 2.} We prove the propagation of the regularity uniformly in $\eps$ in weighted Sobolev spaces. This is motivated by the fact \eqref{asyOP}.  We expect that the Uehling-Uhlenbeck operator \eqref{OPQUU} will behave like a diffusive operator when $\eps$ is sufficiently small. Thus the $L^2$ framework to prove the propagation of regularity is reasonable. To implement the main idea, we develop some tools as follows. 
	
	\quad $\bullet$\underline{\,Explicit formula for the change of variable.} As we state it in the before, to kill the singularity, we will use the Taylor expansion of $f(v')-f(v)$ up to the second order. Technically we will meet the intermediate points $\kappa(v)=\kappa v'+(1-\kappa)v, 
	\iota(v_*)=\iota v'_*+(1-\iota)v_*$ where $\kappa, \iota \in[0,1]$. As a result,
	the change of variables $v \to \kappa(v)$ and $v_{*} \to \iota(v_*)$ are compulsory. Since now our kernel is not a factorized form $\Phi(|v-v_{*}|) b(\cos\theta)$
	of relative velocity and deviation angle, rough treatment is not enough. For this reason, we explicitly compute the change of variables $v \to \kappa(v)$ and $v_{*} \to \iota(v_*)$ in Lemma \ref{usual-change} and carefully use it in Lemma \ref{UPQ3R}  for  $B^{\eps}_{3}$ and Lemma \ref{Technical-lemma} for  $B^{\eps}_{1}$.
	
	\quad $\bullet$\underline{\,Coercivity estimate and the cancellation lemma.} As we explain it in the before, the operator $Q_1$ is supposed to produce the dissipation.   Indeed, if $f \geq 0$,
	\ben \nonumber 
	&& \lr{Q_1(f, \pa^{\alpha}f) , \pa^\alpha f} 
	\\ \label{into-two-terms} &=& - \f12 \int B^{\eps}_1 f_* ((\pa^\alpha f)' - (\pa^\alpha f))^2
	\mathrm{d}v \mathrm{d}v_* \mathrm{d}\sigma + \f12 \int B^{\eps}_1 f_*(((\pa^\alpha f)^2)' - (\pa^\alpha f)^2) \mathrm{d}v \mathrm{d}v_* \mathrm{d}\sigma
	\\ \label{non-negativity-solution} & \leq & \f12 \int B^{\eps}_1 f_*(((\pa^\alpha f)^2)' - (\pa^\alpha f)^2) \mathrm{d}v \mathrm{d}v_* \mathrm{d}\sigma.
	\een
	The first term in \eqref{into-two-terms} is non-positive and corresponds to the coercitity of the operator.  Unfortunately, since we consider a general interaction potential, we are unable to obtain an explicit description of the coercivity which is related closely to the one from the Landau collision operator. As a result, we only make full use of the sign. To treat the second term in   \eqref{into-two-terms}, we establish the cancellation lemma(see Lemma \ref{cancellem} for details) to balance the loss of the derivative. 
	
	\quad $\bullet$\underline{\,Estimate the collision operators from the cutoff and non-cutoff perspectives.} To estimate the collision operator $Q_{UU}$ in weighted Sobolev spaces $H^{N}_{l}$, we first  remind that $Q_1, Q_2, Q_3$ and $R$ behave quite differently and 
	each of them has its own difficulty. Moreover, since the kernels $B^{\eps}_{i}$ (where $i=1,2,3$) cannot be expressed in the product form $\Phi(|v-v_{*}|) b(\cos\theta)$, it will lead to numerous technical difficulties in the analysis.  Our main approach is based on the integration of two distinct perspectives: the cutoff view and the non-cutoff view. These enable us to balance the regularity to get the uniform-in-$\eps$ estimates. We refer readers to  Sect. 2 and Sect. 3 for details. 
	
	\quad $\bullet$\underline{\,Integration by parts formulas for the penultimate order terms.} Since we have no explicit description for the dissipation mechanism, the main obstruction to prove the propagation of regularity uniformly in $\eps$ lies in the estimates for the penultimate order terms. To bound penultimate order terms,  we have to sacrifice the regularity to kill the singular factor. To balance the regularity, we borrow the idea from \cite{chaturvedi2021stability} to establish the integration by parts formulas(see Lemma \ref{integration-by-parts-formula} for details).

	\underline{\it Step 3.} According to the Sobolev embedding theorem, the uniform-in-$\eps$ estimate indicates that the $L^\infty$ upper bound of the solution can be constrained by the Sobolev norm of the initial data. This in particular evokes the continuity argument to extend the lifespan $T_*$ to be $O(1)$ which is independent of $\eps$.

	\subsection{Organization of the paper} Section 2 and Section  3 aim to obtain a precise energy estimate of $Q_{UU}^{\eps}$ in the space $H^{N}_{l}$ through a comprehensive analysis of the bi-linear operators $Q_1, Q_2, Q_3$, and the tri-linear operator $R$. In Section 4, we prove the results in Theorems \ref{gwp-F-D} and \ref{lwp-B-E} that hold uniformly in $\eps$.
	Finally, Section 5 contains the proof of Theorem \ref{mainthm}.

	\setcounter{equation}{0}
	%
	%
	
	\section{Analysis of  Uehling-Uhlenbeck operator}
	In this section, we will examine the upper bounds of $Q$ and $R$, and investigate the commutator estimates between these operators and the weight function  $W_{l}$. The operators will be considered from the perspective of both angular cutoff view and angular non-cutoff view.
	
	%
	%
	%
	%

	\subsection{Some elementary facts} In this subsection, we will introduce some fundamental formulas that are commonly employed in the analysis of the Boltzmann operator. These formulas are particularly useful for studying the Boltzmann operator from the perspective of angular non-cutoff view.
	
	\subsubsection{Taylor expansion} When evaluating the difference $f'-f$ (or $f'_*-f_*$) before and after collision, various Taylor expansions are often used. We first introduce the order-1 expansion as 
	\ben
	\label{Taylor1-order-1}
	f'-f=\int_0^1 (\na f)(\kappa(v)) \cdot (v'-v)\mathrm{d}\kappa, \quad 
	f'_*-f_*= \int_0^1 (\na f)(\iota(v_*)) \cdot (v'_*-v_*)\mathrm{d}\iota,
	\een
	where for $\kappa, \iota \in[0,1]$, the intermediate points are defined as
	\ben\label{Defkappav} \kappa(v)=\kappa v'+(1-\kappa)v,  \quad 
	\iota(v_*)=\iota v'_*+(1-\iota)v_*. \een 
	Observing that 
	$|v'-v| = |v'_*-v_*| = |v-v_*| \sin\f{\theta}{2}$,  we have 
	\[f'-f\sim C(\na f) \theta; \quad f'_*-f_*\sim C(\na f) \theta.\]

	As we emphasize in the introduction,  expansion of $f'-f$ up to the second order is compulsory. We have
	\ben
	&&\label{Taylor1}
	f'-f=(\na f)(v)\cdot(v'-v)+\int_0^1(1-\kappa) (\na^2 f)(\kappa(v)):(v'-v)\otimes(v'-v)\mathrm{d}\kappa;
	\\
	&&\label{Taylor2} f'-f=(\na f)(v')\cdot(v'-v)-\int_0^1 \kappa(\na^2 f)(\kappa(v)):(v'-v)\otimes(v'-v)\mathrm{d}\kappa.
	\een
	Thanks to the symmetry property of $\sigma$-integral, the first terms in the formulas can be computed as follows
	\ben \label{cancell1} \int B(|v-v_*|, \f{v-v_*}{|v-v_*|}\cdot\sigma) (v'-v) \mathrm{d}\sigma = \int B(|v-v_*|, \f{v-v_*}{|v-v_*|}\cdot\sigma)  \sin^{2}\f{\theta}{2} (v_* - v)  \mathrm{d}\sigma, 
	\\
	\label{cancell2} \int  B(|v-v_*|,\f{v-v_*}{|v-v_*|} \cdot \sigma)  (v'-v) h(v') \mathrm{d}\sigma \mathrm{d}v =0. \een
	We remark that the formula \eqref{cancell1} holds for fixed $v, v_*$ and \eqref{cancell2} holds for fixed $v_*$. Therefore, 
	\eqref{Taylor1} and \eqref{cancell1} lead to $O(\theta^2)$ for the quantity $\int B g_*h(f'-f)\mathrm{d}\sigma \mathrm{d}v_* \mathrm{d}v$; so do \eqref{Taylor2}  and \eqref{cancell2} for $\int B g_*h'(f'-f)\mathrm{d}\sigma \mathrm{d}v_* \mathrm{d}v$. 
	
	\subsubsection{Momentum transfer} We claim that the kernels $B^{\eps}_{i}$ defined in \eqref{DefB1}-\eqref{DefB3}   satisfy the estimate: 
	\ben \label{order-2-cancellation}
	\int B^{\eps} \sin^{2} \f{\theta}{2} \mathrm{d}\sigma \leq \int (B^{\eps}_{1} + |B^{\eps}_2| + B^{\eps}_3) \sin^{2} \f{\theta}{2} \mathrm{d}\sigma \lesssim I_{3} |v-v_{*}|^{-3}.
	\een
	Indeed, for $B^{\eps}_1$, using the change of variable $r = \eps^{-1} |v-v_*| \sin(\theta/2)$, we have
	\ben \label{order-2-cancelation-B1}
	\int B^{\eps}_1 \sin^{2} \f{\theta}{2} \mathrm{d}\sigma 
	= 8 \pi \int_{0}^{\pi/2}    
	\eps^{-4} |v-v_*|  \sin^3(\theta/2) 	
	\hphi^2\left( \eps^{-1} |v-v_*| \sin(\theta/2) \right)
	\mathrm{d}\sin(\theta/2)
	\\ \nonumber
	= 8 \pi \int_{0}^{2^{-1/2}}    
	\eps^{-4}  
	\hphi^2 \left(  \eps^{-1} |v-v_*| t \right)
	t^{3} \mathrm{d}t 
	= 8 \pi |v-v_*|^{-3} \int_{0}^{2^{-1/2} \eps^{-1} |v-v_*|}    
	\hphi^2 (r)  r^{3} \mathrm{d}r 
	\leq  8 \pi I_{3} |v-v_*|^{-3}.
	\een
	For $B^{\eps}_3$, using the fact $\sqrt{2}/2 \leq \cos \f{\theta}{2}$ for $0 \leq \theta \leq \pi/2$ and
	the change of variable $r = \eps^{-1} |v-v_*| \cos(\theta/2)$, we can 
	similarly get that $\int B^{\eps}_3 \sin^{2} \f{\theta}{2} \mathrm{d}\sigma \leq 8 \pi I_{3} |v-v_*|^{-3}$. For $B^{\eps}_2$, the desired result follows the fact that $|B^{\eps}_2| \leq B^{\eps}_1 + B^{\eps}_3$.

	\subsubsection{Estimates for the Riesz potentials} We list the following lemma without proof.
	\begin{lem} It holds that
		\ben \label{minus-1}
		\int |g_* h f||v-v_*|^{-1} 
		\mathrm{d}v \mathrm{d}v_* \lesssim \|g\|_{L^1 \cap L^2} \|h\|_{L^2} \|f\|_{L^2}.
		\een
		Let $\delta>0, s_{1}, s_2, s_3 \geq 0, s_{1}+ s_2+ s_3 = \f12 + \delta$, then
		\ben \label{minus-2}
		\int |g_* h f||v-v_*|^{-2} 
		\mathrm{d}v \mathrm{d}v_* \lesssim_{\delta} \|g\|_{H^{s_1}} \|h\|_{H^{s_2}} \|f\|_{H^{s_3}}.
		\een
		Let $\delta>0, s_{1}, s_2 \geq 0, s_{1}+ s_2 = \f12 + \delta$, then
		\ben \label{2-2-minus-1-s1-s2} 
		\int   |v-v_{*}|^{-1} |g_*|^{2}   |h|^{2} 
		\mathrm{d}v \mathrm{d}v_{*} \lesssim_{\delta} \|g\|_{H^{s_1}}^2 \|h\|_{H^{s_2}}^2.
		\een
	\end{lem}
	%
	
	\subsection{A change of variable} In order to deal with the intermediate variables $\kappa(v)$ and $\iota(v_*)$ defined in
	\eqref{Defkappav}, we derive a useful formula involving the change of variable $v \to \kappa(v)$ and $v_{*} \to \iota(v_*)$. It is quite important for the estimates of the integrals involving the kernels $B^{\eps}_i(i=1,2,3)$.
	
	\begin{lem}\label{usual-change} For $\kappa \in [0, 2]$, let us define
		\ben \label{general-Jacobean}
		\psi_{\kappa}(\theta) \colonequals  (\cos^{2}\frac{\theta}{2}+(1-\kappa)^{2}\sin^{2}\frac{\theta}{2})^{-1/2}.
		\een
		For any $0 \leq \kappa \leq 1, v_{*} \in \R^{3}$, it holds that
		\ben \label{change-of-variable}
		\int_{\R^{3}} \int_{\SS^{2}_{+}}  B(|v-v_{*}|, \cos\theta)
		f(\kappa(v))  \mathrm{d}v \mathrm{d}\sigma
		= \int_{\R^{3}} \int_{\SS^{2}_{+}} B(|v-v_{*}|\psi_{\kappa}(\theta), \cos\theta)
		f(v) \psi_{\kappa}^{3}(\theta) \mathrm{d}v \mathrm{d}\sigma.
		\een
	Here $\SS^{2}_{+}\colonequals\{\sigma\in\SS^2|(v-v_{*}) \cdot \sigma \geq 0\}$,
		For any $0 \leq \kappa, \iota \leq 1$, it holds that
		\ben \label{change-of-variable-2}
		&& \int_{\R^{3}} \int_{\R^{3}} \int_{\SS^{2}_{+}}  B(|v-v_{*}|, \cos\theta)
		g(\iota(v_{*})) f(\kappa(v))  \mathrm{d}v \mathrm{d}v_*  \mathrm{d}\sigma 
		\\ \nonumber &=& \int_{\R^{3}} \int_{\R^{3}} \int_{\SS^{2}_{+}} B(|v-v_{*}|\psi_{\kappa+\iota}(\theta), \cos\theta)
		g(v_{*}) f(v) \psi_{\kappa+\iota}^{3}(\theta) \mathrm{d}v \mathrm{d}v_*  \mathrm{d}\sigma .
		\een
		
	\end{lem}
	\begin{proof} Recalling \eqref{Defkappav},
		we set $\cos\beta_{\kappa}\colonequals\sigma\cdot (\kappa(v)-v_{*})/|\kappa(v)-v_{*}|$. 
		To express $\beta_{\kappa}$ in terms of $\theta$, we notice that
		$ \kappa(v) - v_{*} = v^{\prime}-v_{*} + (\kappa-1)(v^{\prime}-v),$ which implies that
		\beno |\kappa(v) - v_{*}|^{2} = |v^{\prime}-v_{*}|^{2} + (\kappa-1)^{2}|v^{\prime}-v_{*}|^{2} = (\cos^{2}\frac{\theta}{2}+(1-\kappa)^{2}\sin^{2}\frac{\theta}{2})|v-v_{*}|^{2}=\psi_{\kappa}(\theta)^2|v-v_{*}|^{2}.\eeno
		From this together with fact that 
		$\left( \kappa(v) - v_{*} \right) \cdot \sigma = \left( \cos^{2}\frac{\theta}{2}+(\kappa-1)\sin^{2}\frac{\theta}{2} \right)|v-v_{*}|$, 
		we have \[
		\cos \beta_{\kappa}  =  \frac{ \cos^{2}\frac{\theta}{2}+(\kappa-1)\sin^{2}\frac{\theta}{2}}
		{\left( \cos^{2}\frac{\theta}{2}+(1-\kappa)^{2}\sin^{2}\frac{\theta}{2} \right)^{1/2} } = \varphi_{\kappa}(\sin\frac{\theta}{2}),
		\]
		where 
		$\varphi_{\kappa}(x) =
		\frac{ 1 - x^{2}+(\kappa-1)x^2}
		{\left( 1 - x^{2} +(1-\kappa)^{2}x^{2} \right)^{1/2} }.$
		The above relation yields that  if $0 \leq \theta \leq \f{\pi}{2}$, then 
		$0 \leq \beta_{\kappa} \leq \delta_{\kappa} \colonequals  
		\arccos( \f{\sqrt{2}} {2} \frac{ \kappa }
		{\sqrt{1 + (1-\kappa)^{2}  }  }) $
		is a bijection.  
		
		Now we are in a position to prove \eqref{change-of-variable}. By the fact that 
		\ben
		\label{Jacobi} \det (\frac{\partial u}{\partial v})  = (1-\frac{\kappa}{2})^2 \left( (1-\frac{\kappa}{2})+\frac{\kappa}{2} \cos\theta \right) \colonequals  \alpha_{\kappa}(\theta),
		\een  
		we get that
		\beno
		\int_{\R^{3}} \int_{\SS^{2}_{+}}  B(|v-v_{*}|, \cos\theta)
		f(\kappa(v))  \mathrm{d}v \mathrm{d}\sigma
		= 2 \pi \int_{\R^{3}} \int_{0}^{\delta_{\kappa}}  B(|v-v_{*}|\psi_{\kappa}(\theta), \cos\theta)
		f(v)  \alpha_{\kappa}^{-1}(\theta) \sin \beta_{\kappa} \mathrm{d}v \mathrm{d}\beta_{\kappa}.\eeno
		Then the desired result follows the computation  
		\[\sin \beta_{\kappa} \mathrm{d}\beta_{\kappa} = - \mathrm{d} \cos \beta_{\kappa} =  -\f{1}{4} \varphi^{\prime}_{\kappa}(\sin\f{\theta}{2}) \sin^{-1}\f{\theta}{2}
		\sin\theta \mathrm{d}\theta,\quad -\f{1}{4}  \varphi^{\prime}_{\kappa}(\sin\f{\theta}{2}) \sin^{-1}\f{\theta}{2} \alpha_{\kappa}^{-1}(\theta) =   \psi_{\kappa}^{3}(\theta).\]

		As for \eqref{change-of-variable-2}, the case  $\kappa = \iota =1$ is obviously given by the natural change of variable $(v, v_*, \sigma) \to (v', v'_*, \sigma')$ where $\sigma' = (v - v_*)/|v - v_*|$.   If  $\kappa + \iota < 2$,
		we can similarly repeat the above derivation with $\kappa$ replaced by $\kappa + \iota$. Indeed, one can derive 
		\beno
		\det (\frac{\partial (\kappa(v), \iota(v_*))}{\partial (v, v_*)}) =  \alpha_{\kappa+\iota}(\theta), \quad |v-v_{*}| =  |\kappa(v) - \iota(v_*)| \psi_{\kappa+\iota}(\theta).
		\eeno 
		Let $\beta_{\kappa+\iota}$ be the angle between $\kappa(v) - \iota(v_*)$ and $\sigma$, then $
		\cos \beta_{\kappa+\iota}  = \varphi_{\kappa+\iota}(\sin\frac{\theta}{2})$.
		If  $\kappa + \iota < 2$, then $\delta_{\kappa+\iota}>0$ 
		and
		the function:
		$
		\theta \in [0, \f{\pi}{2}]  \to 
		\beta_{\kappa+\iota} \in [0,\delta_{\kappa + \iota}] 
		$
		is a bijection. These facts are enough to obtain \eqref{change-of-variable-2} for $\kappa + \iota < 2$.
	\end{proof} 	
	
	\subsection{Integrals involving $B^{\eps}_3$} We first derive the upper bound of the integrals involving $B^{\eps}_3$ from the cutoff perspective. We remark that in this situation the estimates depend on $\eps$.
	\begin{lem} 
		Let $a \geq 0$, then
		\ben \label{B-3-a-order-velocity}
		\int B^{\eps}_3 |v-v_*|^a |g_* h| \mathrm{d}V \leq 8 \pi (\sqrt{2})^{(a-1)_{+}} \eps^{a-3} I_a 
		\|g\|_{L^1}
		\|h\|_{L^1}.
		\een
	\end{lem}
	\begin{proof}
		As $\sqrt{2}/2 \leq \cos (\theta/2) \leq 1$, we have
		\ben\quad \label{B3-integral-independent-of-v}
		&&\int |z|^{a}  
		B^{\eps}_3(z, \sigma) \mathrm{d}\sigma
		= 8 \pi \eps^{-4}|z|^{a+1}
		\int_0^{\pi/2} \hat{\phi}^2(\eps^{-1}|z|\cos \f{\theta}{2})
		\cos \f{\theta}{2}\mathrm{d}\cos \f{\theta}{2}
		\\ \nonumber  &&\leq  8 \pi (\sqrt{2})^{(a-1)_{+}}\eps^{-4}|z|^{a+1}
		\int_0^{\pi/2} \hat{\phi}^2(\eps^{-1}|z|\cos \f{\theta}{2})
		\cos^a \f{\theta}{2}\mathrm{d}\cos \f{\theta}{2}
		= 8 \pi (\sqrt{2})^{(a-1)_{+}} \eps^{a-3}
		\int_{\eps^{-1}|z|/\sqrt{2}}^{\eps^{-1}|z|} \hat{\phi}^2(t)
		t^a \mathrm{d}t
		\\ \nonumber &&\leq 8 \pi (\sqrt{2})^{(a-1)_{+}} \eps^{a-3} I_a
		\lesssim   \eps^{a-3}
		\int_0^{\infty} \hat{\phi}^2(r)
		r^a \mathrm{d}r,
		\een
		which implies \eqref{B-3-a-order-velocity}.
	\end{proof}
	
	By taking $a=0$ and replacing $\cos\f{\theta}{2}$ by $\sin\f{\theta}{2}$  in \eqref{B3-integral-independent-of-v}, we have
	\ben \label{upper-bound-of-A-eps}
	A^{\eps} \colonequals  \sup_{z \in \R^3} \int
	B^{\eps}(z, \sigma) \mathrm{d}\sigma \leq 2\sup_{z \in \R^3} \int
	B^{\eps}_1(z, \sigma) \mathrm{d}\sigma + 2\sup_{z \in \R^3} \int
	B^{\eps}_3(z, \sigma) \mathrm{d}\sigma
	\lesssim \eps^{-3} 
	I_{0}.
	\een
	The above inequality  shows that the $L^{\infty}$-norm of $\int
	B^{\eps}_1(\cdot, \sigma) \mathrm{d}\sigma$ and $\int
	B^{\eps}_3(\cdot, \sigma) \mathrm{d}\sigma$ is bounded by $\eps^{-3} 
	I_{0}$ which tends to $\infty$ as  $\eps \to 0$. Considering the $L^{1}$-norm of $\int
	B^{\eps}_3(\cdot, \sigma) \mathrm{d}\sigma$, we find it is bounded uniformly in $\eps$  by the following computation:
	\ben \label{bounded-l1-norm-simple}\qquad
	\iint B^{\eps}_3(z, \sigma) \mathrm{d}\sigma \mathrm{d}z = 4 \pi \int_{0}^{\infty} \int_{\SS^{2}_{+}} \eps^{-4}r^{3}  \hphi^2(\eps^{-1}r \cos(\theta/2))
	\mathrm{d}r \mathrm{d}\sigma 
	=4 \pi I_{3} \int_{\SS^{2}_{+}} \cos^{-4}(\theta/2)
	\mathrm{d}\sigma = 16\pi^{2}  I_{3}.
	\een
	Here $\SS^{2}_{+}$ stands for $0 \leq \theta \leq \pi/2$.

	Based on the above uniform $L^{1}$ upper bound, we can easily obtain
	uniform-in-$\eps$ estimates for various integrals involving $B^{\eps}_3$ 
	with the change of variables in \eqref{change-of-variable-2}.
	\begin{lem}\label{UPQ3R} 
		Fix $\kappa \in [0,1]$, either $u = \kappa(v_{*})$ or $u = \kappa(v)$. Then
		\ben \label{general-L1}
		\int B^{\eps}_3 |g(u)| \mathrm{d}V \lesssim I_{3} \|g\|_{L^1}.
		\een
		As a direct result, fix an integer $k \geq 2$ and $\iota_{i}, \kappa_{i} \in [0,1]$ for $1 \leq i \leq k$, 
		let $u_i \in \{\iota_{i}(v_{*}), \kappa_{i}(v) : 1\leq  i \leq k\}$ for $1 \leq i \leq k$, 
		then 
		\ben \label{general-k-functions-infty-norm}
		\int B^{\eps}_3
		\prod_{i=1}^{k}
		|f_i(u_i)| \mathrm{d}V 
		\lesssim I_{3} \prod_{i=1}^{k} \|f_i\|_{X_i},
		\een
		where two of $X_i$ are taken by $L^2$-norm and the others are taken by $L^\infty$-norm. Let $0 \leq s_i < \f32$ for $1 \leq i \leq k$ and $\sum_{i=1}^{k} s_i = \f{3k}{2} - 3$, then
		\ben \label{general-k-functions-h-norm}
		\int B^{\eps}_3
		\prod_{i=1}^{k}
		|f_i(u_i)| \mathrm{d}V 
		\lesssim_{s_1, \cdots, s_k} I_{3} \prod_{i=1}^{k} \|f_i\|_{H^{s_i}},
		\een
	\end{lem}
	
	\begin{proof}  Applying \eqref{change-of-variable-2}, we have
		$ \int B^{\eps}_3 |g(u)| \mathrm{d}V = \int  J_{\kappa, \eps}(v-v_{*}) |g(v)|  \mathrm{d}v \mathrm{d}v_{*}$,   
		where 
		\ben \label{def-J-kappa-eps}
		J_{\kappa, \eps}(z) \colonequals  \int_{\SS^{2}_{+}} \eps^{-4}|z|\psi^{4}_{\kappa}(\theta) \hphi^2(\eps^{-1}|z|\psi_{\kappa}(\theta)\cos(\theta/2))
		\mathrm{d}\sigma.
		\een
		Similarly to \eqref{bounded-l1-norm-simple},
		it is easy to see that $L^{1}$-norm of $J_{\kappa, \eps}(z)$
		is bounded (uniformly in $\kappa, \eps$)  as follows:
		\ben \label{bounded-l1-norm}
		\|J_{\kappa, \eps}\|_{L^{1}}  &=& 4 \pi \int_{0}^{\infty} \int_{\SS^{2}_{+}} \eps^{-4}r^{3} \psi^4_{\kappa}(\theta) \hphi^2(\eps^{-1}r\psi_{\kappa}(\theta)\cos(\theta/2))
		\mathrm{d}r \mathrm{d}\sigma 
		\\ \nonumber 
		&=& 4 \pi I_{3} \int_{\SS^{2}_{+}} \cos^{-4}(\theta/2)
		\mathrm{d}\sigma = 16\pi^{2}  I_{3},
		\een
		which yields \eqref{general-L1}. As a direct result, \eqref{general-k-functions-infty-norm} is easily followed by  H{\"o}lder's inequality.

		To prove \eqref{general-k-functions-h-norm},  for $2 \leq p_i < \infty$ and $\sum_{i=1}^{k} p_i^{-1} = 1$, we have
		\beno 
		\int B^{\eps}_3
		\prod_{i=1}^{k}
		|f_i(u_i)| \mathrm{d}V  \lesssim \prod_{i=1}^{k} \left(\int B^{\eps}_3
		|f_i(u_i)|^{p_i} \mathrm{d}V \right)^{1/p_i}
		\lesssim I_{3} \prod_{i=1}^{k} \|f_i\|_{L^{p_i}}
		\lesssim_{s_1, \cdots, s_k} I_{3} \prod_{i=1}^{k} \|f_i\|_{H^{s_i}},
		\eeno
		where $\f{s_i}{3} = \f{1}{2} - \f{1}{p_i}$ thanks to the Sobolev embedding theorem.
	\end{proof}
	
	With the estimates in Lemma \ref{UPQ3R}, we derive  upper bounds of $Q_3$ in weighted Sobolev.

	\begin{prop}\label{Q-3} Let $l \geq 0, \delta>0$.
		For $0 \leq  s_1, s_2, s_3$ with $s_1 + s_2 + s_3= \f{3}{2} + \delta$,
		\ben \label{Q-3-upper-bound}
		|\lr{Q_3(g,h), W_{l}f}| \lesssim_{l, \delta} I_{3} \|W_{l}g\|_{H^{s_1}} \|W_{l}h\|_{H^{s_1}} \|f\|_{H^{s_3}}.
		\een
	\end{prop}
	\begin{proof} 	For $l \geq 0$ and $0 \leq \iota_{1}, \kappa_{1}, \iota_{2}, \kappa_{2} \leq 1$,
		it is easy to check that
		\ben \label{weight-formula-1}
		W_l(\kappa_1(v)) + W_{l}(\iota_1(v_*)) \lesssim_l W_l(\kappa_2(v)) + W_{l}(\iota_2(v_*)).
		\een
		As a result, we get that
		\ben \label{direct-weighted}
		|\lr{Q_3(g,h), W_{l}f}| \lesssim_{l} \int B^{\eps}_3 (|(W_{l}g)'_* (W_{l}h)'| + |(W_{l}g)_* W_{l}h|) |f| \mathrm{d}V.
		\een
		Then the desired result follows from   \eqref{general-k-functions-h-norm}.
	\end{proof}
	
	By taking $\delta = \f12$ in Proposition \ref{Q-3},
	we easily close the energy estimate for $Q_3$ in $H^{N}_{l}$.
	\begin{lem}\label{Q3EN} Let $l \geq 0, N \geq 2$ and $m = |\alpha| \le N$. Then 
		\beno
		\sum_{\alpha_1+\alpha_2 = \alpha}|\lr{Q_3(\pa^{\alpha_1}g,\pa^{\alpha_2}f) W_{l}, W_{l}\pa^\alpha f}| \lesssim_{N,l} I_{3} \|g\|_{H^{N}_{l}} \|f\|_{H^{N}_{l}}^2.
		\eeno
	\end{lem}
	\begin{proof} 
		If $|\alpha_1| \ge 2$, then $|\alpha_2| \leq m -2$. Take $\delta = \f12$ in Proposition \ref{Q-3}.
		We take $s_1 = s_3=0$ and $s_2 = 2$ in \eqref{Q-3-upper-bound}
		to get
		\beno
		|\lr{Q_3(\pa^{\alpha_1}g,\pa^{\alpha_2}f) W_{l}, W_{l}\pa^\alpha f }|
		\lesssim_{l} I_{3} \|\pa^{\alpha_1}g\|_{L^{2}_l} 	\|\pa^{\alpha_2}f\|_{H^{2}_l} 
		\|\pa^\alpha f\|_{L^2_{l}} 	\lesssim_{l} I_{3} \|g\|_{H^{m}_l} 	\|f\|_{H^{m}_l}^2.
		\eeno
		If $|\alpha_1| = 1$, then $|\alpha_2| \leq m - 1$. Then the desired result follows by taking  $s_1 = s_2 = 1$ and  $s_3 =0$ in \eqref{Q-3-upper-bound} 
		Similarly argument can be applied to the case that  $|\alpha_1| = 0$. We complete the proof of the lemma.
	\end{proof}

	
	Relying on more regularity, we can get the weighted upper bound of $Q_3$ with a small factor $\eps^\vartheta$. Such estimates will be used in the last section to derive the asymptotic formula in Theorem \ref{mainthm}.
	\begin{prop} 
		Let $\vartheta \in [0,1]$, then
		\ben \label{Q-3-vartheta}
		|\lr{Q_3(g,h), W_l f}| \lesssim \eps^\vartheta I_{3+\vartheta} \|W_l g\|_{H^{\f34 + \f{\vartheta}{2}}} \|W_l h\|_{H^{\f34 + \f{\vartheta}{2}}} \|f\|_{L^2}.
		\een
	\end{prop} 
	\begin{proof} Recalling \eqref{direct-weighted}, by H{\"o}lder's inequality and the change of variable \eqref{change-of-variable-2}($\kappa =\iota =1$), 
		we have
		\[ |\lr{Q_3(g,h), W_{l}f}| 
		\lesssim \left( \int 
		B^{\eps}_3  |v-v_{*}|^{-\vartheta} |(W_{l}g)_*|^4  \mathrm{d}V\right)^{1/4} \left( \int 
		B^{\eps}_3  |v-v_{*}|^{-\vartheta} |W_{l}h|^4 \mathrm{d}V\right)^{1/4} \left( \int 
		B^{\eps}_3 |v-v_{*}|^{\vartheta} |f|^2 \mathrm{d}V\right)^{1/2}
		\]\[  \le \| 
		(J_{0,\eps}  |\cdot|^{-\vartheta}) * (W_{l}^4g^4)  \|_{L^1}^{1/4}
		\| 
		(J_{0,\eps}  |\cdot|^{-\vartheta}) *  (W_{l}^4h^4)  \|_{L^1}
		\| 
		(J_{0,\eps}  |\cdot|^{\vartheta}) * f^2  \|_{L^1}^{1/2},\]
		where we use the notation \eqref{def-J-kappa-eps}.
		Similarly to \eqref{bounded-l1-norm}, we derive that
		\ben \label{L1-vartheta}
		\|J_{0, \eps}|\cdot|^{\vartheta}\|_{L^{1}} = 4 \pi \int_{0}^{\infty} \int_{\SS^{2}_{+}} \eps^{-4}r^{3+\vartheta}  \hphi^2(\eps^{-1}r\cos(\theta/2))
		\mathrm{d}r \mathrm{d}\sigma 
		\lesssim  \eps^{\vartheta} I_{3+\vartheta},
		\een
		from which together with the Hardy's inequality $\int |v-v_{*}|^{-2\vartheta}|F(v)|^4  \mathrm{d}v \lesssim \|F^2\|_{H^{\vartheta}}^2$, we get 
		\beno 
		|\lr{Q_3(g,h), W_{l}f}| \lesssim \eps^{\vartheta} I_{3+\vartheta}
		\|W_{l}^2g^2  \|_{H^{\vartheta}}^{1/2}
		\|W_{l}^2h^2  \|_{H^{\vartheta}}^{1/2}
		\|f  \|_{L^2}.
		\eeno
		Using the fact that $\|F^2\|_{H^{\vartheta}} \lesssim \|F\|_{H^{3/4+\vartheta/2}}^2$, we conclude the desire result \eqref{Q-3-vartheta}.
	\end{proof}
	
	\subsection{Cancellation Lemma}  
	In this subsection, we prove the  cancellation lemma for $Q_1$ and $Q_2$ which  is used to transfer the regularity from one function to the other. 
	
	\begin{lem}[Cancellation Lemma]\label{cancellem} Let $\delta>0$ and  $a,b,c \geq 0$ verifying that $a+b+c =\f32+\delta$. For $i=1,2$, and  functions $g, h, f$,  we set
		$\mathcal{I}_i \colonequals \int B^{\eps}_i(|v-v_*|,\cos\theta)g_*\big((hf)'-hf\big)\mathrm{d}V.$ Then

		(i). For $\mathcal{I}_1$, it holds that
		\ben \label{convolution-form}
		\mathcal{I}_1=\int (J_{\eps} * g)(v) h(v)f(v)\mathrm{d}v.
		\een
		where $J_{\eps}(u)= 8 \pi \int_{\f{\sqrt{2}}2}^1 \eps^{-4} |u| \hphi^2 (\eps^{-1}|u|r)r \mathrm{d}r$ and
		$\|J_{\eps}\|_{L^1}= 16 \pi^{2} I_{3}$.  As a result,  we have
		\ben \label{Q1-cancealltion-g-h-f}
		|\mathcal{I}_1| 
		\lesssim_{\delta} I_{3} \|g\|_{H^{a}} \|h\|_{H^{b}}\|f\|_{H^{c}}.
		\een
		(ii). For $\mathcal{I}_2$, it holds that
		\ben \label{convolution-form-B2}
		\mathcal{I}_2=\int (K_{\eps} * g)(v) h(v)f(v)\mathrm{d}v,
		\een
		where $K_{\eps}(u) = K_{\eps,1}(u) + K_{\eps,2}(u)$ with
		\ben \label{kernel-K-1}
		K_{\eps,1}(u) &=& 
		16 \pi \int_0^{\f{\sqrt{2}}2} \eps^{-4} |u| \hphi(\eps^{-1}|u|r) \left( \hphi(\eps^{-1}|u|) - \hphi(\eps^{-1}|u|\sqrt{1-r^2}) \right) 
		r \mathrm{d}r,
		\\ \label{kernel-K-2}
		K_{\eps,2}(u)   &=& 
		16 \pi \int_{\f{\sqrt{2}}2}^{1} \eps^{-4} |u| \hphi(\eps^{-1}|u|r) \hphi(\eps^{-1}|u|) r \mathrm{d}r.
		\een
		Moreover, $\|K_{\eps}\|_{L^1} \leq  64 \pi^{2} (I_{3} + I^{\prime}_{3})$ which implies that 
		\ben  \label{Q2-cancealltion-g-h-f}
		|\mathcal{I}_2|
		\lesssim_{\delta} (I_{3} + I^{\prime}_{3})  \|g\|_{H^{a}} \|h\|_{H^{b}}\|f\|_{H^{c}}.
		\een
		In general,  for $0 \leq \vartheta \leq 1$,
		$\|K_{\eps}|\cdot|^{\vartheta}\|_{L^1} \lesssim \eps^{\vartheta} (I_{3+\vartheta} + I^{\prime}_{3+\vartheta})$ and 
		\ben  \label{eps-vartheta-B2-cancellation}
		|\mathcal{I}_2|
		\lesssim_{\delta}  \eps^{\vartheta} (I_{3+\vartheta} + I^{\prime}_{3+\vartheta}) \|g\|_{H^{\f{3}{4} + \f{\vartheta}{2} + \delta}}  \|h\| _{H^{\f{3}{4} + \f{\vartheta}{2} + \delta}} \|f\|_{L^{2}}.
		\een
	\end{lem}
	
	\begin{proof} 
		We first prove the estimate of $\mathcal{I}_1$. By \eqref{change-of-variable}, we have
		\beno
		\mathcal{I}_1&=& 2 \pi  \int \big(B^{\eps}_1(\f{|v-v_*|}{\cos(\theta/2)},\cos\theta)(\cos(\theta/2))^{-3}-B^{\eps}_1(|v-v_*|,\cos\theta)\big)  g_* hf \sin\theta \mathrm{d}\theta \mathrm{d}v_* \mathrm{d}v \\ 
		&=& 8 \pi \int \int_0^{\f{\sqrt{2}}2} \eps^{-4} |v-v_*| \big[\hphi^2(\eps^{-1}|v-v_*|\f{r}{\sqrt{1-r^2}})(1-r^2)^{-2}-\hphi^2(\eps^{-1}|v-v_*|r)\big]  g_* hf r \mathrm{d}r \mathrm{d}v_* \mathrm{d}v.
		\eeno
		By the change of variable $\mathcal{r}\colonequals \f{r}{\sqrt{1-r^2}}$ which implies $(1-r^2)^{-2} r \mathrm{d}r = \mathcal{r}  \mathrm{d}\mathcal{r} $, we get 
		\beno
		\mathcal{I}_1
		= 8 \pi \int \int_{\f{\sqrt{2}}2}^{1} \eps^{-4} |v-v_*| \hphi^2(\eps^{-1}|v-v_*|r)  g_* hf r \mathrm{d}r \mathrm{d}v_* \mathrm{d}v,
		\eeno
		which is exactly \eqref{convolution-form}. Since $J_{\eps}$ is radial, let $\mathcal{r}\colonequals |u|$,
		\beno
		\|J_{\eps}\|_{L^1} = 32 \pi^{2} \int_{0}^\infty \int_{\f{\sqrt{2}}2}^1 \eps^{-4} \mathcal{r}^{3} \hphi^2 (\eps^{-1} \mathcal{r} r) 
		r \mathrm{d}r \mathrm{d}\mathcal{r}
		=32 \pi^{2} \left(\int_{0}^\infty  s^{3} \hphi^2 (s) \mathrm{d}s \right) \left( \int_{\f{\sqrt{2}}2}^1 r^{-3} dr \right)
		= 16 \pi^{2} I_{3}. 
		\eeno

		%
		\smallskip
		
		We turn to the estimate of $\mathcal{I}_2$. Following the same argument in the above, we can get the formula \eqref{convolution-form-B2} with $K_{\eps}$ as the sum of \eqref{kernel-K-1} and \eqref{kernel-K-2}. Let us compute the $L^{1}$-norm of $K_{\eps,1}$ and $K_{\eps,2}$.
		Let $\mathcal{r}\colonequals |u|$, then
		\beno
		&&\|K_{\eps,2}\|_{L^1} = 64 \pi^{2} \int_{0}^\infty \left|\int_{\f{\sqrt{2}}2}^1 \eps^{-4} \mathcal{r}^{3} \hphi(\eps^{-1}r\mathcal{r}) \hphi(\eps^{-1}\mathcal{r}) r \mathrm{d}r \right| \mathrm{d}\mathcal{r}
		\\ 
		&&\le  64 \pi^{2} \int_{\f{\sqrt{2}}2}^1  \left(\int_{0}^\infty  \eps^{-4} \mathcal{r}^{3} \hphi^{2}(\eps^{-1}r\mathcal{r}) \mathrm{d}\mathcal{r} \right)^{1/2} \left(\int_{0}^\infty  \eps^{-4} \mathcal{r}^{3} \hphi^{2}(\eps^{-1}\mathcal{r}) \mathrm{d}\mathcal{r} \right)^{1/2} r \mathrm{d}r  \leq 32 \pi^{2}  I_{3}, 
		\eeno
		where for fixed $r$, we used the change of variables $\mathcal{r} \to \eps^{-1} r\mathcal{r} $ and $\mathcal{r} \to \eps^{-1} \mathcal{r}$.
		
		To kill the singularity at $r=0$ in $K_{\eps,1}$, by Taylor expansion, we have 
		\beno 
		K_{\eps,1}(u) = 
		16 \pi \int_0^{\f{\sqrt{2}}2} \int_{0}^{1} \eps^{-5} |u|^{2} \hphi(\eps^{-1}|u|r)  
		(\hphi)^{\prime}(\eps^{-1}|u|(\tau + (1-\tau)\sqrt{1-r^2})) (1-\sqrt{1-r^2}) 
		r \mathrm{d}r \mathrm{d}\tau,
		\eeno
		which implies that
		\beno
		&&\|K_{\eps,1}\|_{L^1} = 64 \pi^{2} \int_{0}^\infty \left|
		\int_0^{\f{\sqrt{2}}2} \int_{0}^{1} \eps^{-5} \mathcal{r}^{4} \hphi(\eps^{-1}\mathcal{r}r)  
		(\hphi)^{\prime}(\eps^{-1}\mathcal{r}(\tau + (1-\tau)\sqrt{1-r^2})) (1-\sqrt{1-r^2})
		r \mathrm{d}r \mathrm{d}\tau \right| \mathrm{d}\mathcal{r}
		\\ 
		&&\le  64 \pi^{2} \int_0^{\f{\sqrt{2}}2} \int_{0}^{1} \left(\int_{0}^\infty  \eps^{-4} \mathcal{r}^{3} \hphi^{2}(\eps^{-1}r\mathcal{r}) \mathrm{d}\mathcal{r}\right)^{1/2} \left(\int_{0}^\infty  \eps^{-6} \mathcal{r}^{5} |(\hphi)^{\prime}(\eps^{-1}\mathcal{r}(\tau + (1-\tau)\sqrt{1-r^2}))|^{2} \mathrm{d}\mathcal{r} \right)^{1/2} \\
		&&\times
		(1-\sqrt{1-r^2})                                  r \mathrm{d}r \mathrm{d}\tau
		\le 64 \pi^{2} \left(\int_{0}^\infty   s^{3} \hphi^{2}(s) \mathrm{d}s \right)^{1/2} \left(\int_{0}^\infty  s^{5} |(\hphi)^{\prime}(s)|^{2} \mathrm{d}s \right)^{1/2}   \int_0^{\f{\sqrt{2}}2} \int_{0}^{1} (1-\sqrt{1-r^2}) r^{-1}
		\\&&\times  (\tau + (1-\tau)\sqrt{1-r^2})^{-3} dr \mathrm{d}\tau
		\leq 
		8 \sqrt{2} \pi^{2} (I_{3} + I^{\prime}_{3}), 
		\eeno
		where  the estimates 
		$1-\sqrt{1-r^2} \leq r^{2}/2$, $\sqrt{2}/2 \leq \tau + (1-\tau) \sqrt{1-r^2} \leq 1$ are used. Now we  have 
		\beno
		\|K_{\eps}\|_{L^1} \leq		\|K_{\eps,1}\|_{L^1} + \|K_{\eps,2}\|_{L^1} \leq 8 \sqrt{2} \pi^{2} (I_{3} + I^{\prime}_{3}) + 32 \pi^{2}  I_{3} \leq 64 \pi^{2} (I_{3} + I^{\prime}_{3}),
		\eeno
		which gives \eqref{Q2-cancealltion-g-h-f}.
		The same argument can be applied to get that	$\|K_{\eps}|\cdot|^{\vartheta}\|_{L^1} \lesssim \eps^{\vartheta} (I_{3+\vartheta} + I^{\prime}_{3+\vartheta})$, which gives \eqref{eps-vartheta-B2-cancellation}. We end the proof.
	\end{proof}

	\subsection{Integrals involving $B^{\eps}_1$} We shall use \eqref{change-of-variable-2} to give the estimate of the integrals involving $B^{\eps}_1$.

	\begin{lem}\label{Technical-lemma} Let $a, b \in \R$ with $b \geq 0$. If $0 \leq \kappa, \iota \leq 1$, then
		\ben \qquad\label{Q1-result-general}
		\int B^{\eps}_1  | g(\iota(v_*)) h(\kappa(v))|  |v-v_*|^{a+b} \sin^{b-1}(\theta/2) \mathrm{d}V \leq (\sqrt{2})^{(a+1)_{+}} 8 \pi 
		\eps^{b-3}  I_{b}  \int |v-v_*|^{a} | g_* h| \mathrm{d}v_* \mathrm{d}v.
		\een
	\end{lem}
	\begin{proof} Applying \eqref{change-of-variable-2}, we have 
		\beno 
		&& \int B^{\eps}_1  | g(\iota(v_*)) h(\kappa(v))|  |v-v_*|^{a+b} \sin^{b-1}(\theta/2) \mathrm{d}V 
		\\ &=& 
		\int    
		\eps^{-4} |v-v_*|^{a+b+1} \psi^{a+b+4}_{\kappa+\iota}(\theta) \sin^{b-1}(\theta/2) 	
		\hphi^2\bigg( \eps^{-1} |v-v_*| \psi_{\kappa+\iota}(\theta) \sin(\theta/2)
		\bigg) | g_* h|
		\mathrm{d}V.
		\eeno 
		Recalling \eqref{general-Jacobean}, if we set $\mathcal{r}: = \psi_{\kappa+\iota}(\theta) \sin(\theta/2)$ then $\mathrm{d} \mathcal{r} = \psi_{\kappa+\iota}^{3}(\theta) \mathrm{d}\sin(\theta/2)$. Since $1 \leq \psi_{\kappa+\iota} \leq \sqrt{2}$, we have
		\ben \label{explicit-computation}
		&& \int    
		\eps^{-4}  \psi^{a+b+4}_{\kappa+\iota}(\theta) \sin^{b-1}(\theta/2) 	
		\hphi^2\bigg( \eps^{-1} |v-v_*| \psi_{\kappa+\iota}(\theta) \sin(\theta/2)
		\bigg)  \mathrm{d}\sigma 
		\\ \nonumber 
		&=& 8 \pi \int_{0}^{\pi/2}    
		\eps^{-4}  \psi^{a+b+4}_{\kappa+\iota}(\theta) \sin^b(\theta/2) 	
		\hphi^2\bigg( \eps^{-1} |v-v_*| \psi_{\kappa+\iota}(\theta) \sin(\theta/2)
		\bigg)  \mathrm{d}\sin(\theta/2)
		\\ \nonumber
		&=& 8 \pi \int_{0}^{(2-2(\kappa+\iota)+(\kappa+\iota)^{2})^{-1/2}}  
		\eps^{-4}  \psi^{a+1}_{\kappa+\iota}(\theta) 	
		\hphi^2\bigg( \eps^{-1} |v-v_*| \mathcal{r}
		\bigg)  \mathcal{r}^{b} \mathrm{d}\mathcal{r} 
		\\ \nonumber
		&\leq& 
		(\sqrt{2})^{(a+1)_{+}} 8 \pi \int_{0}^{1}    
		\eps^{-4}  
		\hphi^2\bigg( \eps^{-1} |v-v_*| \mathcal{r}
		\bigg)  \mathcal{r}^{b} \mathrm{d}\mathcal{r} \leq (\sqrt{2})^{(a+1)_{+}} 8 \pi
		\eps^{b-3}  I_{b} |v-v_*|^{-b-1},
		\een
		which yields \eqref{Q1-result-general}.
	\end{proof}

	\begin{rmk} \label{singular-factor-B2-B3} If we borrow the same idea used here to the estimate of integrals involving $B^{\eps}_3$,    we will use the change of variable: $\mathcal{r} = \psi_{\kappa+\iota}(\theta) \cos(\theta/2)$   which indicates that $\mathrm{d} \mathcal{r} = (1-\kappa-\iota)^2\psi_{\kappa+\iota}^{3}(\theta) \mathrm{d}\cos(\theta/2)$. As a result, the ending estimate will have a singular factor $(1-\kappa-\iota)^{-2}$ as $\kappa + \iota \to 1$. For this reason, we always avoid the change of variable $v \to \kappa(v)$ or $v_* \to \iota(v_*)$ for integrals involving $B^{\eps}_3$. 
	\end{rmk}
	
	\begin{rmk} Lemma \ref{Technical-lemma} and its proof are highly versatile, making them applicable to the majority of integrals involving $B^{\eps}_1$ that will arise in this article. To facilitate future reference, we provide several examples of their use below.
		
		$\bullet$ Using \eqref{change-of-variable} and the computation in \eqref{explicit-computation}, we have
		\ben \label{B1-sigma-vStar}
		\int
		B^{\eps}_1(|v-v_{*}|, \cos\theta)  f_*'  \mathrm{d}\sigma \mathrm{d}v_*  =
		\int
		B^{\eps}_1(|v-v_{*}|\psi_{1}(\theta), \cos\theta)  f_* \psi_{1}^{3}(\theta) \mathrm{d}\sigma \mathrm{d}v_*
		\lesssim \eps^{-3}  I_{0}
		\| f \|_{L^1} .
		\een

		$\bullet$  By taking $a=0$ in \eqref{Q1-result-general}, for $b \geq 0$,
		we have
		\ben \label{general-result-cutoff}
		\int B^{\eps}_1 | g(\iota(v_*)) h(\kappa(v))|  |v-v_*|^b 
		\sin^{b-1} \f{\theta}{2} \mathrm{d}V \leq 8 \sqrt{2} \pi  \eps^{b-3} I_{b} \|g\|_{L^1}
		\|h\|_{L^1}.
		\een

		$\bullet$ Let $0 \leq \vartheta \leq 1$.
		Let $c, d \in \R$ with $d \geq 2+\vartheta$. Thanks to the fact that
		$|v-v_*|^c \sin^d(\theta/2) \leq |v-v_*|^{c - 3-\vartheta} |v-v_*|^{3+\vartheta} \sin^{2+\vartheta}(\theta/2)$, if we 
		take $a = c - 3-\vartheta$ and $b = 3+\vartheta$ in
		\eqref{Q1-result-general},  then
		\ben \label{Q1-result-with-factor-a-minus3-with-small-facfor}\qquad
		\int B^{\eps}_1 | g(\iota(v_*)) h(\kappa(v))| |v-v_*|^c \sin^d(\theta/2) \mathrm{d}V \leq 8 \pi (\sqrt{2})^{(c-2-\vartheta)_+} \eps^\vartheta I_{3+\vartheta} \int |v-v_*|^{c-3-\vartheta} | g_* h| \mathrm{d}v_* \mathrm{d}v. 
		\een
		In particular, if $\vartheta = 0, d = 2$ in \eqref{Q1-result-with-factor-a-minus3-with-small-facfor}, then we get that
		\ben \label{Q1-result-with-factor-a-minus3-type-2}\qquad
		\int B^{\eps}_1  | g(\iota(v_*)) h(\kappa(v))| |v-v_*|^c \sin^2(\theta/2) \mathrm{d}V \leq 8 \pi (\sqrt{2})^{(c-2)_{+}} I_{3} \int |v-v_*|^{-3+c} | g_* h| \mathrm{d}v_* \mathrm{d}v.
		\een
	\end{rmk}

	Using \eqref{Q1-result-with-factor-a-minus3-type-2} properly, we derive the following flexible estimate allowing some balance of weight. 
	\begin{lem} 
		For $0 \leq \iota, \kappa_{1}, \kappa_{2} \leq 1$, it holds that
		\ben \label{general-B1-h-f}
		\int B^{\eps}_{1} |g(\iota(v_*)) h(\kappa_1(v)) f(\kappa_2(v))| |v' - v|^{2} \mathrm{d}V
		\lesssim I_{3} (\|g\|_{L^1_{a_1}} + \|g\|_{L^2}) \|h\|_{L^2_{a_2}} \|f\|_{L^2_{a_3}},
		\een
		where $-1 \leq a_1 \leq 0 \leq a_2, a_3$ and $a_1+ a_2 + a_3 = 0$.
	\end{lem}
	\begin{proof} We divide the integration domain into two parts: $\mathcal{U} \colonequals  \{(v, v_*, \sigma) : |v-v_*| \leq 1, (v-v_*) \cdot \sigma \geq 0\}$ and $\mathcal{F} \colonequals  \{(v, v_*, \sigma) : |v-v_*| \ge 1, (v-v_*) \cdot \sigma \geq 0\}$ and denote the associated integrals by $\mathcal{I}_{\leq}$ and $\mathcal{I}_{\geq}$ accordingly.
		
		$\bullet$	In the domain $\mathcal{U}$, one has $|\kappa(v) - \iota(v_*)| \leq 1$ for all $0 \leq \iota, \kappa\leq 1$. 
		Using \eqref{Q1-result-with-factor-a-minus3-type-2}, we will have
		\beno 
		\mathcal{I}_{\leq} &\leq& \left(\int \mathrm{1}_{\mathcal{U}} B^{\eps}_{1}  |g(\iota(v_*)) h^2(\kappa_1(v))| |v' - v|^{2} \mathrm{d}V \right)^{1/2}
		\left(	\int  \mathrm{1}_{\mathcal{U}}  B^{\eps}_{1}|g(\iota(v_*))  f^2(\kappa_2(v))| |v' - v|^{2} \mathrm{d}V \right)^{1/2}
		\\ &\lesssim& I_{3} \left(\int \mathrm{1}_{|v-v_*| \leq 1} |g_* h^2| |v - v_*|^{-1} \mathrm{d}v \mathrm{d}v_* \right)^{1/2}
		\left(\int \mathrm{1}_{|v-v_*| \leq 1} |g_* f^2| |v - v_*|^{-1} \mathrm{d}v \mathrm{d}v_* \right)^{1/2} 
		\\ &\lesssim& I_{3} \|g\|_{L^2}\|h\|_{L^2}\|f\|_{L^2}.
		\eeno
		
		$\bullet$ In the domain $\mathcal{F}$, one has $|\kappa(v) - \iota(v_*)| \geq \sqrt{2}/2$ and
		$ W_l(v-v_*)\sim |v - v_*|^l$. Then we get that
		\beno
		|v - v_*|^{-1} &\lesssim& |v - v_*|^{a_1} = |v - v_*|^{-a_2 - a_3} \sim
		|\kappa_1(v) - \iota(v_*)|^{-a_2} |\kappa_2(v) - \iota(v_*)|^{-a_3}
		\\	&\lesssim& W_{a_1}(\iota(v_*))  W_{a_2}(\kappa_1(v)) W_{a_3}(\kappa_2(v)),
		\eeno
		which gives
		\beno 
		&&\mathcal{I}_{\geq} \lesssim  \int \mathrm{1}_{\mathcal{F}} B^{\eps}_{1}  |(W_{a_1}g)(\iota(v_*)) (W_{a_2}h)(\kappa_1(v)) (W_{a_3}f)(\kappa_2(v))| |v - v_*|^{3}\sin^2\f{\theta}{2} \mathrm{d}V
		\\	&&\lesssim\left(\int B^{\eps}_{1} |(W_{a_1}g)(\iota(v_*)) (W_{a_2}h)^2(\kappa_1(v))| |v - v_*|^{3}\sin^2\f{\theta}{2} \mathrm{d}V \right)^{1/2}
		\bigg(	\int B^{\eps}_{1}  |(W_{a_1}g)(\iota(v_*)) (W_{a_3}f)^2(\kappa_2(v))|\\&&\times |v - v_*|^{3}\sin^2\f{\theta}{2} \mathrm{d}V \bigg)^{1/2}
		\lesssim I_{3} \left(\int  |(W_{a_1}g)_* (W_{a_2}h)^2| \mathrm{d}v \mathrm{d}v_* \right)^{1/2}
		\left(\int |(W_{a_1}g)_* (W_{a_3}f)^2| \mathrm{d}v \mathrm{d}v_* \right)^{1/2} 
		\\	 &&\lesssim I_{3} \|g\|_{L^1_{a_1}} \|h\|_{L^2_{a_2}}\|f\|_{L^2_{a_2}}.
		\eeno

		We complete the proof of the lemma by patching together these two estimates. 
	\end{proof}
	
	\subsection{Upper bounds from the cutoff perspective} 	
	We will prove several upper bounds of
	the operator $Q$ and $R$. All the estimates will depend heavily on the parameter $\eps$.
	
	\begin{lem} \label{Q-inner-product-upper-bound-L2-L2} It holds that
		\ben \label{rough-upper-Q}
		|\lr{ Q (g, h), f}| \lesssim (\eps^{-3} I_{0} + I_{3})\|g\|_{L^{1} \cap L^2}	\|h\|_{L^2} \|f\|_{L^2},
		\\
		\label{rough-upper-R-rho-L-infty}
		|\lr{ R (g, h, \rho), f}| \lesssim    (I_{0} + \eps^{3} I_{3}) \|g\|_{L^{1} \cap L^2}	\|h\|_{L^2} \|\rho\|_{L^\infty} \|f\|_{L^2},
		\\
		\label{rough-upper-R-rho-L2}
		|\lr{ R (g, h, \rho), f}| \lesssim    (I_{0} + \eps^{3} I_{3})
		\|g\|_{L^{1} \cap L^\infty}	\|h\|_{L^2 \cap L^\infty} \|\rho\|_{L^2} \|f\|_{L^2}.
		\een
	\end{lem}
	\begin{proof} Observing that $B^{\eps} \leq 2 B^{\eps}_1 + 2 B^{\eps}_3$, we get 
		\ben \label{Q-inner-product}\qquad
		|\lr{ Q (g, h), f}| =  |\int B^{\eps} g_* h (f' - f) \mathrm{d}V|
		\leq 2\int B^{\eps}_1 |g_* h (|f'| + |f|)| \mathrm{d}V + 2\int B^{\eps}_3 |g_* h (|f'| + |f|)| \mathrm{d}V.
		\een

		Applying  \eqref{general-result-cutoff} with  $\iota=\kappa=b=0$ and noting 
		$0 \leq \sin \f{\theta}{2} \leq 1$, we have
		\ben \label{B1-g-star-hf}
		\int B^{\eps}_1 |g_* h f| \mathrm{d}V \lesssim 	
		\eps^{-3} I_{0} \|g\|_{L^1}
		\|h f\|_{L^1} \lesssim 	
		\eps^{-3} I_{0} \|g\|_{L^1}
		\|h\|_{L^2} 	\|f\|_{L^2}.
		\een
		Again by \eqref{general-result-cutoff} with $\iota=b=0$, $\kappa=0$ or $1$, we get the same bound for $\int B^{\eps}_1 |g_* h f'|  \mathrm{d}V$. 
		
		Next, by taking $a=0$ in \eqref{B-3-a-order-velocity}, we have
		\ben \label{B3-g-star-hf}
		\int B^{\eps}_3 |g_* h f| \mathrm{d}V \lesssim 	
		\eps^{-3} I_{0} \|g\|_{L^1}
		\|h f\|_{L^1} \lesssim 	
		\eps^{-3} I_{0} \|g\|_{L^1}
		\|h\|_{L^2} 	\|f\|_{L^2}.
		\een
		It is not difficult to check that  
		\ben\label{B3-g-star-h-f-prime}\int B^{\eps}_3 |g_* h f'|  \mathrm{d}V\le (I_0\eps^{-3} I_{3})^{\f12} \|g\|_{L^2}
		\|h\|_{L^2} 	\|f\|_{L^2}.\een We conclude  \eqref{rough-upper-Q} by 
		patching together the above estimates.
		
		For the cubic term, we first observe that 
		\ben \label{R-inner-product}
		|\lr{ R (g, h, \rho), f}| =  \eps^{3}|\int B^{\eps} g_* h (\rho' + \rho'_*) (f' - f) \mathrm{d}V|.
		\een Then   \eqref{rough-upper-R-rho-L-infty} follows \eqref{rough-upper-Q} by imposing $L^\infty$ norm on $\rho$.

		To prove \eqref{rough-upper-R-rho-L2}, by using $B^{\eps} \leq 2 B^{\eps}_1 + 2 B^{\eps}_3$, we get that
		\beno
		&&|\lr{ R (g, h, \rho), f}| =  \eps^{3} |\int B^{\eps} g_* h (\rho' + \rho'_*) (f' - f) \mathrm{d}V|
		\leq  2 \eps^{3} \|h\|_{L^\infty} \int B^{\eps}_1 |g_* \rho' (|f'| + |f|)| \mathrm{d}V \\ &&+ 2 \eps^{3}
		\int B^{\eps}_1 |g_* h \rho'_* (|f'| + |f|)| \mathrm{d}V
		+ 2 \eps^{3} \|g\|_{L^\infty} \|h\|_{L^\infty} \int B^{\eps}_3  (|\rho'| + |\rho'_*|) (|f'| + |f|)| \mathrm{d}V.	
		\eeno
		By the Cauchy-Schwartz inequality,
		using \eqref{general-result-cutoff} and \eqref{general-k-functions-infty-norm}, we conclude the desired result.	\end{proof}

	In the next, we prove the 
	commutator estimates from the cutoff perspective. 
	\begin{lem}\label{commutator-Q-cutoff} 
		Let $l \geq 2$, then
		\ben \label{rough-commutator-Q}
		|\lr{Q(g,h)W_l-Q(g,hW_l),f}| \leq C_{l} (\eps^{-3} I_{0} + I_{3}) \|g\|_{L^2_{l}}
		\|h\|_{L^2_l} \|f\|_{L^2},
		\\ \label{rough-commutator-R-rho-L-infty}
		|\lr{R(g, h, \rho)W_l - R(g, hW_l, \rho),f}|  \leq C_{l} (I_{0} + \eps^{3}  I_{3}) \|g\|_{L^2_{l}}
		\|h\|_{L^2_l} \|\rho\|_{L^\infty} \|f\|_{L^2},
		\\ \label{rough-commutator-R}
		|\lr{R(g, h, \rho)W_l - R(g, hW_l, \rho),f}| \leq C_{l} (I_{0} + \eps^{3}  I_{3}) \|g\|_{L^2_{l} \cap L^\infty_l}
		\|h\|_{L^2_l \cap L^\infty_l} \|\rho\|_{L^2} \|f\|_{L^2}.
		\een
	\end{lem}

	\begin{proof}It is easy to compute that 
		\ben\label{Commutator-formula-Q}
		\mathcal{I}\colonequals \lr{Q(g,h)W_l-Q(g,hW_l),f}=\int B^{\eps} g_* h f' (W_l'-W_l)\mathrm{d}V.
		\een
		From the fact that $|\na^{k} W_{l}| \lesssim_{l,k} W_{l-k}$ with $l \in \R, k \in \N$, 
		\eqref{weight-formula-1}  yields that
		\ben \label{weight-diff-into-two-cases}
		|W_{l}' - W_{l}| \lesssim_{l} 
		\mathrm{1}_{|v-v_*| \leq 1} W_{l} + \mathrm{1}_{|v-v_*| \geq 1} (W_{l-1} + (W_{l-1})_*) |v-v_*| \sin\f{\theta}{2}.
		\een
		Then we have
		\beno
		|\mathcal{I}|
		\lesssim_{l}  \int B^{\eps} |g_* W_l h  f'  | \mathrm{d}V +
		\int B^{\eps} |g_* W_{l-1} h f'  | |v-v_*|  \mathrm{d}V 
		+
		\int B^{\eps}   \mathrm{1}_{|v-v_*| \geq 1} |v-v_*| \sin\f{\theta}{2}
		|(W_{l-1}g)_* h f'  | \mathrm{d}V.
		\eeno
		
		Thanks to \eqref{B3-g-star-h-f-prime}, the first term can be bounded as
		$\int B^{\eps} |g_* W_l h  f'  | \mathrm{d}V 	\lesssim  (\eps^{-3} I_{0} + I_{3}) \|g\|_{L^1 \cap L^2}
		\|h\|_{L^2_l}  \|f\|_{L^2}$.
		For the second term, using
		$B^{\eps} \leq 2B^{\eps}_1 + 2B^{\eps}_3$, it suffices to estimate $\int B^{\eps}_1 |g_* W_{l-1} h f'  |  |v-v_*|  \mathrm{d}V$ and $
		\int B^{\eps}_3 |g_* W_{l-1} h f'  |  |v-v_*| \mathrm{d}V$. By the Cauchy-Schwartz inequality, \eqref{general-result-cutoff}, \eqref{B-3-a-order-velocity} and \eqref{general-L1} imply that 
		\beno \int B^{\eps} |g_* W_{l-1} h f'  | |v-v_*|  \mathrm{d}V\lesssim (\eps^{-2}I_1+(\eps^{-1} I_2 I_{3})^{1/2}) \|g\|_{L^1 \cap L^2}
		\|h\|_{L^2_l} \|f\|_{L^2}\lesssim  (\eps^{-3} I_{0} + I_{3}) \|g\|_{L^1 \cap L^2}
		\|h\|_{L^2_l}  \|f\|_{L^2},\eeno 
		where we use the fact that $\eps^{a-3}I_a\le \eps^{-3} I_{0} + I_{3}$ for $0\le a\le3$ thanks to the interpolation. The similar argument can be applied to the third term and we   conclude the desired result  \eqref{rough-commutator-Q}.
		
		For the cubic term, we first have
		\ben \label{Commutator-formula-R}
		\mathcal{K} \colonequals  \lr{R(g, h, \rho)W_l - R(g, hW_l, \rho),f}= \pm  \eps^{3} \int B^{\eps} g_* h (\rho' + \rho'_*) f' (W_l'-W_l)\mathrm{d}V.
		\een
		By comparing the structure of \eqref{Commutator-formula-Q} and that of \eqref{Commutator-formula-R},  \eqref{rough-commutator-R-rho-L-infty} easily follows 
		\eqref{rough-upper-R-rho-L-infty}.
		
		It remains to derive \eqref{rough-commutator-R}. Note that
		\beno 
		\eps^{-3}|\mathcal{K}| \lesssim  \int B^{\eps}_1 |g_* h| (|\rho'| + |\rho'_*|) |f'| |W_l'-W_l|\mathrm{d}V
		+ \int B^{\eps}_3 |g_* h| (|\rho'| + |\rho'_*|) |f'| |W_l'-W_l|\mathrm{d}V.
		\eeno
		Since $|W_l'-W_l| \lesssim W_l + (W_l)_*$,   \eqref{general-k-functions-infty-norm} implies that 
		$\int B^{\eps}_3 |g_* h| (|\rho'| + |\rho'_*|) |f'| |W_l'-W_l|\mathrm{d}V \lesssim I_{3} \|g\|_{ L^\infty_l}
		\|h\|_{ L^\infty_l} \|\rho\|_{L^2} \|f\|_{L^2}
		.$
		
		We use \eqref{weight-diff-into-two-cases} to expand the first term by
		\beno
		&& \int B^{\eps}_1 |g_* h| (|\rho'| + |\rho'_*|) |f'| |W_l'-W_l|\mathrm{d}V
		\lesssim \int B^{\eps}_1 |g_* W_l h| (|\rho'| + |\rho'_*|) |f'| \mathrm{d}V 
		+
		\int B^{\eps}_1 |g_* W_{l-1}h| (|\rho'| + |\rho'_*|) \\	&&\times |f'| |v-v_*| \mathrm{d}V 
		+
		\int B^{\eps}_1     
		|(W_{l-1}g)_* h \rho'_* f'| |v-v_*| \mathrm{d}V
		+
		\int B^{\eps}_1   \mathrm{1}_{|v-v_*| \geq 1} |v-v_*| \sin\f{\theta}{2}
		|(W_{l-1}g)_* h \rho' f'| \mathrm{d}V
		.
		\eeno
		By repeatedly using  \eqref{general-result-cutoff} and Cauchy-Schwartz inequality, we are led to the desired result 
		\eqref{rough-commutator-R} and then complete
		the proof of the lemma.
	\end{proof}

	Lemma \ref{Q-inner-product-upper-bound-L2-L2} and \ref{commutator-Q-cutoff} together yield the following  upper bounds in weighted $L^p$ spaces.
	
	\begin{prop} \label{Q-inner-product-weighted-L2-L2} Let $l \geq 2$,
		\ben \label{rough-upper-Q-weighted}
		|\lr{ W_l Q (g, h), f}|  \leq   C_{l} (\eps^{-3} I_{0} + I_{3}) \|g\|_{L^{2}_l}	\|h\|_{L^2_l} \|f\|_{L^2},
		\\ \label{rough-upper-R-weighted-rho-L-infty}
		|\lr{ W_l R (g, h, \rho), f}|  \leq   C_{l} (I_{0} + \eps^{3} I_{3}) \|g\|_{L^{2}_l}	\|h\|_{L^2_l} \|\rho\|_{L^\infty} \|f\|_{L^2},
		\\\label{rough-upper-R-weighted}
		|\lr{ R (g, h, \rho), W_l f}| \leq   C_{l} (I_{0} + \eps^{3} I_{3})
		\|g\|_{L^{2}_l \cap L^\infty_l}	\|h\|_{L^{2}_l \cap L^\infty_l} \|\rho\|_{L^2} \|f\|_{L^2}
		.
		\een
	\end{prop}
	
	In the next subsection, we will do the energy estimate in  the weighted Sobolev space $H^{N}_{l}$. In particular, we have to estimate the typical term: 
	$\lr{R(\pa^{\alpha_1}f,\pa^{\alpha_2}f, \pa^{\alpha_3}f) W_{l}, W_{l}\pa^\alpha f}$ for $\alpha_1+\alpha_2+\alpha_3=\alpha$.	
	It can divided into four cases: Case 1: $\alpha_1= \alpha$;
	Case 2: $\alpha_2= \alpha$; Case 3: $\alpha_3= \alpha$;  Case 4: $|\alpha| \geq 1$ and $|\alpha_1|, |\alpha_2|, |\alpha_3| \le |\alpha|-1$. 
	We will apply \eqref{rough-upper-R-weighted-rho-L-infty} to  Case 1 and 2 and apply \eqref{rough-upper-R-weighted} to  Case 3. For Case 4, to balance the regularity of $g, h$ and $\rho$, we need additional  upper bound of $R$ which can be stated as follows.

	\begin{prop} \label{R-rough} Let $l \geq 0$, then
		\ben \label{R-h1-h1-h1-l2}
		|\lr{R (g, h, \rho), W_{l}f}|	
		\leq C_{l} (I_{0} + \eps^{3}  I_{3})   \|g\|_{H^1_2}\|h\|_{H^{1}_{l}} 
		\|\rho\|_{H^{1}_{l}} \|f\|_{L^2}.
		\een
	\end{prop} 
	\begin{proof}
		We  recall that
		$\lr{\eps^{-3}R (g, h, \rho), W_{l}f} = \int B g_* h
		\left( \rho' + \rho'_* \right) 
		((W_{l}f)' - W_{l}f)
		\mathrm{d}V 
		.$
		Obviously, using \eqref{weight-formula-1} and $B^{\eps} \leq 2B^{\eps}_1 + 2B^{\eps}_3$, we may get that
		\ben \label{directly-into-B1}
		|	\lr{\eps^{-3}R (g, h, \rho), W_{l}f}|
		&\lesssim_{l}&  \int B^{\eps}_1  |g_* W_{l}h|
		\left( |(W_{l}\rho)'| + |(W_{l}\rho)'_*| \right) 
		(|f'| + |f|)
		\mathrm{d}V 		
		\\ \nonumber 
		&&+ \int B^{\eps}_3  |g_* W_{l}h|
		\left( |(W_{l}\rho)'| + |(W_{l}\rho)'_*| \right) 
		(|f'| + |f|)
		\mathrm{d}V 
		.		
		\een		
		
		By \eqref{general-k-functions-h-norm}, we first have
		\ben \label{Q3-line-bound-H1-H1-H1}
		\int B^{\eps}_3  |g_* W_{l}h|
		\left( |(W_{l}\rho)'| + |(W_{l}\rho)'_*| \right) 
		(|f'| + |f|)
		\mathrm{d}V 		 \lesssim 
		I_{3}\|g\|_{H^{1}} \|  h\|_{H^{1}_{l}} 
		\| \rho\|_{H^{1}_{l}} \|f\|_{L^2}.
		\een
		
		By H{\"o}lder's inequality and  \eqref{general-result-cutoff} with  $b=0$,
		the integral containing $B^{\eps}_1$ in \eqref{directly-into-B1} is bounded by 		
		\ben \nonumber 
		&&\left(	\int B^{\eps}_1  |g_*| |W_{l}h|^{4}
		\mathrm{d}V 	\right)^{1/4}
		\left(	\int B^{\eps}_1  |g_*| |(W_{l}\rho)'|^{4} 	\mathrm{d}V 	\right)^{1/4}
		\left(
		\int B^{\eps}_1  |g_*| (|f'| + |f|)^2
		\mathrm{d}V 	\right)^{1/2}
		\\ \nonumber 
		&&+ \left(	\int B^{\eps}_1   |W_{l}h|^{2} |(W_{l}\rho)'_*|^{2} 
		\mathrm{d}V 	\right)^{1/2}
		\left(
		\int B^{\eps}_1  |g_*|^2 (|f'| + |f|)^2
		\mathrm{d}V 	\right)^{1/2}
		\lesssim  \eps^{-3}I_{0} 
		\|g\|_{L^{2}_2} \|W_l h\|_{H^{3/4}} \|W_l \rho\|_{H^{3/4}}   \|f\|_{L^{2}},
		\een
		where we use the Sobolev embedding theorem. Patching together all the estimates, we get \eqref{R-h1-h1-h1-l2}.
	\end{proof}

	We now apply  
	Propositions \ref{Q-inner-product-weighted-L2-L2} and \ref{R-rough} to get the following energy estimate of $R$ in $H^{N}_{l}$ spaces.
	\begin{lem}\label{UpR1R2} Let $N, l \geq 2$. Then
		\beno
		\big|\sum_{|\alpha| \leq N} \lr{W_{l} \pa^\alpha R(f, f, f) , W_{l}\pa^\alpha f} \big| \lesssim_{N,l} (I_{0} + \eps^{3}  I_{3})  
		\|f\|_{H^{N}_{l}}^4.
		\eeno
	\end{lem}
	\begin{proof} It is reduced to the consideration of $\lr{R(\pa^{\alpha_1}f,\pa^{\alpha_2}f, \pa^{\alpha_3}f) W_{l}, W_{l}\pa^\alpha f}$ for $\alpha_1+\alpha_2+\alpha_3=\alpha$.	
		There are four cases: Case 1: $\alpha_1= \alpha$;
		Case 2: $\alpha_2= \alpha$; Case 3: $\alpha_3= \alpha$;  Case 4: $|\alpha| \geq 1$ and $|\alpha_1|, |\alpha_2|, |\alpha_3| \le |\alpha|-1$. 
		The  desired result follows by applying \eqref{rough-upper-R-weighted-rho-L-infty} to Case 1 and Case 2,  \eqref{rough-upper-R-weighted} to Case 3 and  \eqref{R-h1-h1-h1-l2} to Case 4.
	\end{proof}

	\subsection{Upper bounds from the non-cutoff perspective}   In this subsection, we will give the upper bounds of the operators $Q, R$ from the non-cutoff perspective. Roughly speaking, we expand the Taylor expansion up to the second order to kill the singular factor $\eps^{-3}$.  
	We start with the uniform-in-$\eps$ estimate of $Q_1$ in the following Proposition.
	
	\begin{prop}\label{UpQ1Sob} It holds that
		\ben \label{Q1-g-h-f-L2-L2-H2}
		|\lr{Q_1(g,h),f}|\lesssim I_{3} (\|g\|_{L^{1}} + \|g\|_{L^{2}})
		\|h\|_{L^2}}\|f\|_{H^{2}, 
		\\ \label{Q1-g-h-f-L2-H2-L2}
		|\lr{Q_1(g,h),f}|\lesssim I_{3} (\|g\|_{L^{1}} + \|g\|_{L^{2}})\|h\|_{H^2}}\|f\|_{L^{2}. 
		\een
	\end{prop}
	\begin{proof} It is easy to see $\lr{Q_{1}(g,h), f} = \int B^{\eps}_{1} g_* h (f' -f) \mathrm{d}V.$
		By \eqref{Taylor1} for $f'-f$, \eqref{Q1-g-h-f-L2-L2-H2} follows by applying  \eqref{cancell1}, \eqref{order-2-cancellation} and \eqref{minus-2} for the first order and applying \eqref{general-B1-h-f} for the second order.
		For \eqref{Q1-g-h-f-L2-H2-L2}, we use Cancellation Lemma to 
		transfer the regularity from $f$ to $h$ through
		\beno \lr{Q_{1}(g,h), f} = \int B^{\eps}_{1}g_*\big(h - h'\big)
		f'\mathrm{d}V +\int B^{\eps}_{1}g_*\big((h f)'- hf\big)\mathrm{d}V.
		\eeno
		The former term is dealt with by \eqref{Taylor2}, \eqref{cancell2}, and \eqref{general-B1-h-f}. The latter term is estimated by \eqref{Q1-cancealltion-g-h-f}. Then we conclude \eqref{Q1-g-h-f-L2-H2-L2}.
	\end{proof}

	We next derive upper bounds of $Q_2$.
	\begin{prop}\label{UpQ2Sob} Let $\delta>0, a \geq 0, b \geq 1, a+b=\f32+\delta$. Then
		\ben \label{g-h-total-all-on-f}
		|\lr{Q_2(g,h),f}|\lesssim_{\delta} I_{3} \|g\|_{L^{2}}\|h\|_{L^2}\|f\|_{H^{\f32+\delta}},
		\\
		\label{g-h-total-3over2}
		|\lr{Q_2(g,h),f}|\lesssim_{\delta} (I_{3} + I^{\prime}_{3}) \|g\|_{H^{a}}\|h\|_{H^b}\|f\|_{L^{2}},
		\\ \label{g-h-total-vartheta}
		|\lr{Q_2(g,h),W_{l}f}|\lesssim \eps^{\vartheta} (I_{3+\vartheta} + I^{\prime}_{3+\vartheta}) \|g\|_{H^{2}_{l}}\|h\|_{H^2_{l}}\|f\|_{L^{2}}.
		\een
	\end{prop}
	\begin{proof}  We begin with the estimate of \eqref{g-h-total-all-on-f}. 
		We first observe that
		\beno |\int B^{\eps}_{2}g_* h (f' - f)
		\mathrm{d}V| \leq 2 \left(  \int B^{\eps}_{1} g^{2}_* (f' - f)^{2}
		\mathrm{d}V \right)^{1/2}  \left(  \int B^{\eps}_{3} 
		h^{2}\mathrm{d}V \right)^{1/2}. 
		\eeno	
		By the 
		order-1 Taylor expansion \eqref{Taylor1-order-1},   \eqref{Q1-result-with-factor-a-minus3-type-2} and \eqref{2-2-minus-1-s1-s2} imply that
		\ben\label{2-2-difference}
		\int B^{\eps}_{1} g^{2}_* (f' - f)^{2}
		\mathrm{d}V  &\lesssim& I_{3}
		\int \int_0^1  B^{\eps}_{1}  |g_*|^{2} |(\na f)(\kappa(v))|^{2}  |v'-v|^{2}
		\mathrm{d}\kappa \mathrm{d}V 
		\lesssim_{\delta} I_{3} \|g\|_{H^a}^2 \|f\|_{H^{b}}^2.
		\een
		By \eqref{general-L1}, the integral containing $B^{\eps}_3$ is bounded by $ I_{3} \|h\|_{L^{2}}^2$. Then \eqref{g-h-total-all-on-f} follows.
		
		As for \eqref{g-h-total-3over2}, we first have
		$\lr{Q_{2}(g,h), f} = \int B^{\eps}_{2}g_*\big(h - h'\big)
		f'\mathrm{d}V + \int B^{\eps}_{2} g_*\big((h f)'- hf\big)\mathrm{d}V.$
		Following the above argument for \eqref{g-h-total-all-on-f},  
		the term involving $h - h'$ is bounded by $I_{3} \|g\|_{H^{a}}\|h\|_{H^b}\|f\|_{L^{2}}$. The term involving $(h f)'- hf$ is estimated by
		\eqref{Q2-cancealltion-g-h-f}. These lead to   \eqref{g-h-total-3over2}.
		
		We turn to the estimate of \eqref{g-h-total-vartheta}. We observe that
		$ \lr{Q_{2}(g,h), W_{l}f} = \int B^{\eps}_{2}g_*\big(h - h'\big)
		(W_{l}f)'\mathrm{d}V + \int B^{\eps}_{2} g_*\big((W_{l} h f)'- W_{l}hf\big)\mathrm{d}V.$ We only focus on the first term since
		\eqref{eps-vartheta-B2-cancellation} implies that 
		$|\int B^{\eps}_{2} g_*\big((W_{l} h f)'- W_{l}hf\big)\mathrm{d}V| \lesssim 
		\eps^{\vartheta} (I_{3+\vartheta} + I^{\prime}_{3+\vartheta}) \|g\|_{H^{2}}  \|h\| _{H^{2}_{l}} \|f\|_{L^{2}}.$ We split it into two cases.
		
		$\bullet$ If $\vartheta \geq 1/2$, by
		\eqref{Taylor2} and \eqref{cancell2},  \eqref{weight-formula-1} together with the Cauchy-Schwartz inequality will lead to that
		\beno
		&& |\int B^{\eps}_{2}g_*\big(h - h'\big)
		(W_{l}f)'\mathrm{d}V| \lesssim  
		\int |B^{\eps}_{2} (W_{l}g)_* W_{l}(\kappa(v))(\na^2 h)(\kappa(v)) f'| |v'-v|^2 
		\mathrm{d}V \mathrm{d}\kappa
		\\
		&\lesssim&  \left( \int B^{\eps}_{1} |(W_{l}g)_* W_{l}(\kappa(v))(\na^2 h)(\kappa(v))|^2 |v'-v|^4 |v-v_*|^{-\vartheta} 
		\mathrm{d}V \mathrm{d}\kappa \right)^{1/2}
		\left( \int B^{\eps}_{3}|f'|^2 |v-v_*|^{\vartheta} 
		\mathrm{d}V \right)^{1/2}
		.
		\eeno
		Using the change of variable \eqref{change-of-variable-2} with $\kappa =\iota =1$,   \eqref{L1-vartheta} yields that
		\beno
		\int B^{\eps}_{3}|f'|^2 |v-v_*|^{\vartheta} 
		\mathrm{d}V = \int B^{\eps}_{3} |f|^2 |v-v_*|^{\vartheta} 
		\mathrm{d}V = \| 
		(J_{0,\eps}  |\cdot|^{\vartheta}) * f^2  \|_{L^1} \lesssim 
		\eps^{\vartheta} I_{3+\vartheta} \|f\|_{L^{2}}^2.
		\eeno Then we can conclude \eqref{g-h-total-vartheta} since
		\eqref{Q1-result-with-factor-a-minus3-with-small-facfor} and  the Hardy's inequality yield that 
		\beno
		&& \int B^{\eps}_{1} |(W_{l}g)_* W_{l}(\kappa(v))(\na^2 h)(\kappa(v))|^2 |v'-v|^4 |v-v_*|^{-\vartheta} 
		\mathrm{d}V \mathrm{d}\kappa 
		\\ &&\lesssim 
		\eps^{\vartheta} I_{3+\vartheta} \int |v-v_*|^{1-2\vartheta}  |(W_{l}g)_* W_{l} \na^2 h|^2  \mathrm{d}v_* \mathrm{d}v \lesssim 
		\eps^{\vartheta} I_{3+\vartheta} \|g\|_{H^{\vartheta - 1/2}_{l}}^2
		\|h\|_{H^{2}_{l}}^2 \lesssim 
		\eps^{\vartheta} I_{3+\vartheta} \|g\|_{H^{1/2}_{l}}^2
		\|h\|_{H^{2}_{l}}^2
		.\eeno

		$\bullet$ If $\vartheta \leq 1/2$, Thanks to the Taylor expansion, 
		\eqref{weight-formula-1} will imply that\beno
		&& |\int B^{\eps}_{2}g_*\big(h - h'\big)
		(W_{l}f)'\mathrm{d}V| \lesssim  
		\int |B^{\eps}_{2} (W_{l}g)_* W_{l}(\kappa(v))(\na h)(\kappa(v)) f'| |v'-v|
		\mathrm{d}V \mathrm{d}\kappa
		\\
		&&\lesssim   \left( \int B^{\eps}_{1} |(W_{l}g)_* W_{l}(\kappa(v))(\na h)(\kappa(v))|^2 |v'-v|^2 \sin\f{\theta}{2} |v-v_*|^{-\vartheta} 
		\mathrm{d}V \mathrm{d}\kappa \right)^{1/2}
		\left( \int B^{\eps}_{3} \sin^{-1}\f{\theta}{2} |f'|^2 |v-v_*|^{\vartheta} 
		\mathrm{d}V \right)^{1/2}
		.
		\eeno
		We first use \eqref{change-of-variable-2} and \eqref{L1-vartheta} to derive that
		$\int B^{\eps}_{3}|f'|^2 |v-v_*|^{\vartheta} 
		\mathrm{d}V = \int B^{\eps}_{3} \sin^{-1}\f{\theta}{2} |f|^2 |v-v_*|^{\vartheta} 
		\mathrm{d}V = \| 
		(K_{\eps} |\cdot|^{\vartheta}) * f^2  \|_{L^1},$
		where
		$K_{\eps}(z) \colonequals   \int_{\SS^{2}_{+}} \eps^{-4}|z|  \hphi^2(\eps^{-1}|z| \cos(\theta/2)) \sin^{-1}\f{\theta}{2}
		\mathrm{d}\sigma$.
		It is easy to see that
		\beno
		\|K_{\eps}|\cdot|^{\vartheta}\|_{L^{1}} 
		\lesssim  \eps^{\vartheta} I_{3+\vartheta} \int_{\SS^{2}_{+}} \cos^{-4-\vartheta}(\theta/2) \sin^{-1}\f{\theta}{2}
		\mathrm{d}\sigma \lesssim  \eps^{\vartheta} I_{3+\vartheta},
		\eeno
		which yields  $\| 
		(K_{\eps} |\cdot|^{\vartheta}) * f^2  \|_{L^1} \lesssim 
		\eps^{\vartheta} I_{3+\vartheta} \|f\|_{L^{2}}^2.$
		Following the similar argument used in the before, 
		we have
		\beno
		&& \int B^{\eps}_{1} |(W_{l}g)_* W_{l}(\kappa(v))(\na h)(\kappa(v))|^2 |v'-v|^2 \sin\f{\theta}{2} |v-v_*|^{-\vartheta} 
		\mathrm{d}V \mathrm{d}\kappa 
		\lesssim 
		\eps^{\vartheta} I_{3+\vartheta} \|g\|_{H^{1}_{l}}^2
		\|h\|_{H^{2}_{l}}^2
		.
		\eeno		
		Patching together the above estimates, we arrive at  \eqref{g-h-total-vartheta}.		
	\end{proof}
	
	
	Thanks to Propositions \ref{UpQ1Sob}, \ref{UpQ2Sob} and \ref{Q-3}, we get the following 
	upper bounds of $Q$ uniformly in $\eps$.
	
	\begin{thm} 
		It holds that
		\ben \label{l2-h2-l2}
		|\lr{Q(g,h), f}| \lesssim I_{3}(\|g\|_{L^1} + 	\|g\|_{L^2}) \|h\|_{L^2}\|f\|_{H^2}	 ,
		\\ \label{l2-l2-h2}
		|\lr{Q(g,h), f}| \lesssim (I_{3} + I^{\prime}_{3})(\|g\|_{L^1} + 	\|g\|_{L^2})	 \|h\|_{H^2}\|f\|_{L^2}.
		\een
	\end{thm}
	\smallskip
	
	As a direct application, we get the following upper bounds of $R$ with the small factor $\eps^{3}$.
	\begin{thm} 
		It holds that
		\ben \label{R-general-rough-versionh-h2-h2-h2}
		|\lr{R (g, h, \rho), f}|	
		\lesssim \eps^{3} (I_{3} + I^{\prime}_{3}) \|g\|_{H^{2}_2} 	\|h\|_{H^{2}} 
		\|\rho\|_{H^{2}} \|f\|_{L^2}.
		\een
	\end{thm} 
	\begin{proof}
		We observe that 
		$|\lr{R (g, h, \rho), f}|= \eps^{3} |\lr{Q(g,h),\rho f}+\lr{Q(g,hf),\rho} +\lr{Q(g\rho,h),f}+\lr{Q(hf, \rho),g}+\lr{Q(hf, g),\rho} + \mathcal{I}_1 + \mathcal{I}_2|$, 
		where  $\mathcal{I}_1 = 
		\int B g_*(h-h')(\rho_*'-\rho_*)f' \mathrm{d}V$ and 
		$\mathcal{I}_2 = 
		\int B (hf)_*(g\rho-(g\rho)')\mathrm{d}V.$

		For the terms containing $Q$, we appropriately use \eqref{l2-h2-l2} and \eqref{l2-l2-h2} to get that 
		\beno
		&& |\lr{Q(g,h),\rho f}|+|\lr{Q(g,hf),\rho}| +|\lr{Q(g\rho,h),f}|+|\lr{Q(hf, \rho),g}|+|\lr{Q(hf, g),\rho}|
		\\ &\lesssim&
		(I_{3} + I^{\prime}_{3}) \|g\|_{H^{2}_2} 	\|h\|_{H^{2}} 
		\|\rho\|_{H^{2}} \|f\|_{L^2}.
		\eeno
		We first separate $\mathcal{I}_1$  into two parts:
		\beno
		|\mathcal{I}_1| \leq 
		2 \int B^{\eps}_1 |g_*(h-h')(\rho_*'-\rho_*)f'| \mathrm{d}V + 2 \int B^{\eps}_3 |g_*(h-h')(\rho_*'-\rho_*)f'| \mathrm{d}V
		.
		\eeno
		By \eqref{general-k-functions-h-norm}, the integral containing $B^{\eps}_3$ is bounded by $I_{3} \|g\|_{H^{1}} 	\|h\|_{H^{1}} 
		\|\rho\|_{H^{1}} \|f\|_{L^2}$. Thanks to \eqref{2-2-difference},  the integral containing $B^{\eps}_1$ is bounded by 
		\beno
		\left(	\int B^{\eps}_1 g_*^2 (h-h')^2 \mathrm{d}V \right)^{1/2}
		\left(	\int B^{\eps}_1 (\rho'-\rho)^2 f^2_* \mathrm{d}V \right)^{1/2} \lesssim I_{3}
		\|g\|_{L^{2}} 	\|h\|_{H^{2}} 
		\|\rho\|_{H^{2}} \|f\|_{L^2}.
		\eeno
		
		As for $\mathcal{I}_2$, we separate it into three parts according to the fact that $B = B^{\eps}_1 + B^{\eps}_2 +B^{\eps}_3$. Then applying \eqref{Q1-cancealltion-g-h-f}, \eqref{Q2-cancealltion-g-h-f} and \eqref{general-k-functions-h-norm} to  the integrals containing $
		B^{\eps}_1, B^{\eps}_2$ and $B^{\eps}_3$ respectively, we arrive at
		\beno
		|\mathcal{I}_2|  \lesssim (I_{3} + I^{\prime}_{3}) \|g\|_{H^{2}} 	\|h\|_{H^{2}} 
		\|\rho\|_{H^{2}} \|f\|_{L^2}.
		\eeno
		
		We complete the proof by patching together the above estimates.
	\end{proof}

	In the forthcoming Lemma,
	we  give another estimate of the commutator between the operator $R$ and the weight function $W_l$. Comparing to
	\eqref{rough-commutator-R-rho-L-infty} and \eqref{rough-commutator-R}, here we can get rid of $I_{0}$ and keep the small factor $\eps^{3}  I_{3}$
	by imposing more regularity on the involved functions.

	\begin{lem} 
		Let $l \geq 2$. Then
		\ben \label{R-commutator-rough-version}
		|\lr{R(g, h, \rho) W_{l}-R(g, W_{l}h, \rho), f}\big| \leq \eps^{3} C_{l} I_{3} \|g\|_{H^2_l} \|h\|_{H^2_l} \|\rho\|_{H^2} \|f\|_{L^2}.
		\een
	\end{lem} 
	\begin{proof} Recalling \eqref{Commutator-formula-R},
		it is easy to check that 
		\beno |\mathcal{K}| &=& \eps^{3} |\int  B^{\eps} g_*h (\rho' + \rho'_*) f'( W_{l}' -  W_{l})\mathrm{d}V| =  \eps^{3} |\int  B g'_* h' (\rho + \rho_*) f( W_{l}' -  W_{l})\mathrm{d}V|
		\\ &\leq& \eps^{3} |\int  B^{\eps} g_* h (\rho + \rho_*) f( W_{l}' -  W_{l})\mathrm{d}V|
		+ \eps^{3} |\int  B^{\eps} (g'_*-g_*) h' (\rho + \rho_*) f ( W_{l}' -  W_{l})\mathrm{d}V|
		\\ &&+ \eps^{3} |\int  B^{\eps} g_* (h'-h) (\rho + \rho_*) f( W_{l}' -  W_{l})\mathrm{d}V|
		\colonequals  \eps^{3} \mathcal{K}_1 + \eps^{3} \mathcal{K}_2 + \eps^{3} \mathcal{K}_3
		.\eeno
		For $\mathcal{K}_1$,  by \eqref{Taylor1}, \eqref{cancell1} and
		\eqref{order-2-cancellation},   \eqref{minus-2} and \eqref{minus-1} imply that
		$|\mathcal{K}_1| \lesssim I_{3}\|g\|_{H^1_l} \|h\|_{L^2_l} \|\rho\|_{L^\infty} \|f\|_{L^2}.$
		For $\mathcal{K}_2$, we use $\mathcal{K}_2 \leq 2 \mathcal{K}_2^1 + 2 \mathcal{K}_2^3$ where  
		$\mathcal{K}_2^1 \colonequals 	\int  B^{\eps}_1 |(g'_*-g_*) h' (\rho + \rho_*) f ( W_{l}' -  W_{l})| \mathrm{d}V$ and 
		$\mathcal{K}_2^3 \colonequals 	\int  B^{\eps}_3 |(g'_*-g_*) h' (\rho + \rho_*) f ( W_{l}' -  W_{l})|\mathrm{d}V$
		.
		By order-1 Taylor expansion \eqref{Taylor1-order-1},   \eqref{general-B1-h-f} yields that
		\beno
		|\mathcal{K}_2^1| \lesssim \|\rho\|_{L^\infty}	\int  B^{\eps}_1 | (\na g)(\iota(v_*)) h'  f | W_{l-1}(\iota(v_*)) W_{l-1}(v')| |v'-v|^2 \mathrm{d}V
		\lesssim I_{3} \|g\|_{H^1_l} \|h\|_{L^2_l} \|\rho\|_{L^\infty} \|f\|_{L^2}.
		\eeno
		Using \eqref{weight-formula-1} and
		\eqref{general-k-functions-infty-norm}, 
		we have
		$|\mathcal{K}_2^3| \lesssim I_{3} \|g\|_{L^\infty_l} \|h\|_{L^2_l} \|\rho\|_{L^\infty} \|f\|_{L^2}.$
		As a result, by Sobolev embedding theorem, we get that
		$|\mathcal{K}_2| \lesssim I_{3} \|g\|_{H^2_l} \|h\|_{L^2_l} \|\rho\|_{H^2} \|f\|_{L^2}.$
		Similarly, we can also  get $|\mathcal{K}_3| \lesssim I_{3} \|g\|_{L^2_l} \|h\|_{H^2_l} \|\rho\|_{H^2} \|f\|_{L^2}.$ These are enough to conclude \eqref{R-commutator-rough-version}.
	\end{proof}

	By the upper bound estimate \eqref{R-general-rough-versionh-h2-h2-h2} and the commutator estimate \eqref{R-commutator-rough-version}, we have the following
	upper bound estimate for $R$ in weighted Sobolev spaces. 
	
	\begin{prop} \label{weighted-R-rough}
		Let $l \geq 2$, then
		\ben \label{R-general-rough-versionh-h2-h2-h2-weighted}
		|\lr{R (g, h, \rho), W_l f}|	
		\leq \eps^{3} C_{l} (I_{3} + I^{\prime}_{3}) \|g\|_{H^{2}_l} 	\|h\|_{H^{2}_l} 
		\|\rho\|_{H^{2}} \|f\|_{L^2}.
		\een
	\end{prop}
	
	Note that the small factor $\eps^{3}$ is kept in \eqref{R-general-rough-versionh-h2-h2-h2-weighted} which  
	shows that $R$ is smaller term if the involved functions are regular enough. Proposition \ref{weighted-R-rough} 
	will be used in the last section to derive the asymptotic formula in Theorem \ref{mainthm}.

		\setcounter{equation}{0}
		\section{Uniform upper bounds in weighted Sobolev space} 
		In this section, we will close energy estimate for the Uehling-Uhlenbeck operator
		$Q_{UU}^{\eps}$ uniformly in $\eps$
		in the weighted Sobolev space $H^{N}_{l}$. 	More precisely, 
		we will derive 
		\begin{thm}  \label{UU-in-H-N-l-uniform-in-eps} 
			Let $N, l \ge 2$. Let $f \geq 0$, then
			\ben \label{energy-in-H-N-l-uniform-in-eps}
			\sum_{|\alpha| \leq  N}  \lr{W_{l} \pa^{\alpha} Q_{UU}^{\eps}(f), W_{l} \pa^{\alpha}f}
			\lesssim_{N,l} (I_{0} + I_{3} + I^{\prime}_{3})  \|f\|_{H^{N}_{l}}^2
			(\|f\|_{H^{N}_{l}} + \|f\|_{H^{N}_{l}}^2).
			\een
		\end{thm} 
		\begin{proof}
			Recalling \eqref{OPUU} and \eqref{OPQ}, $Q_{UU}^{\eps}(f) = Q(f,f) + R(f,f,f)$ and $Q(f,f) = Q_1(f,f) + Q_2(f,f) + Q_3(f,f)$, we write
			\beno
			\pa^{\alpha} Q_{UU}^{\eps}(f) = \pa^{\alpha} Q(f,f) + \pa^{\alpha} R(f,f,f),
			\eeno
			and
			\ben \label{highest-order}
			\pa^{\alpha} Q(f,f) &=&   Q(f, \pa^{\alpha}f) 
			\\ \label{pul-order} &&+
			\sum_{\alpha_1+\alpha_2 = \alpha, |\alpha_1| =  1} C^{\alpha_1}_{\alpha}
			Q_1(\pa^{\alpha_1}f, \pa^{\alpha_2}f) 
			\\ \label{less-than-minus2-Q1}	 &&+
			\sum_{\alpha_1+\alpha_2 = \alpha, |\alpha_1| \ge  2} C^{\alpha_1}_{\alpha}
			Q_1(\pa^{\alpha_1}f, \pa^{\alpha_2}f) 
			\\ \label{less-than-minus1-Q2} &&+
			\sum_{\alpha_1+\alpha_2 = \alpha, |\alpha_1| \ge  1} C^{\alpha_1}_{\alpha}
			Q_2(\pa^{\alpha_1}f, \pa^{\alpha_2}f) 
			\\ \label{less-than-minus1-Q3}  &&+
			\sum_{\alpha_1+\alpha_2 = \alpha, |\alpha_1| \ge  1} C^{\alpha_1}_{\alpha}
			Q_3(\pa^{\alpha_1}f, \pa^{\alpha_2}f) 
			.
			\een
			
			Let us give the estimates term by term. We first observe that the term $Q(f, \pa^{\alpha}f)$ in \eqref{highest-order} can be deal by the forthcoming coercivity estimate \eqref{H-2-result}. More precisely, using $f \geq 0$,
			\ben \label{coercivity-application}
			\lr{Q(f, \pa^{\alpha}f) W_l, W_l \pa^{\alpha}f } \leq C_l (I_{3} + I^{\prime}_{3}) \|f\|_{H^2_l}\|\pa^{\alpha}f\|_{L^2_l}^2.
			\een
			The term in \eqref{pul-order} will be given by Lemma \ref{Q1EN-g-order-1}. This is referred to as the penultimate order term. The terms in \eqref{less-than-minus2-Q1} and in \eqref{less-than-minus1-Q2}  will be done in Lemma \ref{Q1EN} and   	
			in Lemma \ref{Q2EN} respectively.  Then we conclude \eqref{energy-in-H-N-l-uniform-in-eps} since the terms $Q_3(\pa^{\alpha_1}f, \pa^{\alpha_2}f)$ in \eqref{less-than-minus1-Q3} and    $\pa^{\alpha} R(f,f,f)$ have been estimated  by Lemma \ref{Q3EN} and  Lemma \ref{UpR1R2}. 
		\end{proof}
		
		In the following theorem, we derive two types of coercivity estimates of $Q$. 
		
		\begin{thm}[Coercivity estimate of $Q$]
			\label{Lowerbounds}  Let $l\ge2$. Let $g \geq 0$. Then  
			\ben 
			\label{L-infty-result}
			-\lr{Q(g,f)W_l,fW_l} &\ge& \f18
			\int (B^{\eps}+B^{\eps}_{1})g_*((fW_l)'-fW_l)^2\mathrm{d}V
			\\ \nonumber &&- C_l (I_{3} + I^{\prime}_{3}) \|g\|_{L^1 \cap L^\infty}\|f\|_{L^2_{l}}^2
			-C_l (I_{3} + I^{\prime}_{3}) \|g\|_{L^2_{l}}\|f\|_{L^1 \cap L^\infty}\|f\|_{L^2_{l}},
			\\ 
			\label{H-2-result}
			-\lr{Q(g,f)W_l,fW_l} &\ge& \f18
			\int (B^{\eps}+B^{\eps}_{1})g_*((fW_l)'-fW_l)^2\mathrm{d}V
			-C_l (I_{3} + I^{\prime}_{3}) \|g\|_{H^2_l}\|f\|_{L^2_l}^2.
			\een
		\end{thm} 
		\begin{proof} We observe that
			$-\lr{Q(g,f)W_l,fW_l}=-\lr{Q(g,fW_l),fW_l}-\lr{Q(g,f)W_l-Q(g,fW_l),fW_l}$ 
			and \\
			$-\lr{Q(g,fW_l),fW_l}=\f12 \int B^{\eps}g_*((fW_l)'-fW_l)^2\mathrm{d}V
			+\f12 \int B^{\eps}g_*((f^2W_l^2)-(f^2W_l^2)')\mathrm{d}V
			.$
			We first notice that by \eqref{commutator-Q-for-highest-order-l-infty}  the commutator can be estimate as follows:  
			for $0<\eta<1$,
			\beno
			&&|\lr{Q(g,f)W_l-Q(g,fW_l),fW_l}| 
			\\ &\lesssim_{l}& \eta \int B^{\eps}_1 g_* ((fW_l)'-fW_l)^2 \mathrm{d}V + \eta^{-1} I_{3} \|g\|_{L^1 \cap L^\infty}\|f\|_{L^2_{l}}^2 +  I_{3} \|g\|_{L^2_{l}}\|f\|_{L^1 \cap L^\infty}\|f\|_{L^2_{l}}
			.
			\eeno
			Then
			Lemma \ref{cancellem}  and \eqref{general-k-functions-infty-norm} imply that
			\beno 
			\big|\int B^{\eps}g_*((f^2W_l^2)-(f^2W_l^2)')\mathrm{d}V\big| \lesssim  (I_{3} + I^{\prime}_{3}) \|g\|_{L^\infty}\|f\|_{L^2_l}^2
			.
			\eeno
			
			By taking $\eta$ small enough,
			we  derive that
			\beno -\lr{Q(g,f)W_l,fW_l} &\ge& \f38
			\bigg(\int B^{\eps}g_*((fW_l)'-fW_l)^2\mathrm{d}V\bigg) \\&&- C_l (I_{3} + I^{\prime}_{3}) \|g\|_{L^1 \cap L^\infty}\|f\|_{L^2_{l}}^2
			-C_l (I_{3} + I^{\prime}_{3}) \|g\|_{L^2_{l}}\|f\|_{L^1 \cap L^\infty}\|f\|_{L^2_{l}}
			.\eeno
			Using the fact $B^{\eps} \ge \f12 B^{\eps}_1 - B^{\eps}_3$, we get
			\beno \int B^{\eps}g_*((fW_l)'-fW_l)^2\mathrm{d}V\ge \f12\int B^{\eps}_1g_*((fW_l)'-fW_l)^2\mathrm{d}V-\int B^{\eps}_3g_*((fW_l)'-fW_l)^2\mathrm{d}V.\eeno
			By \eqref{general-k-functions-infty-norm}, we have
			$\int B^{\eps}_3g_*((fW_l)'-fW_l)^2\mathrm{d}V \lesssim I_{3}  \|g\|_{L^\infty}\|f\|_{L^2_l}^2$
			. Then
			\eqref{L-infty-result} follows. 
			If  
			\eqref{commutator-Q-for-highest-order} is applied to $\lr{Q(g,f)W_l-Q(g,fW_l),fW_l}$, we will get \eqref{H-2-result}.
		\end{proof}
		\begin{rmk} \label{Q1-still-holds}
			Noting that $B^{\eps}_1 \ge \f12 B^{\eps}-B^{\eps}_3$,
			the estimates in Theorem \ref{Lowerbounds}  still hold if $Q$ is replaced by $Q_1$.
		\end{rmk}

		\subsection{Commutator estimates and weighted upper bounds of $Q$ from the angular non-cutoff perspective} 
		We have the following 
		commutator estimates that are uniform in $\eps$.

		\begin{lem}\label{COMT2W} For $i=1,2,3$, we define
			$\mathcal{I}_i\colonequals \lr{Q_i(g,h)W_l-Q_i(g,hW_l),f}$   
			$\mathcal{I} \colonequals \lr{Q(g,h)W_l-Q(g,hW_l),f}.$
			\begin{enumerate}\item 
				If $g\ge0$, then for $\eta>0$, then
				\ben \label{commutator-Q1-for-highest-order}
				&&|\mathcal{I}_1|  \lesssim_{l}  \eta \int B^{\eps}_1g_*(f'-f)^2\mathrm{d}V + \eta^{-1} I_{3} \|g\|_{L^1 \cap L^2}\|h\|_{L^2_{l-1}}^2
				+ I_{3}(\|g\|_{L^2_l}\|h\|_{L^2}+ \|g\|_{L^{1} \cap H^{1}}\|h\|_{L^2_{l-1}}) \|f\|_{L^2},
				\\ \label{commutator-Q1-for-highest-order-l-infty}
				&&|\mathcal{I}_1|   \lesssim_{l}  \eta \int B^{\eps}_1g_*(f'-f)^2\mathrm{d}V + \eta^{-1} I_{3} \|g\|_{L^1 \cap L^2}\|h\|_{L^2_{l-1}}^2 
				+ I_{3} (\|g\|_{L^2_l}\|h\|_{L^2}+\|g\|_{L^{1} \cap L^{\infty}}\|h\|_{L^2_{l-1}}) \|f\|_{L^2}.
				\een
				In general, it holds that
				\ben \label{commu-Q1-with-Wl}
				|\mathcal{I}_1|
				\lesssim_{l}  I_{3} (\|g\|_{L^2_2}\|h\|_{H^1_{l-1}}+\|g\|_{L^2_l}\|h\|_{H^1_1} )\|f\|_{L^2}.
				\een
				\item  For $(a,b)=(1,0)$ or $(0,1)$,  then 
				\ben \label{commu-Q2-with-Wl}
				|\mathcal{I}_2| \lesssim_{l} I_{3}  (\|g\|_{H^a}\|h\|_{H^b_{l-1}}+ \|g\|_{H^a_{l-1}}\|h\|_{H^b} )\|f\|_{L^2},
				\\ \label{commu-Q2-with-Wl-l-infty}
				|\mathcal{I}_2| \lesssim_{l} I_{3}   (\|g\|_{L^1 \cap L^\infty}\|h\|_{L^2_{l-1}} + \|g\|_{L^2_{l-1}}\|h\|_{L^1 \cap L^\infty})\|f\|_{L^2}.
				\een
				\item   For $\delta>0$ and $c,d \geq 0$ with $c+d = \f32 + \delta$, then
				\ben  \label{commu-Q3-with-Wl}
				|I_3| \lesssim_{l, \delta} I_{3}   (\|g\|_{H^c}\|h\|_{H^d_{l}} + \|g\|_{H^c_l}\|h\|_{H^d})\|f\|_{L^2},
				\\
				\label{commu-Q3-with-Wl-L-infty}
				|I_3| \lesssim_{l, \delta} I_{3}   (\|g\|_{L^\infty}\|h\|_{L^2_l} + \|g\|_{L^2_l}\|h\|_{L^\infty})\|f\|_{L^2}.
				\een

				\item  If $g\ge0$, then for $\eta>0$,  
				\ben \label{commutator-Q-for-highest-order}
				|\mathcal{I}| \lesssim_{l} \eta \int B^{\eps}g_*(f'-f)^2\mathrm{d}V 
				+ \eta^{-1} I_{3}   \|g\|_{L^2_2}\|h\|_{L^2_{l-1}}^2 + I_{3}   \|g\|_{H^2_l}\|h\|_{L^2_l} \|f\|_{L^2}.
				\een
				In general,
				\ben \label{commutator-Q-for-highest-order-l-infty}
				|\mathcal{I}| &\lesssim_{l}& \eta \int B^{\eps}_1g_*(f'-f)^2\mathrm{d}V + \eta^{-1} I_{3}   \|g\|_{L^1 \cap L^2}\|h\|_{L^2_{l-1}}^2 \\ \nonumber &&+ I_{3} (\|g\|_{L^1 \cap L^\infty}\|h\|_{L^2_{l}} + I_{3}   \|g\|_{L^2_{l}}\|h\|_{L^1 \cap L^\infty})\|f\|_{L^2}.
				\een
			\end{enumerate}
		\end{lem}
		\begin{proof}
			It is easy to compute that \ben\label{CQ1}\mathcal{I}_i=\int B^{\eps}_ig_*hf'(W_l'-W_l)\mathrm{d}V.\een

			{\bf \underline{\it Step 1: Estimate of $\mathcal{I}_1$.}}
			We first consider the case that $g\ge0$. By the Taylor expansion \eqref{Taylor1} for $W_l'-W_l$, we have  $\mathcal{I}_1= \mathcal{I}_1^1 + \mathcal{I}_1^2$ where
			$ \mathcal{I}_1^1  \colonequals  \int B^{\eps}_1g_*h f' (\na W_l)(v)\cdot (v'-v)\mathrm{d}V$ and 
			$\mathcal{I}_1^2  \colonequals 	\int B^{\eps}_1g_*h f' (1-\kappa) (\na^2 W_{l})(\kappa(v)):(v'-v)\otimes(v'-v)\mathrm{d}\kappa \mathrm{d}V.$
			For $\mathcal{I}_1^1$, 
			we have
			\beno  \mathcal{I}_1^1 = \int B^{\eps}_1g_*h(f'-f) (\na W_l)(v)\cdot (v'-v)\mathrm{d}V+\int B^{\eps}_1g_*hf (\na W_l)(v)\cdot (v'-v)\mathrm{d}V
			\colonequals \mathcal{I}_1^{1,1}+\mathcal{I}_1^{1,2}
			.\eeno
			\eqref{general-B1-h-f} implies that
			\beno |\mathcal{I}_1^{1,1}|&\lesssim& \bigg(\int B^{\eps}_1g_*(f'-f)^2\mathrm{d}V\bigg)^{\f12}
			\bigg( \int B^{\eps}_1g_*|hW_{l-1}|^2|v'-v|^2\mathrm{d}V\bigg)^{\f12}\\
			&\lesssim& \eta \int B^{\eps}_1g_*(f'-f)^2\mathrm{d}V + \eta^{-1} I_{3} \|g\|_{L^{1} \cap L^2}\|h\|_{L^2_{l-1}}^2.\eeno
			Using \eqref{cancell1},   \eqref{order-2-cancellation} and \eqref{minus-2},
			we have
			\beno 
			|\mathcal{I}_1^{1,2}|&=&\big|\int  B^{\eps}_1g_*hf (\na W_l)(v)\cdot (v-v_*) \sin^2(\theta/2)
			\mathrm{d}V\big| \\ 
			&\lesssim&  I_{3} \int |v-v_*|^{-2} g_*|hW_{l-1}|f\mathrm{d}v_* \mathrm{d}v
			\lesssim I_{3} \|g\|_{H^1}\|h\|_{L^2_{l-1}}\|f\|_{L^2} \text{ or }
			I_{3} \|g\|_{L^2 \cap L^{\infty}} \|h\|_{L^2_{l-1}}\|f\|_{L^2} 
			.\eeno
			For $\mathcal{I}_1^2$, by applying  \eqref{general-B1-h-f}, one has 
			\ben\qquad \label{I-1-2-upper-bound}
			|\mathcal{I}_1^2| \lesssim \int B^{\eps}_1 |g_*||h||f'||v'-v|^2 (W_{l-2} + (W_{l-2})_* )  \mathrm{d}V \lesssim I_{3} (\|g\|_{L^2_l}\|h\|_{L^2}+\|g\|_{L^1 \cap L^2}\|h\|_{L^2_{l-2}})\|f\|_{L^2}.
			\een
			Patching together the above estimates, we get	\eqref{commutator-Q1-for-highest-order} and \eqref{commutator-Q1-for-highest-order-l-infty} by using $\|g\|_{L^2 } \lesssim \|g\|_{L^1} + \|g\|_{L^{\infty}}$.

			We next turn to the general case. Using \eqref{Taylor2} for $W_l'-W_l$, we have $\mathcal{I}_1 =\mathcal{I}_1^3+\mathcal{I}_1^4$ where
			$\mathcal{I}_1^3 \colonequals	\int B^{\eps}_1g_*h f' (\na W_{l})(v')\cdot(v'-v) \mathrm{d}V$ and 
			$\mathcal{I}_1^4 \colonequals	-\int B^{\eps}_1g_*h f' \kappa(\na^2 W_{l})(\kappa(v)):(v'-v)\otimes(v'-v)\mathrm{d}\kappa \mathrm{d}V.$
			For $\mathcal{I}_1^3$, thanks to \eqref{cancell2}, we observe that
			\beno |\mathcal{I}_1^3| = \big|\int B^{\eps}_1g_*(h-h')f'(\na W_l)(v')(v'-v)\mathrm{d}V\big|.\eeno
			Then by order-1 Taylor expansion for $h-h'$,  suitably using \eqref{general-B1-h-f}, we 
			get
			\beno |\mathcal{I}_1^3|
			&\lesssim&  \int B^{\eps}_1 |g_*||(\na h)(\kappa(v))||f'||v'-v|^2(W_{l-1}(v_*)+W_{l-1}(\kappa(v))) \mathrm{d}\kappa \mathrm{d}V 
			\\ &\lesssim&  I_{3} \|g\|_{L^2_{l}} \|h\|_{H^1_{1}} \|f\|_{L^2} + I_{3} \|g\|_{L^2_2}\|h\|_{H^1_{l-1}}\|f\|_{L^2}.
			\eeno
			By comparing $\mathcal{I}_1^2$ and $\mathcal{I}_1^4$, it is not difficult to see that  $\mathcal{I}_1^4$ is also bounded by \eqref{I-1-2-upper-bound}. We get
			\eqref{commu-Q1-with-Wl}.
			
			\smallskip
			
			{\bf \underline{\it Step 2: Estimate of $\mathcal{I}_2$.}} Thanks to the Taylor expansion \eqref{Taylor1} for $W_l'-W_l$, we deduce that
			\beno |\mathcal{I}_2|\lesssim \int |B^{\eps}_2| |g_*||h||f'| |v-v'|
			\big(W_{l-1}+(W_{l-1})_*\big) \mathrm{d}V. \eeno
			By  the estimates \eqref{Q1-result-with-factor-a-minus3-type-2} and \eqref{general-L1} yield that
			\beno
			\mathcal{I} &\lesssim&
			\bigg(\int |B^{\eps}_1| (|g_*|^{2}|W_{l-1} h|^2 + |(W_{l-1}g)_*|^{2}|h |^2)|v-v'|^2\mathrm{d}V\bigg)^{\f12}
			\bigg(\int |B^{\eps}_3| |f'|^2\mathrm{d}V\bigg)^{\f12}
			\\  &\lesssim& I_{3}
			\bigg(\int  (|g_*|^{2}|W_{l-1} h|^2 + |(W_{l-1}g)_*|^{2}|h |^2)|v-v_{*}|^{-1}
			\mathrm{d}v \mathrm{d}v_{*} \bigg)^{\f12} \|f\|_{L^2}.
			\eeno
			Now we can use \eqref{2-2-minus-1-s1-s2} to get \eqref{commu-Q2-with-Wl} and  use \eqref{minus-1} to get \eqref{commu-Q2-with-Wl-l-infty}.

			{\bf \underline{\it Step 3: Estimate of $I_3$.}} Since $W_l' + W_l \lesssim  W_l +	(W_l)_*$,  \eqref{general-k-functions-infty-norm} implies \eqref{commu-Q3-with-Wl-L-infty}. Similarly, using
			\eqref{general-k-functions-infty-norm} and	\eqref{general-k-functions-h-norm},  
			one can easily get \eqref{commu-Q3-with-Wl}.

			Obviously, \eqref{commutator-Q-for-highest-order} and \eqref{commutator-Q-for-highest-order-l-infty} follow
			\eqref{commutator-Q1-for-highest-order}, \eqref{commu-Q2-with-Wl}, \eqref{commu-Q3-with-Wl}  and
			\eqref{commutator-Q1-for-highest-order-l-infty}, \eqref{commu-Q2-with-Wl-l-infty}, \eqref{commu-Q3-with-Wl-L-infty} respectively.  
			%
		\end{proof}
		
		With the upper bounds in the previous section and the commutator estimates in Lemma \ref{COMT2W},
		we are ready to state the following weighted upper bounds for $Q_1$ and $Q_2$.
		
		\begin{prop} 
			Let $l \geq 2$, $(a,b)=(1,0)$ or $(0,1)$, then
			\ben \label{Q1-h1-h2-l2-weighted}
			|\lr{ Q_1(g, h), W_l f}|	
			\lesssim_{l} I_{3} \|g\|_{L^{2}_l} 	\|h\|_{H^{2}_l} 
			\|f\|_{L^2},
			\\
			\label{Q2-h1-h1-l2-weighted}
			|\lr{ Q_2(g, h), W_l f}|	
			\lesssim_{l} (I_{3} + I^{\prime}_{3}) \|g\|_{H^{a}_l} 	\|h\|_{H^{b+1}_l} 
			\|f\|_{L^2}.
			\een
		\end{prop}
		\begin{proof}  It is not difficult to conclude that \eqref{Q1-h1-h2-l2-weighted} and  \eqref{Q2-h1-h1-l2-weighted} follow \eqref{Q1-g-h-f-L2-H2-L2},\eqref{commu-Q1-with-Wl}   \eqref{g-h-total-3over2} and \eqref{commu-Q2-with-Wl}. 
		\end{proof}
		
		Now we are ready to state the weighted upper bound of the operator $Q_{UU}^{\eps}$.
		
		\begin{thm} 
			Let $N, l \geq 2$, then
			\ben \label{UU-general-rough-versionh-h2-h2-h2-weighted}
			\| Q_{UU}^{\eps}(f) \|_{L^{2}_{l}} \leq  C_{l} (I_{3} + I^{\prime}_{3}) \| f  \|_{H^{2}_{l}}^2
			+ \eps^{3} C_{l} (I_{3} + I^{\prime}_{3}) \| f \|_{H^{2}_{l}}^3, 
			\\ \label{lose-order-2-derivative}
			\| Q_{UU}^{\eps}(f) \|_{H^{N-2}_{l}}  \leq  C_{l,N} (I_{3} + I^{\prime}_{3}) \| f \|_{H^{N}_{l}}^2
			+ \eps^{3} C_{l,N} (I_{3} + I^{\prime}_{3}) \| f \|_{H^{N}_{l}}^3.
			\een
		\end{thm}
		\begin{proof}
			Recalling $Q_{UU}^{\eps}(f) = Q_1(f, f) + Q_2(f, f) + Q_3(f, f) + R(f, f, f)$,
			by \eqref{Q1-h1-h2-l2-weighted}, \eqref{Q2-h1-h1-l2-weighted}, \eqref{Q-3-upper-bound} and
			\eqref{R-general-rough-versionh-h2-h2-h2-weighted}, we derive
			\eqref{UU-general-rough-versionh-h2-h2-h2-weighted}. Of course, \eqref{lose-order-2-derivative} is a direct result of
			\eqref{UU-general-rough-versionh-h2-h2-h2-weighted}.
		\end{proof}

		%
		
		Note that \eqref{Q2-h1-h1-l2-weighted} allows us to consider $\lr{Q_2(\pa^{\alpha_1}g, \pa^{\alpha_2}f) W_{l}, W_{l}\pa^\alpha f}$ with $|\alpha_1|\ge 1$.

		\begin{lem}\label{Q2EN} Let $1 \le m = |\alpha| \le N$, then \beno
			\sum_{\alpha_1+\alpha_2 = \alpha, |\alpha_1|\ge1}|\lr{Q_2(\pa^{\alpha_1}g,\pa^{\alpha_2}f) W_{l}, W_{l}\pa^\alpha f }| \lesssim_{N,l} (I_{3} + I^{\prime}_{3}) \|g\|_{H^{N}_{l}} \|f\|_{H^{N}_{l}}^2.
			\eeno
		\end{lem}
		\begin{proof} 
			If $|\alpha_1| \ge 2$, then $|\alpha_2| \leq m -2$. By taking $s_1 = 0, s_2 = 2, s_3 =0$ in \eqref{Q-3-upper-bound}
			and $a=0, b=1$ in  \eqref{Q2-h1-h1-l2-weighted},  we get that
			\beno
			|\lr{Q_2(\pa^{\alpha_1}g,\pa^{\alpha_2}f) W_{l}, W_{l}\pa^\alpha f }|
			\lesssim_{l} (I_{3} + I^{\prime}_{3}) \|\pa^{\alpha_1}g\|_{L^{2}_l} 	\|\pa^{\alpha_2}f\|_{H^{2}_l} 
			\|\pa^\alpha f\|_{L^2_{l}} 	\lesssim_{l} (I_{3} + I^{\prime}_{3}) \|g\|_{H^{m}_l} 	\|f\|_{H^{m}_l}^2.
			\eeno
			If $|\alpha_1| = 1$, then $|\alpha_2| \leq m - 1$. By taking $s_1 = 1, s_2 = 1, s_3 =0$ in \eqref{Q-3-upper-bound}
			and $a=1, b=0$ in  \eqref{Q2-h1-h1-l2-weighted},  
			we get
			\beno
			|\lr{Q_2(\pa^{\alpha_1}g,\pa^{\alpha_2}f) W_{l}, W_{l}\pa^\alpha f }|
			\lesssim_{l} (I_{3} + I^{\prime}_{3}) \|\pa^{\alpha_1}g\|_{H^{1}_l} 	\|\pa^{\alpha_2}f\|_{H^{1}_l} 
			\|\pa^\alpha f\|_{L^2_{l}} 	\lesssim_{l} (I_{3} + I^{\prime}_{3}) \|g\|_{H^{2}_l} 	\|f\|_{H^{m}_l}^2.
			\eeno
			These are enough to conclude the desired result.
		\end{proof}

		Note that \eqref{Q1-h1-h2-l2-weighted} allows us to consider $\lr{Q_1(\pa^{\alpha_1}g, \pa^{\alpha_2}f) W_{l}, W_{l}\pa^\alpha f}$ with $|\alpha_1|\ge 2$.
		\begin{lem}\label{Q1EN} Let $ 2 \le m =|\alpha|\le N$, then \beno
			\sum_{\alpha_1+\alpha_2 = \alpha, |\alpha_1| \ge  2}
			|\lr{Q_1(\pa^{\alpha_1}g, \pa^{\alpha_2}f) W_{l}, W_{l}\pa^\alpha f}| \lesssim_{N,l} I_{3} \|g\|_{H^{m}_{l}} \|f\|_{H^{m}_{l}}^2.
			\eeno
		\end{lem}
		\begin{proof} Since $|\alpha_1|\ge 2$, $|\alpha_2| \leq m -2$, using \eqref{Q1-h1-h2-l2-weighted}, we have
			\beno
			|\lr{Q_1(\pa^{\alpha_1}g, \pa^{\alpha_2}f) W_{l}, W_{l}\pa^\alpha f}|
			\lesssim_{l} I_{3} \|\pa^{\alpha_1}g\|_{L^{2}_l} 	\|\pa^{\alpha_2}f\|_{H^{2}_l} 
			\|\pa^\alpha f\|_{L^2_{l}} \lesssim_{l} I_{3} \|g\|_{H^{m}_l} 	\|f\|_{H^{m}_l}^2.
			\eeno The desired result follows by summation.
		\end{proof}


		\subsection{The penultimate order terms} In this subsection, we deal with the penultimate order terms $ \lr{Q_1(\pa^{\alpha_1}g,\pa^{\alpha_2}f) W_{l}, W_{l}\pa^\alpha f}$ where $|\alpha_1| = 1$.
		The forthcoming Lemma is motivated by  \cite{chaturvedi2021stability}. For the ease of interested readers, we reproduce its proof in a clearer way. The proof is based on
		some basic property of Boltzmann-type integral and integration by parts. 
		\begin{lem}\label{integration-by-parts-formula} 
			For simplicity, let $\tilde{\pa}_{i} \colonequals \pa^{e_{i}} + \pa^{e_{i}}_*$ where $\pa^{e_{i}} = \pa_{v_{i}}, \pa^{e_{i}}_* = \pa_{(v_{i})_*}$ for a unit index $|e_{i}| =1$. For a general function depending on the  variables
			$g = g(v, v_*, v', v'_*)$, let $g' = g(v', v'_*, v, v_*)$.
			For a general kernel $B=B(|v-v_{*}|,\cos\theta)$, and a general function $G = G(v)$,
			it holds that
			\beno
			\int B g (G^{\prime} - G) \tilde{\pa}_{i} G \mathrm{d}V
			= 
			\f{1}{4} \int B \tilde{\pa}_{i} g (G^{\prime} - G)^{2} \mathrm{d}V+ \f{1}{2}
			\int B (g - g') (G^{\prime} - G) \tilde{\pa}_{i} G  \mathrm{d}V.
			\eeno
		\end{lem}
		\begin{proof} We first observe that 
			\beno
			\tilde{\pa}_{i} B = 0, \quad \tilde{\pa}_{i} G  = \pa^{e_{i}} G, \quad \tilde{\pa}_{i} G_*  = (\pa^{e_{i}} G)_*, \quad \tilde{\pa}_{i} G'  = (\pa^{e_{i}} G)', \quad \tilde{\pa}_{i} G'_*  = (\pa^{e_{i}} G)'_*.
			\eeno
			As a result, 
			when  the integration by parts is used, the derivative on $B$ disappears. Moreover, derivatives and taking values are commutative through $\pa^{e_{i}}$
			and $\tilde{\pa}_{i}$. Using these facts, we have
			\ben \label{key-prood-line1}
			B g (G^{\prime} - G) \tilde{\pa}_{i} G = - B \tilde{\pa}_{i} g (G^{\prime} - G)  G - B g ( \tilde{\pa}_{i} G^{\prime} - \tilde{\pa}_{i} G) G
			\\ \label{key-prood-line2}
			= B \tilde{\pa}_{i} g (G - G^{\prime})^2  + B \tilde{\pa}_{i} g (G - G^{\prime}) G^{\prime}  - B g ( \tilde{\pa}_{i} G^{\prime} - \tilde{\pa}_{i} G) G
			\\ \label{key-prood-line3}
			= B \tilde{\pa}_{i} g (G - G^{\prime})^2  - B g ( \tilde{\pa}_{i} G - \tilde{\pa}_{i} G^{\prime}) G^{\prime} - B  g (G - G^{\prime}) \tilde{\pa}_{i} G^{\prime}  + B g ( \tilde{\pa}_{i} G - \tilde{\pa}_{i} G^{\prime} ) G
			\\ \label{key-prood-line4}
			= B \tilde{\pa}_{i} g (G - G^{\prime})^2   - 2 B  g (G - G^{\prime}) \tilde{\pa}_{i} G^{\prime}  + B g \tilde{\pa}_{i} G  (G - G^{\prime})
			\\ \label{key-prood-line5}
			= B \tilde{\pa}_{i} g (G - G^{\prime})^2   - 2 B  g^{\prime} (G^{\prime} - G) \tilde{\pa}_{i} G  - B g   (G^{\prime} -G) \tilde{\pa}_{i} G
			\\ \label{key-prood-line6}
			= B \tilde{\pa}_{i} g (G - G^{\prime})^2   - 2 B (g^{\prime} -g ) (G^{\prime} - G) \tilde{\pa}_{i} G  - 3 B g   (G^{\prime} -G) \tilde{\pa}_{i} G.
			\een
			Here for the equality in \eqref{key-prood-line1}, we use the integration by parts formula to transfer $\tilde{\pa}_{i}$ to other functions. For the equality from \eqref{key-prood-line2} to \eqref{key-prood-line3}, we use the integration by parts formula to deal with the middle term in the line \eqref{key-prood-line2}. For the equality from \eqref{key-prood-line4} to \eqref{key-prood-line5}, we use the change of variable $(v, v_*, \sigma) \to (v', v'_*, \sigma')$  to deal with the last two terms in \eqref{key-prood-line4}. The other changes are easily to verify. We conclude the result by moving the last term of  \eqref{key-prood-line6} to the left of \eqref{key-prood-line1}.
		\end{proof}
		

		Now we are ready to estimate the penultimate order terms.
		\begin{lem}\label{Q1EN-g-order-1}  Let $1 \le m = |\alpha| \le N$, then 
			\beno
			\sum_{\alpha_1+\alpha_2= \alpha, |\alpha_1|=1 }|\lr{Q_1(\pa^{\alpha_1}g,\pa^{\alpha_2}f) W_{l}, W_{l}\pa^\alpha f}| \lesssim_{N,l} I_{3} \|g\|_{H^{2}_{l}} \|f\|_{H^{m}_{l}}^2.
			\eeno
		\end{lem}
		\begin{proof} 
			Let $\mathcal{A}=|\lr{Q_1(\pa^{\alpha_1}g,\pa^{\alpha_2}f) W_{l}, W_{l}\pa^\alpha f }|$.  It is easy to see that  $\mathcal{A}\le |\mathcal{A}_1|+ |\mathcal{A}_2|$ where $\mathcal{A}_1= \lr{Q_1(\pa^{\alpha_1}g,  W_{l}\pa^{\alpha_2}f), W_{l}\pa^\alpha f }$ and $\mathcal{A}_2=\lr{Q_1(\pa^{\alpha_1}g,\pa^{\alpha_2}f) W_{l}-Q_1(\pa^{\alpha_1}f, W_{l}\pa^{\alpha_2}f), W_{l}\pa^\alpha f }.$ 
			By \eqref{commu-Q1-with-Wl} in Lemma \ref{COMT2W}, as $|\alpha_2| = m-1$, we have
			\beno 
			|\mathcal{A}_2|  \lesssim_{l} I_{3} (\| W_{l} \pa^{\alpha_1}g\|_{L^{2}}\|\pa^{\alpha_2}f\|_{H^1_1} + \|\pa^{\alpha_1}g\|_{L^{2}_{2}}\| W_{l} \pa^{\alpha_2}f\|_{H^1_{-1}})\| W_{l}\pa^\alpha f\|_{L^{2}} \lesssim_{l} I_{3} \|g\|_{H^{2}_l} 	\|f\|_{H^{m}_l}^2.
			\eeno
			Since $|\alpha_1|=1$, we write $\alpha_1 = e_{i}$ for some $1 \leq i \leq 3$.  Let $\pa^{e_{i}} = \partial_{v_{i}}$ and $F = \pa^{\alpha_2}f = \pa^{\alpha-e_{i}}f$, then $\pa^\alpha f =\pa^{e_{i}} F$ and
			\ben \label{A-1-into-2-terms} 
			\mathcal{A}_1  = 
			\lr{Q_1(\pa^{e_{i}}g,  W_{l} F), \pa^{e_{i}} ( W_{l} F)} + \lr{Q_1(\pa^{e_{i}}g,  W_{l} F),
				W_{l} \pa^{e_{i}}F -\pa^{e_{i}} ( W_{l} F)}
			.
			\een
			
			$\bullet$ Let $G =  W_{l} F$, then the first term of the r.h.s of \eqref{A-1-into-2-terms} can be written as
			\ben \label{first-trick}
			\lr{Q_1(\pa^{e_{i}}g, G), \pa^{e_{i}}G}
			= \int B^{\eps}_{1} (\pa^{e_{i}}g)^{\prime}_{*} (G^{\prime} - G) \pa^{e_{i}} G \mathrm{d}V
			+ \int B^{\eps}_{1} ((\pa^{e_{i}}g)^{\prime}_{*} - (\pa^{e_{i}}g)_{*}) G \pa^{e_{i}} G \mathrm{d}V
			.
			\een
			For the second term in \eqref{first-trick}, we use  \eqref{Q1-cancealltion-g-h-f} to get 
			\beno 
			|\int B^{\eps}_{1} ((\pa^{e_{i}}g)^{\prime}_{*} - (\pa^{e_{i}}g)_{*}) G \pa^{e_{i}} G \mathrm{d}V| 
			\lesssim I_{3} \|\pa^{e_{i}}g\|_{H^{1}} \|G\|_{H^{1}}\|\pa^{e_{i}}G\|_{L^{2}} \lesssim_{l} I_{3} \|g\|_{H^{2}} 	\|f\|_{H^{m}_l}^2.
			\eeno
			By  Lemma \ref{integration-by-parts-formula}, the first term in \eqref{first-trick} is
			\ben \label{transfer-regularity-trick}
			\f{1}{4} \int B^{\eps}_{1} (\pa^{2e_{i}}g)_{*} ( G^{\prime} - G)^{2} \mathrm{d}V+ \f{1}{2}
			\int B ((\pa^{e_{i}}g)_{*}^{\prime} - (\pa^{e_{i}}g)_{*}) (G^{\prime} - G) \pa^{e_{i}} G \mathrm{d}V.
			\een
			By Taylor expansion to $G$,   \eqref{general-B1-h-f} gives that
			\beno 
			\int B^{\eps}_{1} (\pa^{2e_{i}}g)_{*} (G^{\prime} - G)^{2} \mathrm{d}V  \lesssim  I_{3} \|g\|_{H^2_{2}} \|G\|_{H^{1}} \lesssim_{l} I_{3} \|g\|_{H^{2}_2} 	\|f\|_{H^{m}_l}^2.
			\eeno
			Apply  order-1 Taylor expansion to $\pa^{e_{i}}g$ and $G$ and use  \eqref{general-B1-h-f}, then the second term in \eqref{transfer-regularity-trick} is bounded by
			\beno |\int B^{\eps}_1  |(\na \pa^{e_{i}}g)(\iota(v_*))| |(\na G)(\kappa(v))|
			| \pa^{e_{i}} G| \mathrm{d}V \mathrm{d}\kappa \mathrm{d}\iota| 
			\lesssim_{l} I_{3}
			\|\na \pa^{e_{i}}g\|_{L^2_{2}} \|\na G\|_{L^{2}} \|G\|_{L^{2}} \lesssim I_{3} \|g\|_{H^{2}_2} 	\|f\|_{H^{m}_l}^2.
			\eeno

			$\bullet$  Since $ W_{l} \pa^{e_{i}}F -\pa^{e_{i}} ( W_{l} F) = - F \pa^{e_{i}}  W_{l} $, 
			the second term in the r.h.s of \eqref{A-1-into-2-terms} is
			\beno 
			- \lr{Q_1(\pa^{e_{i}}g,  W_{l} F),
				F \pa^{e_{i}}  W_{l}} &=& \f{1}{2} \int B_1 (\pa^{e_{i}}g)_{*} (( W_{l} F)^{\prime} -  W_{l} F) ((F \pa^{e_{i}}  W_{l})^{\prime} - F \pa^{e_{i}}  W_{l})\mathrm{d}V 
			\\&&+ \f{1}{2}
			\int B_1 ((\pa^{e_{i}}g)_{*}^{\prime}- (\pa^{e_{i}}g)_{*})  W_{l} F ((F \pa^{e_{i}}  W_{l})^{\prime} - F \pa^{e_{i}}  W_{l})\mathrm{d}V.
			\eeno
			By similar argument, the r.h.s can be bounded by 
			$I_{3} \|g\|_{H^{2}_2} 	\|f\|_{H^{m}_l}^2$ and $I_{3} \|g\|_{H^{2}_2} 	\|f\|_{H^{m}_l}^2$.  We are led to  the desired result.
		\end{proof}

		\setcounter{equation}{0}

		\section{Well-posedness and   propagation of regularity}
		
		In this section, we will prove 
		part of Theorem \ref{gwp-F-D} for Fermi-Dirac particles
		as well as Theorem \ref{lwp-B-E} for Bose-Einstein particles.

		\subsection{Fermi-Dirac particles} In this subsection, will prove the first two results in Theorem \ref{gwp-F-D}.
		We start with the mild solution since in this situation the non-negativity can be easily proved. We recall that
		$\pa_t f
		=Q_{UU}^{\eps}(f)$,
		where $Q_{UU}^{\eps}$   is defined by
		\ben\label{OPQUU-eps-FD} Q_{UU}^{\eps}(f) = \int_{\R^3\times\SS^2}
		B^{\eps} \Pi^{\eps}(f) \mathrm{d}\sigma \mathrm{d}v_*, \een
		where 
		\beno
		\Pi^{\eps}(f) \colonequals 
		f_*' f' (1 - \eps^3 f_*)  (1 - \eps^3 f) - f_* f (1 - \eps^3 f_*') (1 - \eps^3f').
		\eeno
		As usual,
		we define the gain and loss terms by
		\ben\label{OPQUU-eps-FD-gain-loss} Q^{\eps,+}(f) = \int_{\R^3\times\SS^2}
		B^{\eps} \Pi^{\eps,+}(f) \mathrm{d}\sigma \mathrm{d}v_*,
		\quad  Q^{\eps,-}(f) = \int_{\R^3\times\SS^2}
		B^{\eps} \Pi^{\eps,-}(f) \mathrm{d}\sigma \mathrm{d}v_*, \een
		where 
		\beno
		\Pi^{\eps,+}(f) \colonequals 
		f_*' f' (1 - \eps^3 f_*)  (1 - \eps^3 f), \quad 
		\Pi^{\eps,-}(f) \colonequals 
		f_* f (1 - \eps^3 f_*') (1 - \eps^3 f').
		\eeno

		We now prove that the initial value problem admits a mild solution if $f_0 \in L^1$ with $0 \leq f_{0} \leq \eps^{-3}$. To do that, we prove the following lemma.
		\begin{lem} 
			Let $f, g \in L^1 \cap L^{\infty}, 0 \leq f, g \leq \eps^{-3}$, then 
			\ben \label{Q-g-L1-F-D}
			\|Q_{UU}^{\eps}(f)\|_{L^1}      
			\lesssim \eps^{-3} I_0  \|f\|_{L^1}^2,
			\\ \label{Q-g-Linfty-F-D}
			\|Q_{UU}^{\eps}(f)\|_{L^\infty}  \lesssim \eps^{-6}  I_{0}
			\| f \|_{L^1} +  \eps^{-6} I_{3},
			\\ \label{Q-g-minus-f-L1-F-D}
			\|Q_{UU}^{\eps}(f) - Q_{UU}^{\eps}(g)\|_{L^1}  \lesssim (\eps^{-3}  I_{0} + I_{3}) (\| f \|_{L^1} + \| g \|_{L^1} +  \eps^{-3}) \| f -g \|_{L^1}.
			\een
		\end{lem}
		\begin{proof} 
			By \eqref{OPQUU-eps-FD} and \eqref{OPQUU-eps-FD-gain-loss}, we have
			\ben \label{Q-into-gain-and-loss}
			Q_{UU}^{\eps}(f) = Q^{\eps,+}(f) - Q^{\eps,-}(f).
			\een
			From \eqref{Q-into-gain-and-loss} and the change of variable \eqref{change-of-variable-2}
			in the special case $\kappa = \iota = 1$, we get that
			\ben \label{F-D-Q-L1-into-gain-and-loss}
			\|Q_{UU}^{\eps}(f)\|_{L^1} \leq \|Q^{\eps,+}(f)\|_{L^1} + \|Q^{\eps,-}(f)\|_{L^1} = 2 \|Q^{\eps,-}(f)\|_{L^1}
			.
			\een  
			Since $0 \leq g \leq \eps^{-3}$,  \eqref{upper-bound-of-A-eps} implies that
			\ben \label{F-D-Q-L1-loss-up}
			\|Q^{\eps,-}(f)\|_{L^1}
			=
			\int
			B^{\eps} 
			f_* f (1 - \eps^3 f_*') (1 - \eps^3 f') \mathrm{d}\sigma \mathrm{d}v_* \mathrm{d}v
			\leq \int
			B^{\eps} 
			f_* f \mathrm{d}\sigma \mathrm{d}v_* \mathrm{d}v
			\lesssim \eps^{-3} I_0  \|g\|_{L^1}^2.
			\een
			This gives  \eqref{Q-g-L1-F-D}.   
			
			From \eqref{Q-into-gain-and-loss}, 
			it holds that
			$\|Q_{UU}^{\eps}(f)\|_{L^\infty} \leq \|Q^{\eps,+}(f)\|_{L^\infty} + \|Q^{\eps,-}(f)\|_{L^\infty}.$ 
			Since $0 \leq f \leq \eps^{-3}$, by using \eqref{upper-bound-of-A-eps}, 
			we obtain that 
			\ben\qquad \label{Q-minus-L-infty}
			Q^{\eps,-}(f)(v) 
			= 
			\int
			B^{\eps} 
			f_* f (1 - \eps^3 f_*') (1 - \eps^3 f') \mathrm{d}\sigma \mathrm{d}v_* 
			\leq f \int
			B^{\eps} 
			f_* \mathrm{d}\sigma \mathrm{d}v_* = A^{\eps} f \|f\|_{L^1}
			\lesssim \eps^{-6} I_0  \|f\|_{L^1}.
			\een
			For $Q^{\eps,+}(f)(v)$, we use $B^{\eps} \leq 2 B^{\eps}_1 + 2 B^{\eps}_3$ to  get
			\ben \label{Q-plus-L-infty}
			Q^{\eps,+}(f)(v)  =
			\int
			B^{\eps} 
			f_*' f' (1 - \eps^3 f_*) (1 - \eps^3 f) \mathrm{d}\sigma \mathrm{d}v_* 
			\leq 2 \eps^{-3}  \int
			B^{\eps}_1 
			f_*'  \mathrm{d}\sigma \mathrm{d}v_*  + 2 \eps^{-6}  \int
			B^{\eps}_3 
			\mathrm{d}\sigma \mathrm{d}v_*.
			\een
			By \eqref{Q-plus-L-infty}, \eqref{bounded-l1-norm-simple} and \eqref{B1-sigma-vStar}, we arrive at
			\ben \label{Q-plus-L-infty-result}
			Q^{\eps,+}(f)(v) \lesssim \eps^{-6}  I_{0}
			\| f \|_{L^1} +  \eps^{-6} I_{3}
			.
			\een
			Patching together \eqref{Q-minus-L-infty} and \eqref{Q-plus-L-infty-result} will give \eqref{Q-g-Linfty-F-D}.
			
			From \eqref{Q-into-gain-and-loss},  
			we notice that
			$Q_{UU}^{\eps}(f) - Q_{UU}^{\eps}(g) = Q^{\eps,+}(f) - Q^{\eps,+}(g) - \left(Q^{\eps,-}(f) - Q^{\eps,-}(g)\right).$
			Recalling \eqref{OPQUU-eps-FD-gain-loss}, by the change of variable \eqref{change-of-variable-2}
			in the special case $\kappa = \iota = 1$, we get that
			\ben  \nonumber 
			&&\|Q_{UU}^{\eps}(f) - Q_{UU}^{\eps}(g)\|_{L^1}  \leq  \|Q^{\eps,+}(f) - Q^{\eps,+}(g)\|_{L^1} + \|Q^{\eps,-}(f) - Q^{\eps,-}(g)\|_{L^1} \\  
			&&\leq
			\int B^{\eps} |\Pi^{\eps,+}(f) - \Pi^{\eps,+}(g)| \mathrm{d}V + 	\int B^{\eps} |\Pi^{\eps,-}(f) - \Pi^{\eps,-}(g)| \mathrm{d}V
			=  2 \int B^{\eps} |\Pi^{\eps,-}(f) - \Pi^{\eps,-}(g)| \mathrm{d}V. \label{L1-upper-by-2-times-loss}
			\een
			
			Observe that
			\beno
			\Pi^{\eps,-}(f) - \Pi^{\eps,-}(g)  &=& (f_* - g_{*}) f (1 - \eps^3 f_*') (1 - \eps^3 f') + g_* (f - g) (1 - \eps^3 f_*') (1 - \eps^3 f') 
			\\ &&+ \eps^3 g_* g (g_*' - f_*') (1 - \eps^3 f') + \eps^3 g_* g  (1 - \eps^3 g_*') (g' - f')
			.
			\eeno
			If $ 0 \leq f, g \leq \eps^{-3}$, then
			$|\Pi^{\eps,-}(f) - \Pi^{\eps,-}(g)| \leq 
			|f_* - g_{*}| f + g_* |f - g| + 
			g |g_*' - f_*'|  +  g_* |g' - f'|
			.$
			Therefore, 
			\beno 
			\int B^{\eps} |\Pi^{\eps,-}(f) - \Pi^{\eps,-}(g)| \mathrm{d}V
			\leq 
			\int
			B^{\eps} (|f_* - g_{*}| f + g_* |f - g|) 
			\mathrm{d}V + 2 \int
			B^{\eps} g_* |g' - f'| 
			\mathrm{d}V
			.	\eeno
			By using \eqref{upper-bound-of-A-eps},  \eqref{general-L1} and \eqref{B1-sigma-vStar}, we get  that
			\beno
			\int B^{\eps} |\Pi^{\eps,-}(f) - \Pi^{\eps,-}(g)| \mathrm{d}V
			\lesssim 
			(\eps^{-3}  I_{0} + I_{3}) (\| f \|_{L^1} + \| g \|_{L^1} + \| g \|_{L^\infty})\| f -g \|_{L^1},
			\eeno
			which yields \eqref{Q-g-minus-f-L1-F-D}.
		\end{proof}

		\begin{proof}[Proof of Theorem \ref{gwp-F-D}: (Well-posedness part)]  We recall that  $\mathcal{A}_{T}  \colonequals L^{\infty}([0,T]; L^{1}(\R^3))$ associated with the norm $\|f\|_{T} \colonequals \sup\limits_{0 \leq t \leq T} \|f(t)\|_{L^{1}}$. We 
			define the operator $J^{\eps}(\cdot)$ on 
			$\mathcal{A}_{T}$ by 
			\beno
			J^{\eps}(f)(t,v) \colonequals f_{0}(v) + \int_{0}^{t} Q_{UU}^{\eps}(|f| \wedge \eps^{-3} )(\tau, v) \mathrm{d}\tau.
			\eeno
			
			$\bullet$ $J^{\eps}(\cdot)$ is a map from  $\mathcal{A}_{T}$ onto  $\mathcal{A}_{T}$. It is easy to check that  
			$\|J^{\eps}(f)\|_{T} \leq  \|f_{0}\|_{L^1} + T \| Q_{UU}^{\eps}(|f| \wedge \eps^{-3} )\|_{T}.$
			By \eqref{Q-g-L1-F-D},  
			we get 
			\ben \label{map-on-itself}
			\|J^{\eps}(f)\|_{T} \leq  \|f_{0}\|_{L^1} + T C_{\eps, \phi} 
			\|f\|_{T}^2,
			\een where   $C_{\eps, \phi} \lesssim \eps^{-3} 
			I_{0}$.
			This means  $J^{\eps}(\cdot)$ is an operator onto 
			$\mathcal{A}_{T}$.

			$\bullet$	Let $\mathcal{B}_{T}  \colonequals \{f \in  \mathcal{A}_{T},   \|f\|_{T} \leq 2  \|f_{0}\|_{L^1} \}$. We want to show that $J^{\eps}(\cdot)$ is contraction on the complete metric space $(\mathcal{B}_{T}, \|\cdot - \cdot\|_{T})$ for small enough $T>0$. By \eqref{map-on-itself}, 
			if $T$ satisfies
			\ben \label{small-condition-1-on-T}
			4 T C_{\eps, \phi} 
			\|f_{0}\|_{L^1} \leq 1,
			\een
			then
			$J^{\eps}(\cdot)$ is an operator onto  
			$\mathcal{B}_{T}$.  Given 
			$f, g \in \mathcal{B}_{T}$, we have 
			\beno
			J^{\eps}(f)(t,v) - J^{\eps}(g)(t,v) =   \int_{0}^{t} \left( Q_{UU}^{\eps}(|f| \wedge \eps^{-3}) -  Q_{UU}^{\eps}(|g| \wedge \eps^{-3}) 
			(\tau, v) \right) \mathrm{d}\tau.
			\eeno
			Similarly to the above, we   have
			$\|J^{\eps}(f) - J^{\eps}(g)\|_{T} \leq  T \| Q_{UU}^{\eps}(|f| \wedge \eps^{-3} ) - Q_{UU}^{\eps}(|g| \wedge \eps^{-3} )\|_{T}.$
			By \eqref{Q-g-minus-f-L1-F-D}, since the function $x \in \R \to |x| \wedge \eps^{-3}$ is Lipschitz continuous with Lipschitz constant
			$1$, 
			we get that
			\beno
			\|J^{\eps}(f) - J^{\eps}(g)\|_{T} \leq T C_{\epsilon, \phi} (\| f \|_{T} + \| g \|_{T} + \eps^{-3})
			\| f- g\|_{T}
			\leq T C_{\epsilon, \phi} (4 \| f_0 \|_{L^1} + \eps^{-3})
			\| f- g\|_{T}.
			\eeno
			where $ C_{\epsilon, \phi} \lesssim\eps^{-3}  I_{0} + I_{3}$.
			If $T$ satisfies
			\ben \label{small-condition-2-on-T}
			T C_{\epsilon, \phi} (4 \| f_0 \|_{L^1} + \eps^{-3}) \leq \f12,
			\een
			then $J^{\eps}(\cdot)$ is a contraction mapping.   By the Fixed Point Theorem,  
			there exists 
			a unique $f^{\eps} \in \mathcal{B}_{T}$, s.t.  $f^{\eps} = J^{\eps}(f^{\eps})$. After a modification on the $v$-null sets,  there is a null set $Z \subset \R^3$, for all $t \in  [0, T]$ and $v \in \R^3 \setminus Z$, 
			\beno 
			f^{\eps}(t,v) = J^{\eps}(f^{\eps}) (t,v) = f_{0}(v) + \int_{0}^{t} Q_{UU}^{\eps}(|f^{\eps}| \wedge \eps^{-3} )(\tau, v) \mathrm{d}\tau.
			\eeno

			$\bullet$ We now utilize the special definition of $J^{\eps}(\cdot)$
			to prove that $0 \leq f^{\eps}(t) \leq \eps^{-3}$ for all  $t \in  [0, T]$ and $v \in \R^3 \setminus Z$. 
			By \eqref{Q-g-Linfty-F-D},
			we can see that for any $v \in \R^3 \setminus Z$ and $t_1, t_2 \in [0, T]$,
			\beno
			|f^{\eps}(t_2,v) - f^{\eps}(t_1,v)| \leq C_{\eps, \phi} (\| f_0 \|_{L^1} + 1) |t_2 - t_1|,
			\eeno
			for some $C_{\eps, \phi} \lesssim  \eps^{-6} \left( I_{0} + I_{3} \right)$. That is, $f^{\eps}(\cdot, v)$ is uniformly continuous w.r.t. $t$ on $[0, T]$
			for any $v \in \R^3 \setminus Z$. Then we have
			\ben \nonumber 
			(-f^{\eps}(t,v))^{+} &=& - \int_{0}^{t} Q_{UU}^{\eps}(|f^{\eps}| \wedge \eps^{-3})(\tau, v)
			\mathrm{1}_{f^{\eps}(\tau,v)<0} \mathrm{d}\tau 
			\leq \int_{0}^{t} Q^{\eps,-}(|f^{\eps}| \wedge \eps^{-3})(\tau, v)
			\mathrm{1}_{f^{\eps}(\tau,v)<0} \mathrm{d}\tau
			\\ \label{non-negativity} 
			&\leq& A^{\eps} \int_{0}^{t} 
			\|f^{\eps}(\tau)\|_{L^1} |f^{\eps}(\tau, v)|
			\mathrm{1}_{f^{\eps}(\tau,v)<0} \mathrm{d}\tau \leq 2 A^{\eps} \|f_0\|_{L^1}
			\int_{0}^{t}  (-f^{\eps}(\tau,v))^{+}  \mathrm{d}\tau.	\een  By Gr\"{o}nwall's inequality, we get $(-f^{\eps}(t,v))^{+} = 0$ on $[0, T]$ and so $f^{\eps}(t,v) \geq 0$. We also have
			\beno
			(f^{\eps}(t,v) - \eps^{-3})^{+} &=& \int_{0}^{t} Q_{UU}^{\eps}(|f^{\eps}| \wedge \eps^{-3})(\tau, v)
			\mathrm{1}_{f^{\eps}(\tau,v)>\eps^{-3}} \mathrm{d}\tau 
			\le \int_{0}^{t} Q^{\eps,+}(|f^{\eps}| \wedge \eps^{-3})(\tau, v)
			\mathrm{1}_{f^{\eps}(\tau,v)>\eps^{-3}} \mathrm{d}\tau
			=0,
			\eeno
			which yields $f^{\eps}(t,v) \leq \eps^{-3}$.

			Recalling   Definition \ref{def-mild-solution} for F-D particles, we already get a mild solution on $[0, T]$. By \eqref{Q-g-L1-F-D} and 
			Fubini's Theorem, we have 
			conservation of mass $\| f^{\eps}(t) \|_{L^1} \colonequals  \| f_0 \|_{L^1}$ for $t \in [0, T]$.
			Note that the lifespan $T$ depends on $\eps, \phi$ and $\| f_0 \|_{L^1}$. Thanks to the conservation of mass, we can continue to  construct the solution on $[T, 2T], [2T, 3T], \cdots$ and get a global mild solution. That is, there is a unique measurable function $f^{\eps} \in \mathcal{A}_{\infty}$  satisfying:
			there is a null set $Z \subset \R^3$ s.t., for all $t \geq  0$ and $v \in \R^3 \setminus Z$, 
			\ben \label{mild-solution-existence}
			f^{\eps}(t,v) = f_{0}(v) + \int_{0}^{t} Q_{UU}^{\eps}(f^{\eps})(\tau, v) \mathrm{d}\tau, \quad 0 \leq f^{\eps}(t,v) \leq \eps^{-3}.
			\een
			Moreover, for all $t \geq  0$, $\|f^{\eps}(t)\|_{L^1} = \|f_0\|_{L^1}.$
		\end{proof}

		\begin{proof}[Proof of Theorem \ref{gwp-F-D}: (Propagation of regularity)] 
			We will  focus on the propagation of the regularity in weighted Sobolev spaces uniformly in $\eps$ for the mild solution. Suppose that
			\ben \label{def-T-max}
			T_M\colonequals\sup_{t>0}\bigg\{t\big| \sup_{0\le s\le t} \|f^{\eps}(t)\|_{H^N_{l}} \le 2 \|f_0\|_{H^N_{l}} \bigg\}. \een
			The main goal is to prove that $T_M = T_M(N,l,\phi, \|f_0\|_{H^N_{l}})$ is 
			strictly positive and independent of $\eps$. 
			
			Let $f_0 \in H^{N}_{l}$ with $N, l \geq 2$. Using \eqref{mild-solution-existence} and  Theorem \ref{UU-in-H-N-l-uniform-in-eps},  the conservation of mass and the upper bound: $f(t)\le \eps^{-3}$ will lead to  that
			\ben \label{hnl-norm-t1-2-t2}\nonumber
			&&\f{1}{2}\|f^{\eps}(t)\|_{H^{N}_{l}}^2 - \f{1}{2}\|f^{\eps}(0)\|_{H^{N}_{l}}^2 
			=  \sum_{|\alpha| \leq  N} \int_{0}^{t}  \lr{W_{l} \pa^{\alpha} Q_{UU}^{\eps}(f^{\eps}(\tau)), W_{l} \pa^{\alpha}f^{\eps} (\tau)}
			\mathrm{d}\tau\\&&\le C_{N,l,\phi}\int_{0}^{t}  \|f^{\eps}(\tau)\|_{H^{N}_{l}}^2(\|f^{\eps}(\tau)\|_{H^{N}_{l}} + \|f^{\eps}(\tau)\|_{H^{N}_{l}}^2)d\tau.
			\een 
			Let  
			$ T_{M}^{*}\colonequals \f{21}{ 4C_{N,l,\phi} ( \|f_0\|_{H^{N}_{l}}^2 +1)^{3/2}}.$ Then by
			Gr\"{o}nwall's inequality, we derive that for $t\in[0,T_M^*]$, 
			\ben \label{bounded-2-times-hnl}
			\sup_{0\le t \le T_{M}^{*}} \|f^{\eps}(t)\|_{H^N_{l}} \le 2 \|f_0\|_{H^N_{l}}, \mbox{which implies that}\,\, T_M\ge T_M^*.
			\een
			
			We now prove \eqref{continuous-in-time} based on \eqref{mild-solution-existence} and \eqref{bounded-2-times-hnl}. By \eqref{mild-solution-existence} and Minkowski's  inequality, we have 
			\beno
			\|f^{\eps}(t_2,v) - f^{\eps}(t_1,v)\|_{H^{N-2}_{l}} \lesssim \int_{t_1}^{t_2} \| Q_{UU}^{\eps}(f^{\eps})(\tau, v) \|_{H^{N-2}_{l}}  \mathrm{d}\tau.
			\eeno
			From \eqref{lose-order-2-derivative} and \eqref{bounded-2-times-hnl}, we are led to that
			\beno
			\|f^{\eps}(t_2,v) - f^{\eps}(t_1,v)\|_{H^{N-2}_{l}} \lesssim \int_{t_1}^{t_2}(\| f^{\eps}(\tau) \|_{H^{N}_{l}}^2
			+ \| f^{\eps}(\tau) \|_{H^{N}_{l}}^3)  \mathrm{d}\tau
			\lesssim (\| f_0 \|_{H^{N}_{l}}^2
			+ \| f_0 \|_{H^{N}_{l}}^3)  (t_2 - t_1),
			\eeno
			which gives \eqref{continuous-in-time}.
		\end{proof}

		\subsection{Bose-Einstein particles}
		In this subsection, we will prove the local well-posedness and propagation of regularity uniformly in $\eps$ for Bose-Einstein particles. The proof will follow the same spirit as we did in the previous subsection. However,  density of Bose-Einstein particles could blow up, unlike Fermi-Dirac particles whose density has a natural   bound $f\leq \eps^{-3}$. This motivates us to  construct the solution in $L^{1} \cap L^{\infty}$ space. 
		
		Recall \eqref{BUU} and \eqref{OPQUU}, we consider
		$\pa_t f
		=Q_{UU}^{\eps}(f)$,
		where $Q_{UU}^{\eps}$ denotes the Uehling-Uhlenbeck operator for Bose-Einstein particles, i.e., $Q_{UU}^{\eps}(f) = \int_{\R^3\times\SS^2}
		B^{\eps} \Phi^{\eps}(f) \mathrm{d}\sigma \mathrm{d}v_*$,
		where 
		\beno
		\Phi^{\eps}(f) \colonequals 
		f_*' f' (1 + \eps^3 f_* + \eps^3 f) - f_* f (1 + \eps^3 f_*' + \eps^3 f')
		.\eeno
		As usual,
		we define the gain and loss terms by
		\beno 
		Q^{\eps,+}(f) = \int_{\R^3\times\SS^2}
		B^{\eps} \Phi^{\eps,+}(f) \mathrm{d}\sigma \mathrm{d}v_*,
		\quad  Q^{\eps,-}(f) = \int_{\R^3\times\SS^2}
		B^{\eps} \Phi^{\eps,-}(f) \mathrm{d}\sigma \mathrm{d}v_*, \eeno
		where 
		\beno
		\Phi^{\eps,+}(f) \colonequals 
		f_*' f' (1 + \eps^3 f_* + \eps^3 f), \quad 
		\Phi^{\eps,-}(f) \colonequals 
		f_* f (1 + \eps^3 f_*' + \eps^3 f').
		\eeno

		We begin with a lemma on the upper bound of collision operator in $L^{1} \cap L^{\infty}$ space.
		\begin{lem} 
			Let $f, g \in L^1 \cap L^{\infty}, 0 \leq f, g$, then 
			\ben \label{Q-g-L1-B-E}
			\|Q_{UU}^{\eps}(f)\|_{L^1}      
			\lesssim \eps^{-3} I_{0} (1 + \|f\|_{L^\infty})   \|f\|_{L^1}^2,
			\\ \label{Q-g-Linfty-B-E}
			\|Q_{UU}^{\eps}(f)\|_{L^\infty}  \lesssim (I_{3} + \eps^{-3}  I_{0}) (1 + \|f\|_{L^\infty})  \|f\|_{L^\infty} (\|f\|_{L^1} + \|f\|_{L^\infty}).
			\een
			Moreover, 
			\ben
			\label{Q-g-minus-f-L1-and-infty-B-E}
			&& \|Q_{UU}^{\eps}(f) - Q_{UU}^{\eps}(g)\|_{L^{1} \cap L^\infty}
			\\ \nonumber &\lesssim&
			(\eps^{-3}  I_{0} + I_{3}) (\| f \|_{L^1 \cap L^\infty} + \| g \|_{L^1 \cap L^\infty} + 1) (\| f \|_{L^1 \cap L^\infty} + \| g \|_{L^1 \cap L^\infty}) \| f -g \|_{L^1 \cap L^\infty}
			.\een
		\end{lem} 
		\begin{proof} We give the estimates term by term.
			Similarly to \eqref{F-D-Q-L1-into-gain-and-loss} and \eqref{F-D-Q-L1-loss-up},  \eqref{Q-g-L1-B-E} follows the facts  
			$\|Q_{UU}^{\eps}(f)\le 2 \|Q^{\eps,-}(f)\|_{L^1}\|_{L^1}$ 
			and  
			$\|Q^{\eps,-}(f)\|_{L^1}\lesssim \eps^{-3} I_{0} (1 + \|f\|_{L^\infty})   \|f\|_{L^1}^2.$
			
			$\bullet$ To derive \eqref{Q-g-Linfty-B-E}, we observe that $Q^{\eps,-}(f)(v)\lesssim \eps^{-3} I_0 (1 + \|f\|_{L^\infty}) \|f\|_{L^\infty} \|f\|_{L^1}$ and 
			\ben \nonumber
			Q^{\eps,+}(f)(v)  
			&\lesssim& (1 + \|f\|_{L^\infty}) \|f\|_{L^\infty} \int
			B^{\eps}_1 
			f_*'  \mathrm{d}\sigma \mathrm{d}v_*  + (1 + \|f\|_{L^\infty}) \|f\|_{L^\infty}^2 \int
			B^{\eps}_3 
			\mathrm{d}\sigma \mathrm{d}v_*
			\\ \label{Q-plus-L-infty-B-E} &\lesssim& \eps^{-3} I_{0} (1 + \|f\|_{L^\infty}) \|f\|_{L^\infty} |f|_{L^1} + I_{3} (1 + \|f\|_{L^\infty}) \|f\|_{L^\infty}^2.
			\een
			Then \eqref{Q-g-Linfty-B-E} follows.
			
			$\bullet$ To show \eqref{Q-g-minus-f-L1-and-infty-B-E}, we notice that 
			$\|Q_{UU}^{\eps}(f) - Q_{UU}^{\eps}(g)\|_{L^1} \leq   2 \int B^{\eps} |\Phi^{\eps,-}(f) - \Phi^{\eps,-}(g)| \mathrm{d}V$ and   
			\ben \label{Phi-minus-f-minus-g}
			|\Phi^{\eps,-}(f) - \Phi^{\eps,-}(g)| \leq (1 + 2\|f\|_{L^{\infty}}) \left(
			|f_* - g_{*}| f + g_* |f - g| \right) + 2\|f - g\|_{L^{\infty}} g g_* 
			,\een for $0< \eps \leq 1$ and $ 0 \leq f, g \in L^{\infty}$.
			From these together with \eqref{upper-bound-of-A-eps}, we get that
			\beno
			\int
			B^{\eps} |\Phi^{\eps,-}(f) - \Phi^{\eps,-}(g)| 
			\mathrm{d}V
			\lesssim \eps^{-3} I_{0} (1 + \|f\|_{L^{\infty}}) (\| f \|_{L^1} + \| g \|_{L^1})\| f -g \|_{L^1}
			+   \eps^{-3} I_{0} \|f - g\|_{L^{\infty}}  \| g \|_{L^1}^2
			,\eeno
			which gives
			\ben
			\label{Q-g-minus-f-L1-B-E}
			\|Q_{UU}^{\eps}(f) - Q_{UU}^{\eps}(g)\|_{L^1}  \lesssim \eps^{-3} I_{0} (1 + \|f\|_{L^{\infty}} + \| g \|_{L^1}) (\| f \|_{L^1} + \| g \|_{L^1})\| f -g \|_{L^1 \cap L^{\infty}}.
			\een
			
			$\bullet$ To get the $L^\infty$ bounds,  \eqref{Phi-minus-f-minus-g} and  \eqref{upper-bound-of-A-eps} yield  that
			\beno
			&&|Q^{\eps,-}(f) - Q^{\eps,-}(g)|
			\leq 
			\int
			B^{\eps} |\Phi^{\eps,-}(f) - \Phi^{\eps,-}(g)| 
			\mathrm{d}v_* \mathrm{d}\sigma
			\\ &\lesssim& \eps^{-3} I_{0} (1 + \|f\|_{L^{\infty}}) (\| f \|_{L^\infty}\| f -g \|_{L^1} + \| f -g \|_{L^\infty}\| g \|_{L^1})
			+  \eps^{-3} I_{0} \|f - g\|_{L^{\infty}} \|g\|_{L^1} \|g\|_{L^{\infty}}.
			\eeno
			For the gain term, similar to \eqref{Phi-minus-f-minus-g}, by exchanging $(v, v_*)$ and $(v', v'_*)$, we have 
			\ben \label{Phi-plus-f-minus-g}\quad
			|\Phi^{\eps,+}(f) - \Phi^{\eps,+}(g)| &\leq& (1 + 2\|f\|_{L^{\infty}}) \left(
			|f_*' - g_{*}'| f' + g_*' |f' - g'| \right) + 2 \|f - g\|_{L^{\infty}} g' g_*' 
			\\ \nonumber 
			&\leq& (1 + 2\|f\|_{L^{\infty}}) \left(
			|f_*' - g_{*}'| \|f\|_{L^{\infty}} + g_*' \|f-g\|_{L^{\infty}} \right) + 2\|f - g\|_{L^{\infty}} \|g\|_{L^{\infty}} g_*' 
			.\een
			Now combining it   with \eqref{Phi-plus-f-minus-g}, \eqref{B1-sigma-vStar} and \eqref{bounded-l1-norm-simple},  we have 
			\ben \label{L-infty-estimate}
			&& \|Q_{UU}^{\eps}(f) - Q_{UU}^{\eps}(g)\|_{L^\infty}
			\\ \nonumber &\lesssim&
			(\eps^{-3}  I_{0} + I_{3}) (\| f \|_{L^1 \cap L^\infty} + \| g \|_{L^1 \cap L^\infty} + 1) (\| f \|_{L^1 \cap L^\infty} + \| g \|_{L^1 \cap L^\infty}) \| f -g \|_{L^1 \cap L^\infty}
			.\een We complete the proof of the lemma.
		\end{proof}
		
		We are now ready to prove Theorem \ref{lwp-B-E}.

		\begin{proof}[Proof of Theorem \ref{lwp-B-E}]  We introduce the function space  $\mathcal{E}_{T}  \colonequals L^{\infty}([0,T]; L^{1}(\R^3)\cap L^\infty(\R^3))$ associated with the norm $\|f\|_{ET} \colonequals \sup_{0 \leq t \leq T} \|f(t)\|_{L^{1}\cap L^\infty}$. 
			Similarly we define the operator $J^{\eps}(\cdot)$ on 
			$\mathcal{E}_{T}$ by  
			\beno J^{\eps}(f)(t,v) \colonequals f_{0}(v) + \int_{0}^{t} Q_{UU}^{\eps}(|f| )(\tau, v) \mathrm{d}\tau.\eeno

			$\bullet$ Given $f \in \mathcal{E}_{T}$, we now check that $J^{\eps}(f) \in \mathcal{E}_{T}$.  
			By \eqref{Q-g-L1-B-E} and \eqref{Q-g-Linfty-B-E}, we derive that 
			\ben \label{map-on-itself-ET}
			\|J^{\eps}(f)\|_{ET} \leq  \|f_{0}\|_{L^1 \cap L^{\infty}} + T C_{\eps, \phi} (1+ \|f\|_{ET})
			\|f\|_{ET}^2,
			\een
			where $C_{\eps, \phi} \lesssim I_{3} + \eps^{-3}  I_{0}$.
			This means  $J^{\eps}(\cdot)$ is an operator on  
			$\mathcal{E}_{T}$.

			$\bullet$  Let $\mathcal{F}_{T}  \colonequals \{f \in  \mathcal{F}_{T},   \|f\|_{ET} \leq 2  \|f_{0}\|_{L^1 \cap L^{\infty}} \}$. We want to show that $J^{\eps}(\cdot)$ is a contraction mapping on the complete metric space $(\mathcal{F}_{T}, \|\cdot - \cdot\|_{ET})$ for  a short time $T>0$. By \eqref{map-on-itself-ET}, 
			if $T$ satisfies
			\ben \label{small-condition-1-on-T-boson}
			4 T C_{\eps, \phi} (1 + 2\|f_{0}\|_{L^1 \cap L^{\infty}})
			\|f_{0}\|_{L^1 \cap L^{\infty}} \leq 1,
			\een
			then
			$J^{\eps}(\cdot)$ is an operator onto  
			$\mathcal{F}_{T}$. Now we  prove the contraction property. Since 
			$||x| - |y|| \leq |x-y|$, for 
			$f, g \in \mathcal{F}_{T}$, 
			we get
			\beno
			\|J^{\eps}(|f|) - J^{\eps}(|g|)\|_{ET} &\leq& T C_{\epsilon, \phi} (\| f \|_{ET} + \| g \|_{ET} + 1) (\| f \|_{ET} + \| g \|_{ET}) \| f -g \|_{L^1 \cap L^\infty}
			\\ &\leq& 4\| f_0 \|_{L^1 \cap L^\infty}   (4\| f_0 \|_{L^1 \cap L^\infty}  + 1) T C_{\epsilon, \phi} \| f -g \|_{ET}.
			\eeno
			If $T$ satisfies
			\ben \label{small-condition-2-on-T-BE}
			4\| f_0 \|_{L^1 \cap L^\infty}   (4\| f_0 \|_{L^1 \cap L^\infty}  + 1) T C_{\epsilon, \phi} \leq \f12,
			\een
			then
			$J^{\eps}(\cdot)$ is a contraction mapping on the complete metric space $(\mathcal{F}_{T}, \|\cdot - \cdot\|_{T})$. By the fixed point theorem,  
			there exists 
			a unique $f \in \mathcal{F}_{T}$, s.t.  $f^{\eps} = J^{\eps}(f^{\eps})$. After a modification on the $v$-null sets,  there is a null set $Z \subset \R^3$, for all $t \in  [0, T]$ and $v \in \R^3 \setminus Z$, 
			\beno  
			f^{\eps}(t,v) = J^{\eps}(f^{\eps}) (t,v) = f_{0}(v) + \int_{0}^{t} Q_{UU}^{\eps}(|f^{\eps}| )(\tau, v) \mathrm{d}\tau.
			\eeno
			Recalling \eqref{Q-g-Linfty-B-E},
			we can see that for any $v \in \R^3 \setminus Z$ and $t_1, t_2 \in [0, T]$,
			\beno
			|f^{\eps}(t_2,v) - f^{\eps}(t_1,v)| \leq C_{\eps, \phi} (\| f_0 \|_{L^1 \cap L^\infty} + 1)\| f_0 \|_{L^1 \cap L^\infty}^2 |t_2 - t_1|,
			\eeno
			for some $C_{\eps, \phi} \lesssim I_{3} + \eps^{-3}  I_{0}$. That is, $f(\cdot, v)$ is uniformly continuous w.r.t. $t$ on $[0, T]$
			for any $v \in \R^3 \setminus Z$. Following the argument in \eqref{non-negativity} we can get that
			$f^{\eps}(t,v) \geq 0$ on $[0, T]$, i.e.,
			\ben \label{mild-solution-BE}
			f^{\eps}(t,v) = f_{0}(v) + \int_{0}^{t} Q_{UU}^{\eps}(f^{\eps})(\tau, v) \mathrm{d}\tau, \quad 0 \leq f^\eps(t,v).
			\een
			Conservation of mass is a direct result of \eqref{Q-g-L1-B-E} and the
			Fubini's Theorem. Moreover that the lifespan $T$  
			depends on $\eps, \phi$ and $\| f_0 \|_{L^1 \cap L^{\infty}}$, i.e., $T=T(\eps, \phi, f_0)$.  In order to extend the solution to a positive time  independent of $\eps$, we will prove the propagation of regularity in weighted Sobolev space $H^{N}_{l}$.

			Assume that $f_0 \in H^{N}_{l}$ with $N, l \geq 2$. Again by \eqref{mild-solution-BE} and Theorem \ref{UU-in-H-N-l-uniform-in-eps}, \eqref{hnl-norm-t1-2-t2} still holds. This shows that the solution verifies the  {\it a priori} estimate:  for any $t\in[0,T_M^*=\f{21}{ 4C_{N,l,\phi} ( \|f_0\|_{H^{N}_{l}}^2 +1)^{3/2}}]$, \beno
			\sup_{0\le t \le T_{M}^{*}} \|f^{\eps}(t)\|_{H^N_{l}} \le 2 \|f_0\|_{H^N_{l}}, \eeno
			which implies the {\it a priori} upper bound: 
			\beno \sup_{0\le t \le T_{M}^{*}} \|f^{\eps}(t)\|_{L^\infty}\le 2C_S\|f_0\|_{H^N_{l}}. \eeno
			This enables to use continuity argument to extend the lifespan from $T=T(\eps, \phi, f_0)$ to $ T_{M}^{*}$ which is independent of $\eps$. We end the proof.
		\end{proof}

		\setcounter{equation}{0}
		
		\section{Asymptotic formula}
		This section is devoted to the proof of Theorem \ref{mainthm}. Let $f^\eps$ and  $f$ be the solutions to \eqref{BUU} and \eqref{landau} respectively with the initial datum $f_0$ 
		on $[0, T^*]$ where $T^*=T^*(N, l, \phi, \|f_0\|_{H^{N+3}_{l+5}})$ given in Theorem \ref{gwp-F-D} and \ref{lwp-B-E}. The solutions are uniformly bounded in $H^{N+3}_{l+5}$. More precisely, 
		\ben \label{upper-bound-3-order-higher}
		\sup_{t \leq T^*}	\|f^\eps(t)\|_{H^{N+3}_{l+5}}  \leq  2 \|f_0\|_{H^{N+3}_{l+5}}, \quad 	\sup_{t \leq T^*}	\|f(t)\|_{H^{N+3}_{l+5}}  \leq  2 \|f_0\|_{H^{N+3}_{l+5}}.
		\een
		

		Let $R^\eps \colonequals \eps^{-\vartheta}(f^\eps-f)$, then it is not difficult to see that  $R^\eps$ verifies the following equation:
		\ben \label{EqR} \pa_t R^\eps&=&Q_1(f^\eps, R^\eps)+Q_1(R^\eps, f)+\eps^{-\vartheta}(Q_1(f,f)-Q_L(f,f))+\eps^{-\vartheta}(Q_2+Q_3)(f^\eps,f^\eps)
		\\ \nonumber 
		&&+ \eps^{-\vartheta}R(f^\eps, f^\eps, f^\eps).\een
		We will prove $R^\eps$ is bounded in $H^{N}_{l}$ over $[0, T^*]$. To this end, 
		we first derive an estimate for the operator difference $Q_1-Q_L$. This is the key point for the asymptotic formula.
		\begin{lem} \label{Q1-and-QL}
			It holds that 
			\ben \label{operator-diff-l2l}
			\|Q_1(g,h)-Q_L(g,h)\|_{L^2_l}
			\lesssim_{l} \eps^{\vartheta} I_{3+\vartheta} \left( \|g\|_{H^3_{l+5}}\|h\|_{H^3}+\|g\|_{H^3_2}\|h\|_{H^3_{l+3}} \right). 
			\een
		\end{lem}
		\begin{proof}
			The estimate is proved in the spirit of \cite{des1}. The proof is divided into several steps.

			\underline{\it Step 1: Reformulation of $Q_1$.} 
			Firstly, given $v, v_* \in \R^3$,
			we introduce an orthonormal basis of $\R^3$:
			\beno  \big(\f{v-v_*}{|v-v_*|},~~ h^1_{v,v_*}, ~~h^2_{v,v_*}\big). \eeno
			Then
			$\sigma= \f{v-v_*}{|v-v_*|}\cos\theta+ (\cos\varphi h^1_{v,v_*}+\sin\varphi h^2_{v,v_*})\sin\theta$,
			which implies
			$ v'=v+\f12A$ and $ v_*'=v_*-\f12A$
			with $$A=-(v-v_*)(1-\cos\theta)+|v-v_*|(\cos\varphi h^1_{v,v_*}+\sin\varphi h^2_{v,v_*})\sin\theta.$$
			Now $Q_1$ can be rewritten as \beno
			Q_1(g,h)= \int_{\R^3}\int_0^\f{\pi}{2}\int_0^{2\pi} B^{\eps}_1\big[g(v_*-\f12A)h(v+\f12A)-g(v_*)h(v)\big]\sin\theta \mathrm{d}\varphi \mathrm{d}\theta \mathrm{d}v_*.\eeno
			
			By the Taylor expansion formula up to order 3:
			\beno g(v_*-\f12A)=g(v_*)-\f12A\cdot\na_{v_*} g(v_*)
			+\f18A\otimes A:\na^2g(v_*)+r_1(v, v_*, \sigma),
			\\ h(v+\f12A)=h(v)+\f12A\cdot\na_{v} h(v)+\f18A\otimes A:\na^2h(v)+r_2(v, v_*, \sigma),  \eeno
			where 
			\beno  |r_1(v, v_*, \sigma)|\lesssim |A|^3\int_0^1 |\na^3 g (\iota(v_*))|\mathrm{d}\iota, \quad 
			|r_2(v, v_*, \sigma)|\lesssim |A|^3\int_0^1 |\na^3 h (\kappa(v))|
			\mathrm{d}\kappa.
			\eeno
			Then we arrive at
			\ben\label{apQ1}  Q_1(g,h)&=& \int_{\R^3}\int_0^\f{\pi}{2} \int_0^{2\pi} \big[
			\f12A \cdot (\na_{v}-\na_{v_*})(g(v_*)h(v)) \notag\\&&\qquad
			+\f18A\otimes A:(\na_v-\na_{v_*})^2(g(v_*)h(v))+\mathcal{R}^1(v, v_*, \sigma)\big] B^{\eps}_1 \sin\theta \mathrm{d}\varphi \mathrm{d}\theta \mathrm{d}v_*,\een
			where 
			\beno \mathcal{R}^1(v, v_*, \sigma) &=& r_1(v, v_*, \sigma)\big(h(v)+\f12A \cdot \na h(v)+\f18A\otimes A:\na^2h(v)+r_2(v, v_*, \sigma) \big)
			\\ &&
			+\f18A\otimes A:\na^2g(v_*)\big(\f12A \cdot \na h(v)+\f18A\otimes A:\na^2h(v)+r_2(v, v_*, \sigma) \big)
			\\ &&-\f12A \cdot 
			\na g(v_*)\big(\f18A\otimes A:\na^2h(v)+r_2(v, v_*, \sigma) \big)+g(v_*)r_2(v, v_*, \sigma).
			\eeno

			\underline{\it Step 2: Reduction of $Q_1$.} 
			According to \eqref{apQ1}, if we define
			$T^\eps(v-v_*) \colonequals \int_0^\f{\pi}{2} \int_0^{2\pi} \f12A B^{\eps}_1 \sin\theta \mathrm{d}\varphi \mathrm{d}\theta$  and
			$ U^\eps(v-v_*)\colonequals \int_0^\f{\pi}{2} \int_0^{2\pi} \f18A\otimes A  B^{\eps}_1\sin\theta \mathrm{d}\theta \mathrm{d}\varphi$,
			then
			\ben\label{colli1}  Q_1(g,h) &=& \int_{\R^3}\Big[T^\eps(v-v_*)\cdot
			(\na_{v}-\na_{v_*})(g(v_*)h(v))  +U^\eps(v-v_*):(\na_v-\na_{v_*})^2(g(v_*)h(v))\Big]\mathrm{d}v_*\nonumber\\&&
			+ \int_{\R^3}\int_0^\f{\pi}{2}\int_0^{2\pi} \mathcal{R}^1(v, v_*, \sigma)B^{\eps}_1\sin\theta \mathrm{d}\varphi \mathrm{d}\theta \mathrm{d}v_*.\een
			
			\noindent\underline{Computation of $T^\eps$.} It is not difficult to compute that
			\ben \label{formula-of-T}
			T^\eps(v-v_*)=-8\pi I_{3} |v-v_*|^{-3}(v-v_*)+\mathcal{R}^2(v-v_*) ,
			\een
			with \ben \label{DEF-R2}
			\mathcal{R}^2(v-v_*) = 8\pi|v-v_*|^{-3}(v-v_*)\int_{\f{\sqrt{2}|v-v_*|}{2\eps}}^\infty \hphi^2(r)r^3 \mathrm{d}r
			.
			\een
			
			\noindent\underline{Computation of $U^\eps$.} We claim that
			\ben\label{colli2} U^\eps(v-v_*)= a(v-v_*)+\mathcal{R}^3(v-v_*) ,\een
			where $a$ is the matrix defined in \eqref{def-a},
			and
			\ben \label{DEF-R3} \mathcal{R}^3(z) &=& 2\pi|z|^{-1} \Pi(z)
			\left(-\int_{\f{\sqrt{2}|z|}{2\eps}}^\infty \hphi^2(r)r^3 \mathrm{d}r-
			\int_0^{\f{\sqrt{2}|z|}{2\eps}}\hphi^2(r)r^3 (\eps r)^2 |z|^{-2}\mathrm{d}r \right)
			\\  \notag &&
			+ \f{\pi}{4} z \otimes z  \int_0^{\pi/2}  (1-\cos\theta)^2 B^{\eps}_1\sin\theta \mathrm{d}\theta.
			\een
			where $\Pi$ is the matrix defined in \eqref{def-a}.
			To see this, by definition, we have
			\beno  U^\eps(v-v_*) &=& \f18\int_0^\f{\pi}{2} \int_0^{2\pi} \big(-(v-v_*)(1-\cos\theta)+|v-v_*|(\cos\varphi h^1_{v,v_*}+\sin\varphi h^2_{v,v_*})\sin\theta\big)\\&&\otimes \big(-(v-v_*)(1-\cos\theta)+|v-v_*|(\cos\varphi h^1_{v,v_*}+\sin\varphi h^2_{v,v_*})\sin\theta\big)  B^{\eps}_1\sin\theta \mathrm{d}\theta \mathrm{d}\varphi \\
			&=&\f{\pi}4(v-v_*)\otimes(v-v_*) \int_0^{\pi/2}  (1-\cos\theta)^2 B^{\eps}_1\sin\theta \mathrm{d}\theta
			\\&& +\f18|v-v_*|^2\int_0^\f{\pi}{2} \int_0^{2\pi} \sin^2\theta(\cos^2\varphi  h^1_{v,v_*}\otimes h^1_{v,v_*}+ \sin^2\varphi h^2_{v,v_*}\otimes h^2_{v,v_*})B^{\eps}_1\sin\theta \mathrm{d}\theta \mathrm{d}\varphi.
			\eeno
			Use the fact that $\f{v-v_*}{|v-v_*|}\otimes \f{v-v_*}{|v-v_*|}+h^1_{v,v_*}\otimes h^1_{v,v_*}+h^2_{v,v_*}\otimes h^2_{v,v_*}=Id$, then
			\beno U^\eps(v-v_*) &=& 2\pi|v-v_*|^{-1}(Id- \f{v-v_*}{|v-v_*|}\otimes \f{v-v_*}{|v-v_*|}) \int_0^{\f{\sqrt{2}|v-v_*|}{2\eps}}\hphi^2(r)r^3(1-(\eps r)^2|v-v_*|^{-2})\mathrm{d}r
			\\&&+ \f{\pi}4(v-v_*)\otimes(v-v_*) \int_0^{\pi/2}  (1-\cos\theta)^2 B^{\eps}_1\sin\theta \mathrm{d}\theta, \eeno
			which is enough to get  \eqref{colli2}.
			
			Since $\na\cdot \big((|x|^2 Id- x\otimes x)f(|x|^2)\big)=-2xf(|x|^2)$, recalling \eqref{def-a}
			and  \eqref{formula-of-T}, one has
			\beno && (\na_v-\na_{v_*})\cdot a(v,v_*)=-8\pi I_{3} |v-v_*|^{-3}(v-v_*)  
			= T^\eps(v,v_*)-\mathcal{R}^2(v-v_*) . \eeno 
			Plugging this and \eqref{colli2} into 
			\eqref{colli1}, we get
			\beno Q_1(g,h)&=& \int_{\R^3} (\na_v-\na_{v_*})\cdot \big[a(v,v_*)(\na_v-\na_{v_*})(g_*h)\big]\mathrm{d}v_*+\int_{\R^3} \mathcal{R}^3(v-v_*) :\big[(\na_v-\na_{v_*})^2  (g_*h)\big]\mathrm{d}v_*\\
			&&+\int_{\R^3}\mathcal{R}^2(v-v_*)  \cdot (\na_v-\na_{v_*})(g_*h)  \mathrm{d}v_*  + \int_{\R^3}\int_0^\f{\pi}{2}\int_0^{2\pi} \mathcal{R}^1(v, v_*, \sigma)B^{\eps}_1\sin\theta \mathrm{d}\theta \mathrm{d}\varphi \mathrm{d}v_*.\eeno
			Note that the first integral term is the Landau operator.
			
			\underline{\it Step 3: Estimate of $Q_1-Q_L$.}  
			From the above equality,   we arrive at
			\beno 
			Q_1(g,h)-Q_L(g,h)  &=& \int \mathcal{R}^1(v, v_*, \sigma)B^{\eps}_1\sin\theta \mathrm{d}\theta \mathrm{d}\varphi \mathrm{d}v_*
			+\int \mathcal{R}^2(v-v_*)  \cdot (\na_v-\na_{v_*})(g_*h) \mathrm{d}v_*\\
			&& +  \int  \mathcal{R}^3(v-v_*) :\big[(\na_v-\na_{v_*})^2  (g_*h)\big]\mathrm{d}v_* \colonequals\sum_{i=1}^3 \mathcal{Q}_i.
			\eeno
			
			\underline{Estimate of $\mathcal{Q}_3$.} Recalling \eqref{DEF-R3} for $\mathcal{R}^3$,  it is easy to see
			\beno |\mathcal{R}^3(z) | 
			\lesssim \eps^{\vartheta}	I_{3+\vartheta} |z|^{-1-\vartheta},
			\eeno
			which gives
			\beno  \|\mathcal{Q}_3\|_{L^2_l} \lesssim  \eps^{\vartheta}	I_{3+\vartheta} \|g\|_{H^2_2}\|h\|_{H^2_l}.\eeno
			
			\underline{Estimate of $\mathcal{Q}_2$.} Recalling \eqref{DEF-R2},
			\beno
			\eps^{-\vartheta}\mathcal{R}^2(x)=
			8\pi |x|^{-3-\vartheta}x\int_{\f{\sqrt{2}|x|}{2\eps}}^\infty\hphi^2(r)r^3\big(\f{|x|}{\eps}\big)^{\vartheta}\mathrm{d}r,\eeno
			it is obvious that $\eps^{-\vartheta}|\mathcal{R}^2(x)|\lesssim I_{3+\vartheta} |x|^{-2-\vartheta}$. 
			If $\vartheta=1$, we claim that $\eps^{-1}\mathcal{R}^2$ is the kernel of a Calderon-Zygmund operator.
			To see this, we first compute directly that  for any $0<R_1<R_2<\infty$,
			\beno 
			\int_{R_1<|x|<R_2} \eps^{-1}\mathcal{R}^2(x) dx=0, \quad \sup_{R>0} \int_{R <|x|< 2R} \eps^{-1} |\mathcal{R}^2(x)| dx \lesssim I_{4}.
			\eeno
			%
			%
			Next we need to check that $\eps^{-1} \mathcal{R}^2$ satisfies H{\"o}rmander's condition:
			\ben
			\label{R2condition} \int_{|x|\ge 2|y|} |\eps^{-1}\mathcal{R}^2(x)-\eps^{-1}\mathcal{R}^2(x-y)|dx\lesssim I_{4}.
			\een 
			Note that 
			\beno
			\eps^{-1}\mathcal{R}^2(x)-\eps^{-1}\mathcal{R}^2(x-y)
			&=& 8\pi |x|^{-3}x\int_{\f{\sqrt{2}|x|}{2\eps}}^\infty\hphi^2(r)r^3 \eps^{-1} \mathrm{d}r
			- 8\pi |x-y|^{-3}(x-y)\int_{\f{\sqrt{2}|x-y|}{2\eps}}^\infty\hphi^2(r)r^3 \eps^{-1} \mathrm{d}r
			\\ &=& 8\pi \left(|x|^{-3}x - |x-y|^{-3}(x-y) \right)
			\int_{\f{\sqrt{2}|x|}{2\eps}}^\infty\hphi^2(r)r^3 \eps^{-1} \mathrm{d}r
			\\ &&+ 8\pi |x-y|^{-3}(x-y) \left(
			\int_{\f{\sqrt{2}|x|}{2\eps}}^\infty\hphi^2(r)r^3 \eps^{-1} \mathrm{d}r - \int_{\f{\sqrt{2}|x-y|}{2\eps}}^\infty\hphi^2(r)r^3 \eps^{-1} \mathrm{d}r \right)
			.		\eeno 
			Under the condition $|x|\ge 2|y|$, it is easy to see that $|x-y| \sim |x|$ and $\big|x|x|^{-3}-(x-y)|x-y|^{-3}\big|\lesssim |x|^{-3}|y|$, then 
			\beno
			|\eps^{-1}\mathcal{R}^2(x)-\eps^{-1}\mathcal{R}^2(x-y)|
			\lesssim  |x|^{-3}|y| 
			\int_{\f{\sqrt{2}|x|}{2\eps}}^\infty \hphi^2(r)r^3 \eps^{-1} \mathrm{d}r
			+ |x|^{-2} \int_{ \f{\sqrt{2}}{2\eps}
				\min\{|x|, |x-y|\}}^{ \f{\sqrt{2}}{2\eps} \max\{ |x|, |x-y|\}} \hphi^2(r)r^3 \eps^{-1} \mathrm{d}r
			.		\eeno 
			For the first term,
			\beno \int_{|x|\ge 2|y|} 
			|x|^{-3}|y| 
			\int_{\f{\sqrt{2}|x|}{2\eps}}^\infty\hphi^2(r)r^3 \eps^{-1} \mathrm{d}r \mathrm{d}x 
			\lesssim I_{4} |y| \int_{|x|\ge 2|y|} |x|^{-4} \mathrm{d}x \lesssim I_{4}. \eeno
			For the second term, 
			\beno  &&\int_{|x|\ge 2|y|} |x|^{-2} 
			\int_{ \f{\sqrt{2}}{2\eps}
				\min\{|x|, |x-y|\}}^{ \f{\sqrt{2}}{2\eps} \max\{ |x|, |x-y|\}} \hphi^2(r)r^3 \eps^{-1} \mathrm{d}r \mathrm{d}x \\
			&\lesssim&
			\int_{|x|\ge 2|y|} |x|^{-3}\bigg(\int_{\f{\sqrt{2}(|x|-|y|) }{2\eps}}^{\f{\sqrt{2}|x| }{2\eps}}+\int_{\f{\sqrt{2}|x| }{2\eps}}^{\f{\sqrt{2}(|x|+|y|) }{2\eps}}\bigg)  \hphi^2(r)r^4  \mathrm{d}r \mathrm{d}x
			\\
			&\lesssim&  \int_{|y|\le \sqrt{2}\eps r}\hphi^2(r)
			r^4 \mathrm{d}r\int_{\sqrt{2}\eps r\le |x|\le2\sqrt{2}\eps r}|x|^{-3} \mathrm{d}x+ \int_{2|y|\le \sqrt{2}\eps r}\hphi^2(r)r^4 \mathrm{d}r\int_{\sqrt{2}/2\eps r\le |x|\le\sqrt{2}\eps r}|x|^{-3} \mathrm{d}x\lesssim I_{4}.\eeno 
			Combining these two estimates will yield \eqref{R2condition}. 
			Thus we have
			\beno  \|\mathcal{Q}_2\|_{L^2_l} \lesssim  \eps^{\vartheta} I_{3+\vartheta} \|g\|_{H^2_2}\|h\|_{H^2_l}.\eeno
			
			\underline{Estimate of  $\mathcal{Q}_3$.}
			Recalling the definition of $\mathcal{R}^1(v, v_*, \sigma)$,  it is bounded by
			\beno 
			&&  |A|^3\bigg(|g_*|\int_0^1
			|\na^3h (\kappa(v))|\mathrm{d}\kappa+|(\na g)_*||\na^2 h|+|(\na^2 g)_*||\na h| +\int_0^1
			|\na^3g (\iota(v_*))| \mathrm{d}\iota |h|\bigg)
			\\&+& |A|^4 \bigg(|(\na g)_*|\int_0^1|\na^3h (\kappa(v))|\mathrm{d}\kappa+|(\na^2 g)_*||\na^2 h|+
			\int_0^1|\na^3g (\iota(v_*))| \mathrm{d}\iota |\na h|\bigg)
			\\ &+& |A|^5 \bigg(|(\na^2 g)_*|\int_0^1|\na^3h (\kappa(v))|\mathrm{d}\kappa +
			\int_0^1|\na^3g (\iota(v_*))| \mathrm{d}\iota |\na^2 h|\bigg)
			\\ &+& |A|^6  \int_0^1 |\na^3g (\iota(v_*))|\mathrm{d}\iota \int_0^1|\na^3h (\kappa(v))|\mathrm{d}\kappa \colonequals \sum_{i=1}^4 \mathcal{R}^1_i
			.\eeno
			Let  $\mathcal{Q}_3^i\colonequals \int_{v_*,\theta,\varphi} \mathcal{R}^1_i(v,v_*) B^{\eps}_1\sin\theta \mathrm{d}\theta \mathrm{d}\varphi \mathrm{d}v_*$. 
			Using the facts $|A|\lesssim  \sin(\theta/2)|v-v_*|$
			and $W_l \lesssim_{l} W_l(\kappa(v))+W_l(\iota(v_*))$, 
			by the C-S inequality and \eqref{Q1-result-with-factor-a-minus3-with-small-facfor},
			we have
			\beno  \sum_{i=1}^4|\lr{\mathcal{Q}_3^iW_l, F}| \lesssim_{l} \eps^{\vartheta}I_{3+\vartheta} \|F\|_{L^2} (\|g\|_{H^3_{l+5}}\|h\|_{H^3}+\|g\|_{H^3_2}\|h\|_{H^3_{l+3}}), \eeno
			which yields that $\|\mathcal{Q}_3\|_{L^2_l} \lesssim_{l}  \eps^{\vartheta} I_{3+\vartheta} (\|g\|_{H^3_{l+5}}\|h\|_{H^3}+\|g\|_{H^3_2}\|h\|_{H^3_{l+3}}).$
			
			The desired result \eqref{operator-diff-l2l} follows by patching together all the estimates.
		\end{proof}
		
		Now we are in a position to prove  the asymptotic expansion in Theorem \ref{mainthm}.
		
		\begin{proof}[Proof of Theorem \ref{mainthm}]

			To derive the asymptotic formula in the theorem, the key point is to give the energy estimates for $R^\eps$ in  the  space $H^{N}_{l}$ for $N \geq 0, l \geq 2$.  Recalling \eqref{EqR}, the proof is divided into several steps.
			\smallskip
			
			\underline{\it Step 1: Estimate of $Q_1(f^\eps, R^\eps)$.}  
			We claim that
			\beno 
			\sum_{m=0}^{N} \sum_{|\alpha|=m}  \lr{ \pa^{\alpha}Q_1(f^\eps, R^\eps) W_{l},  W_{l}\pa^{\alpha}R^\eps}
			\lesssim_{N, l, \phi} 
			\|f^\eps\|_{H^{N+2}_{l}} \|R^\eps\|_{H^{N}_{l}}^2.
			\eeno
			We need to consider $\lr{Q_1( \pa^{\alpha_1}f^\eps, \pa^{\alpha_2} R^\eps) W_{l},  W_{l}\pa^{\alpha}R^\eps}$ for $\alpha_1 + \alpha_2 =\alpha$. If $|\alpha_1| = 0$, we use \eqref{H-2-result} and Remark \ref{Q1-still-holds}. 
			If $|\alpha_1| = 1$, we use Lemma \ref{Q1EN-g-order-1}. If $|\alpha_1| \geq 2$, 
			we use Lemma \ref{Q1EN}.

			\medskip
			
			\underline{\it Step 2: Estimate of $Q_1(R^\eps, f)$.}
			Using \eqref{Q1-h1-h2-l2-weighted}, we have
			\beno 
			\sum_{m=0}^{N} \sum_{|\alpha|=m}  \lr{ \pa^{\alpha}Q_1 (R^\eps, f)  W_{l},  W_{l}\pa^{\alpha}R^\eps}
			\lesssim_{N, l, \phi} 
			\|f \|_{H^{N+2}_{l}} \|R^\eps\|_{H^{N}_{l}}^2.
			\eeno
			
			\medskip

			\underline{\it Step 3: Estimate of $\eps^{-\vartheta}(Q_1(f,f)-Q_L(f,f))$.}
			Note that \eqref{operator-diff-l2l}
			yields 
			\beno \eps^{-\vartheta}\|Q_1(g,h)-Q_L(g,h)\|_{H^m_l}
			\lesssim_{m, l} I_{3+\vartheta} \left( \|g\|_{H^{m+3}_{l+5}}\|h\|_{H^{m+3}}+\|g\|_{H^{m+3}_2}
			\|h\|_{H^{m+3}_{l+3}} \right),\eeno
			which gives
			\beno 
			\sum_{m=0}^{N}\sum_{|\alpha|=m, \alpha_1+\alpha_2=\alpha}  \big| \lr{\eps^{-\vartheta}(Q_1(\pa^{\alpha_1}f,\pa^{\alpha_2}f)-Q_L(\pa^{\alpha_1}f,\pa^{\alpha_2}f)) W_{l}, W_{l}\pa^{\alpha}R^\eps}\big| \lesssim_{N, l, \phi} \|f \|_{H^{N+3}_{l+5}}^2 \|R^\eps\|_{H^{N}_{l}}.
			\eeno

			\medskip

			\underline{\it Step 4: Estimate of $\eps^{-\vartheta}(Q_2+Q_3)(f^\eps,f^\eps) $. }
			We claim
			that
			\beno 
			\sum_{m=0}^{N}\sum_{|\alpha|=m, \alpha_1+\alpha_2=\alpha} |\lr{\eps^{-\vartheta}(Q_2+Q_3)(\pa^{\alpha_1} f^\eps,\pa^{\alpha_2} f^\eps)) W_{l},  W_{l}\pa^{\alpha}R^\eps}|\lesssim_{N, l, \phi}   \|f^\eps \|_{H^{N+2}_{l}}^2 \|R^\eps\|_{H^{N}_{l}}.
			\eeno
			For the $Q_2$ term, 
			by using \eqref{g-h-total-vartheta}, we get
			\beno
			|\lr{\eps^{-\vartheta} Q_2(\pa^{\alpha_1} f^\eps,\pa^{\alpha_2} f^\eps)) W_{l},  W_{l}\pa^{\alpha}R^\eps}|
			| \lesssim_{l}  (I_{3+\vartheta} + I^{\prime}_{3+\vartheta}) \|f^\eps \|_{H^{N+2}_{l}}^2 \|R^\eps\|_{H^{N}_{l}}.
			\eeno 
			For the $Q_3$ term, by using \eqref{Q-3-vartheta}, we get 
			\beno
			|\lr{\eps^{-\vartheta} Q_3(\pa^{\alpha_1} f^\eps,\pa^{\alpha_2} f^\eps)) W_{l},  W_{l}\pa^{\alpha}R^\eps}|
			| \lesssim_{l} I_{3} \|f^\eps \|_{H^{N+2}_{l}}^2 \|R^\eps\|_{H^{N}_{l}}.
			\eeno
			
			\medskip

			\underline{\it Step 5: Estimate of $\eps^{-\vartheta}R(f^\eps, f^\eps, f^\eps)
				$.}  Applying \eqref{R-general-rough-versionh-h2-h2-h2-weighted}, we have
			\beno 
			\sum_{m=0}^{N} \sum_{|\alpha|=m} |\lr{\eps^{-3} \pa^{\alpha}R(f^\eps, f^\eps, f^\eps) W_{l},  W_{l}\pa^{\alpha}R^\eps} | \lesssim_{N, l, \phi}  \|f^\eps \|_{H^{N+2}_{l}}^3 \|R^\eps\|_{H^{N}_{l}}.
			\eeno
			
			\medskip
			
			\underline{\it Step 6: Closure of the energy estimates.}
			Patching together all the estimates in the previous steps, we  arrive at
			\beno \f{d}{\mathrm{d}t} \|R^\eps\|_{H^{N}_{l}}^2 &\lesssim_{N, l, \phi} &	
			\|f^\eps\|_{H^{N}_{l}} \|R^\eps\|_{H^{N}_{l}}^2 + \|f \|_{H^{N+2}_{l}} \|R^\eps\|_{H^{N}_{l}}^2 +  \|f \|_{H^{N+3}_{l+5}}^2 \|R^\eps\|_{H^{N}_{l}} 
			\\ &&+   \|f^\eps \|_{H^{N+2}_{l}}^2 \|R^\eps\|_{H^{N}_{l}} 
			+  \|f^\eps \|_{H^{N+2}_{l}}^3 \|R^\eps\|_{H^{N}_{l}}.
			\eeno
			Using the uniform upper bounds \eqref{upper-bound-3-order-higher}
			of $\|f^{\eps} \|_{H^{N+3}_{l+5}}$ and $\|f \|_{H^{N+3}_{l+5}}$, we get
			\beno \f{d}{\mathrm{d}t} \|R^\eps\|_{H^{N}_{l}} \lesssim_{N, l, \phi} (\|f_0\|_{H^{N+3}_{l+5}} + \|f_0\|_{H^{N+3}_{l+5}}^3)
			(\|R^\eps\|_{H^{N}_{l}} + 1).
			\eeno
			Then \eqref{error-estimate} follows the Gr\"{o}nwall's inequality. 
		\end{proof}
		\setcounter{equation}{0}

		{\bf Acknowledgments.} The work was initiated when M. Pulvirenti visited Tsinghua
		University in 2016. The work was partially supported by National Key Research and Development Program of China under the grant 2021YFA1002100. 
		Ling-Bing He was also supported by NSF of China  
		under the grant 12141102.
		Yu-Long Zhou was  partially supported by NSF of China under the grant 12001552,  Science and Technology Projects in Guangzhou under the grant 202201011144, and
		Youth Talent Support Program of Guangdong Provincial Association for Science and Technology under the grant SKXRC202311,

	\end{document}